\documentstyle{amltd2004}
\begin{document}
\annalsline{159}{2004}
\received{April 17, 2001}
\revised{April 4, 2002}
\startingpage{641}
\def\bye{\end{document}}
 \font\tenrm=cmr10
\def\ritem#1{\item[{\rm #1}]}

\input boxedeps.tex 
\SetepsfEPSFSpecial 
\HideDisplacementBoxes
\def\figin#1#2{
$$
 {\BoxedEPSF{#1.eps scaled
#2}%
}%
$$
\noindent}
\def\nhs{\hskip-8pt}
\def\hiha{\vglue8pt \noindent \hglue10pt \hangindent=48pt\hangafter=1}
\input amssym.def
\input amssym.tex

\catcode`\@=11
\font\twelvemsb=msbm10 scaled 1100
\font\tenmsb=msbm10
\font\ninemsb=msbm10 scaled 800
\newfam\msbfam
\textfont\msbfam=\twelvemsb  \scriptfont\msbfam=\ninemsb
  \scriptscriptfont\msbfam=\ninemsb
\def\msb@{\hexnumber@\msbfam}
\def\Bbb{\relax\ifmmode\let\next\Bbb@\else
 \def\next{\errmessage{Use \string\Bbb\space only in math
mode}}\fi\next}
\def\Bbb@#1{{\Bbb@@{#1}}}
\def\Bbb@@#1{\fam\msbfam#1}
\catcode`\@=12

 \catcode`\@=11
\font\twelveeuf=eufm10 scaled 1100
\font\teneuf=eufm10
\font\nineeuf=eufm7 scaled 1100
\newfam\euffam
\textfont\euffam=\twelveeuf  \scriptfont\euffam=\teneuf
  \scriptscriptfont\euffam=\nineeuf
\def\euf@{\hexnumber@\euffam}
\def\frak{\relax\ifmmode\let\next\frak@\else
 \def\next{\errmessage{Use \string\frak\space only in math
mode}}\fi\next}
\def\frak@#1{{\frak@@{#1}}}
\def\frak@@#1{\fam\euffam#1}
\catcode`\@=12


    \def\anf{$\lower1.2ex\hbox{"}$}

\def\frac#1#2{{#1 \over #2}}

\def\>{>\!\!>}

\def\<{<\!\!<}

\def\into{\hookrightarrow}

\def\ssarr{\hbox to 30pt{\rightarrowfill}}

\def\sarr{\hbox to 40pt{\rightarrowfill}}

\def\arr{\hbox to 60pt{\rightarrowfill}}

\def\larr{\hbox to 60pt{\leftarrowfill}}

\def\Arr{\hbox to 80pt{\rightarrowfill}}

\def\mapdown#1{\Big\downarrow\rlap{$\vcenter{\hbox{$\scriptstyle#1$}}$}}

\def\lmapdown#1{\llap{$\vcenter{\hbox{$\scriptstyle#1$}}$}\Big\downarrow}

\def\mapright#1{\smash{\mathop{\arr}\limits^{#1}}}

\def\mapleft#1{\smash{\mathop{\larr}\limits^{#1}}}

\def\ssmapright#1{\smash{\mathop{\ssarr}\limits_{#1}}}

\def\smapright#1{\smash{\mathop{\sarr}\limits_{#1}}}

{}\def\lmapright#1{\smash{\mathop{\arr}\limits_{#1}}}

\def\Mapright#1{\smash{\mathop{\Arr}\limits^{#1}}}

\def\lMapright#1{\smash{\mathop{\Arr}\limits_{#1}}}

\def\mapup#1{\Big\uparrow\rlap{$\vcenter{\hbox{$\scriptstyle#1$}}$}}

\def\lmapup#1{\llap{$\vcenter{\hbox{$\scriptstyle#1$}}$}\Big\uparrow}

\def\vline{\hbox{\Bigg\vert}}

\def\ad{\mathop{\rm ad}\nolimits}

\def\add{\mathop{\rm add}\nolimits}

\def\arsinh{\mathop{\rm arsinh}\nolimits}

\def\arcosh{\mathop{\rm arcosh}\nolimits}

\def\artanh{\mathop{\rm artanh}\nolimits}

\def\arctanh{\mathop{\rm arctanh}\nolimits}

\def\Ad{\mathop{\rm Ad}\nolimits}

\def\Aff{\mathop{\rm Aff}\nolimits}

\def\algint{\mathop{\rm algint}\nolimits}

\def\Aut{\mathop{\rm Aut}\nolimits}

\def\Bil{\mathop{\rm Bil}\nolimits}

\def\ch{\mathop{\rm char}\nolimits}

\def\card{\mathop{\rm card}\nolimits}

\def\cl{\mathop{\rm cl}\nolimits}

\def\codim{\mathop{\rm codim}\nolimits}

\def\Comp{\mathop{\rm Comp}\nolimits}

\def\comp{\mathop{\rm comp}\nolimits}

\def\compr{\mathop{\rm compr}\nolimits}

\def\cone{\mathop{\rm cone}\nolimits}

\def\conv{\mathop{\rm conv}\nolimits}

\def\Der{\mathop{\rm Der}\nolimits}

\def\deg{\mathop{\rm deg}\nolimits}

\def\der{\mathop{\rm der}\nolimits}

\def\det{\mathop{\rm det}\nolimits}

\def\diag{\mathop{\rm diag}\nolimits}

\def\Diff{\mathop{\rm Diff}\nolimits}

\def\dist{\mathop{\rm dist}\nolimits}

\def\Exp{\mathop{\rm Exp}\nolimits}

\def\ev{\mathop{\rm ev}\nolimits}

\def\epi{\mathop{\rm epi}\nolimits}

\def\Fix{\mathop{\rm Fix}\nolimits}

\def\End{\mathop{\rm End}\nolimits}

\def\Ext{\mathop{\rm Ext}\nolimits}

\def\grad{\mathop{\rm grad}\nolimits}

\def\Gal{\mathop{\rm Gal}\nolimits}

\def\Gl{\mathop{\rm Gl}\nolimits}

\def\GL{\mathop{\rm GL}\nolimits}

\def\hats #1{\hat{\hat{\hbox{$#1$}}}}

\def\Herm{\mathop{\rm Herm}\nolimits}

\def\HSp{\mathop{\rm HSp}\nolimits}

\def\Hol{\mathop{\rm Hol}\nolimits}%

\def\Hom{\mathop{\rm Hom}\nolimits}%

\def\id{\mathop{\rm id}\nolimits} 

\def\im{\mathop{\rm im}\nolimits}

\def\Im{\mathop{\rm Im}\nolimits}

\def\inf{\mathop{\rm inf}\nolimits}

\def\Ind{\mathop{\rm Ind}\nolimits}

\def\Inn{\mathop{\rm Inn}\nolimits}

\def\Int{\mathop{\rm int}\nolimits}

\def\Iso{\mathop{\rm Iso}\nolimits}

\def\Max{\mathop{\rm Max}\nolimits}

\def\Mot{\mathop{\rm Mot}\nolimits}

\def\Mp{\mathop{\rm Mp}\nolimits}

\def\OO{\mathop{\rm O}\nolimits}

\def\op{\mathop{\rm op}\nolimits}

\def\Out{\mathop{\rm Out}\nolimits}

\def\Pol{\mathop{\rm Pol}\nolimits}

\def\Prim{\mathop{\rm Prim}\nolimits}

\def\PGl{\mathop{\rm PGl}\nolimits}

\def\PSl{\mathop{\rm PSl}\nolimits}

\def\PU{\mathop{\rm PU}\nolimits}

\def\rad{\mathop{\rm rad}\nolimits}

\def\rank{\mathop{\rm rank}\nolimits}

\def\reg{\mathop{\rm reg}\nolimits}

\def\resi{\mathop{\rm res}\nolimits}

\def\Re{\mathop{\rm Re}\nolimits}

\def\rk{\mathop{\rm rk}\nolimits}

\def\sgn{\mathop{\rm sgn}\nolimits}

\def\Sl{\mathop{\rm Sl}\nolimits}

\def\SO{\mathop{\rm SO}\nolimits}

\def\span{\mathop{\rm span}\nolimits}

\def\Skew{\mathop{\rm Skew}\nolimits}

\def\Symm{\mathop{\rm Symm}\nolimits}

\def\Sp{\mathop{\rm Sp}\nolimits}

\def\Spec{\mathop{\rm Spec}\nolimits}

\def\Spin{\mathop{\rm Spin}\nolimits}

\def\SS{\mathop{\rm S}\nolimits}

\def\St{\mathop{\rm St}\nolimits}

\def\SU{\mathop{\rm SU}\nolimits}

\def\sup{\mathop{\rm sup}\nolimits}

\def\supp{\mathop{\rm supp}\nolimits}

\def\tr{\mathop{\rm tr}\nolimits}

\def\UU{\mathop{\rm U}\nolimits}

\def\Up{\mathop{\rm Up}\nolimits}

\def\vol{\mathop{\rm vol}\nolimits}

\def\0{{\bf 0}}

\def\1{{\bf 1}}

\def\a{{\frak a}}

\def\aff{{\frak {aff}}}

\def\aut{{\frak {aut}}}

\def\b{{\frak b}}

\def\cc{{\frak c}}

\def\d{{\frak d}}

\def\e{{\frak e}}

\def\f{{\frak f}}

\def\g{{\frak g}}

\def\gl{{\frak {gl}}}

\def\h{{\frak h}}

\def\ii{{\frak i}}

\def\j{{\frak j}}

\def\k{{\frak k}}

\def\l{{\frak l}}

\def\m{{\frak m}}

\def\mot{{\frak {mot}}}

\def\n{{\frak n}}

\def\oo{{\frak o}}

\def\p{{\frak p}}

\def\q{{\frak q}}

\def\r{{\frak r}}

\def\s{{\frak s}}

\def\hsp{{\frak {hsp}}}

\def\osc{{\frak {osc}}}

\def\sp{{\frak {sp}}}

\def\st{{\frak {st}}}

\def\su{{\frak {su}}}

\def\so{{\frak {so}}}

\def\sL{{\frak {sl}}}

\def\t{{\frak t}}

\def\uu{{\frak u}}

\def\uP{{\frak {up}}}

\def\Vir{{\frak {vir}}}

\def\vv{{\frak v}}

\def\w{{\frak w}}

\def\x{{\frak x}}

\def\y{{\frak y}}

\def\z{{\frak z}}

\def\L{\mathop{\bf L}\nolimits}

\def\B{{\Bbb B}} 

\def\C{{\Bbb C}} 

\def\D{{\Bbb D}} 

\def\F{{\Bbb F}} 

\def\H{{\Bbb H}} 

\def\K{{\Bbb K}} 

\def\LL{{\Bbb L}} 

\def\M{{\Bbb M}} 

\def\N{{\Bbb N}} 

\def\OOO{{\Bbb O}} 

\def\P{{\Bbb P}} 

\def\Q{{\Bbb Q}} 

\def\R{{\Bbb R}} 

\def\T{{\Bbb T}} 

\def\Z{{\Bbb Z}}

\def\:{\colon}  

\def\.{{\cdot}}

\def\|{\Vert}

\def\bsk{\bigskip}

\def \de {\, {\buildrel \rm def \over =}} 

\def \di {\diamond}

\def\dom{\mathop{\rm dom}\nolimits}

\def\giantskip{\vskip2\bigskipamount}

\def\gsk{\giantskip}

\def \la {\langle}

\def\msk{\medskip}

\def \ra {\rangle}

\def \res {\!\mid\!\!}

\def\ssubset{\subset\subset}

\def\ssk{\smallskip}

\def\sbr{\smallbreak}

\def\mbr{\medbreak}

\def\bbr{\bigbreak}

\def\giantbreak{\par \ifdim\lastskip<2\bigskipamount \removelastskip

         \penalty-400 \giantskip\fi}

\def\gbr{\giantbreak}

\def\nin{\noindent}

\def\cen{\centerline}

\def\hat{\widehat}

\def\dera #1{{d \over d #1} }

\def\derat#1{{d \over dt} \hbox{\vrule width0.5pt 
                height 5mm depth 3mm${{}\atop{{}\atop{\scriptstyle t=#1}}}$}}

\def\deras#1{{d \over ds} \hbox{\vrule width0.5pt 
                height 5mm depth 3mm${{}\atop{{}\atop{\scriptstyle s=#1}}}$}}

\def\eps{\varepsilon}

\def\phi{\varphi}

\def\epsilon{\varepsilon}

\def\eset{\emptyset}

\def\hb{\hbar}

\def\nin{\noindent}

\def\oline{\overline}

\def\uline{\underline}

\def\pder#1,#2,#3 { {\partial #1 \over \partial #2}(#3)}

\def\pde#1,#2 { {\partial #1 \over \partial #2}}




\def\subeq{\subseteq}

\def\supeq{\supseteq}

\def\sqleq{\sqsubseteq}

\def\sqgeq{\sqsupseteq}

\def\Rarrow{\Rightarrow}

\def\Larrow{\Leftarrow}

\def\tilde{\widetilde}

\def\down{{\downarrow}}

\def\up{{\uparrow}}

\def\what#1{#1\,\,\widehat{ }\,}

\def\Box #1 { \msk\par\nin 
\centerline{
\vbox{\offinterlineskip
\hrule
\hbox{\vrule\strut\hskip1ex\hfil{\smc#1}\hfill\hskip1ex}
\hrule}\vrule}\msk }

 \def\bs{\backslash} 

\title{Holomorphic extensions\\ of representations:\\
 (I) automorphic functions} 
\shorttitle{Holomorphic extensions of representations I} 

 \acknowledgements{The first named author was supported in part by NSF grant DMS-0097314.  The second named author was supported in part
by NSF grant DMS-0070742.}
 \twoauthors{Bernhard Kr\"otz}{Robert J. Stanton}
  
\institutions{The Ohio State University, 
 Columbus, OH\\
{\eightpoint {\it Current address\/}}: University of Oregon, Eugene, OR\\
{\eightpoint {\it E-mail address\/}: kroetz@darkwing.uoregon.edu}
\\ \vglue6pt
The Ohio State University, Columbus, OH\\
{\eightpoint {\it E-mail address\/}:  stanton@math.ohio-state.edu}}
  
 \centerline{\bf Abstract}
\vglue4pt
Let $G$ be a connected, real, semisimple Lie group contained in its 
complexification $G_\C$, and let $K$ be a maximal compact subgroup of $G$. We 
construct a $K_\C$-$G$ double coset domain in $G_\C$, and we show that the 
action of $G$ on the\break $K$-finite vectors of any irreducible unitary 
representation of $G$ has a holomorphic extension to this domain. For the 
resultant holomorphic extension of $K$-finite matrix coefficients we obtain 
estimates of the singularities at the boundary, as well as majorant/minorant 
estimates along the boundary. We\break obtain $L^\infty$ bounds on holomorphically 
extended  automorphic functions on $G/K$ in terms of Sobolev norms, and we use 
these to estimate the  Fourier coefficients of combinations of automorphic 
functions in a number of cases, e.g. of triple products of Maa{\ss} forms. 

\vglue12pt
\intro

Complex analysis played an important role in the classical development of the 
theory of Fourier series. However, even for $\Sl(2,\R)$ contained in 
$\Sl(2,\C)$, complex analysis on  $\Sl(2,\C)$  has had little impact on the 
harmonic analysis of $\Sl(2,\R)$. As the $K$-finite matrix coefficients of an 
irreducible unitary representation of $\Sl(2,\R)$ can be identified with 
classical special functions, such as hypergeometric functions, one knows they 
have holomorphic extensions to some domain. So for any infinite dimensional  
irreducible unitary representation  of $\Sl(2,\R)$, one can expect at most some 
proper subdomain of $\Sl(2,\C)$ to occur. It is less clear that there is a 
universal domain in $\Sl(2,\C)$ to which the action of $G$ on $K$-finite vectors 
of every irreducible unitary representation has holomorphic extension. One goal 
of this paper is to construct such a domain for a real, connected, semisimple 
Lie group $G$ contained in its complexification $G_\C$. It is important to have 
  a maximal domain, and  towards this goal we show that this one is maximal 
in some directions.

Although defined in terms of subgroups of $G_\C$, the domain is natural also 
from the geometric viewpoint. This theme is developed more fully in [KrStII] 
where we show that the quotient of the domain by $K_\C$  is bi-holomorphic to a 
maximal Grauert tube of $G/K$ with the adapted complex structure, and where we 
show that it  also contains a domain bi-holomorphic but not isometric with a 
related bounded symmetric domain. Some implications of this for the harmonic 
analysis of $G/K$ are also developed there.

However, the main goal of this paper is to use the holomorphic extension of 
$K$-finite vectors and their matrix coefficients to obtain estimates involving 
automorphic functions. To our knowledge, Sarnak was the first to use this idea 
in the paper [Sa94]. For example, with it he obtained estimates on the Fourier 
coefficients of polynomials of Maa{\ss} forms for $G = \SO(3,1)$. Sarnak also 
conjectured the size of the exponential decay rate for similar coefficients for 
$\Sl(2,\R)$. Motivated by Sarnak's work,  Bernstein-Reznikov, in [BeRe99], 
verified this conjecture, and in the process  introduced a new technique 
involving $G$-invariant Sobolev norms. As an application of the holomorphic 
extension of representations and with a more representation-theoretic treatment 
of invariant Sobolev norms, we shall verify a uniform version of the conjecture 
for all real rank-one groups. As the representation-theoretic techniques are 
general, we are able also to obtain estimates for the decay rate of Fourier 
coefficients of Rankin-Selberg products of Maa{\ss} forms for $G = \Sl(n,\R)$, and 
to give a conceptually simple proof of results of Good, [Go81a,b], on the growth 
rate of Fourier coefficients of Rankin-Selberg products for co-finite volume 
lattices in  $\Sl(2,\R)$.

It is a pleasure to acknowledge Nolan Wallach's influence on our work by his idea   of viewing automorphic functions as
generalized  matrix coefficients, and to thank Steve Rallis for bringing the 
Bernstein-Reznikov work to our attention, as well as for encouraging us to 
pursue this project. To the referee goes our gratitude for a careful reading of 
our manuscript that resulted in the correction of some oversights, as well as a 
notable improvement of our estimates on automorphic functions for $\Sl(3,\R)$.

\vglue-6pt
\section{The double coset domain}
  \vglue-6pt

To begin we recall some standard structure theory in order to be able to define 
the domain that will be important for the rest of the paper. Any standard 
reference for structure theory, such as [Hel78], is adequate.
  
Let $\g$ be a real, semisimple Lie algebra with a Cartan involution $\theta$.
Denote by $\g=\k\oplus\p$ the associated Cartan decomposition. Take  $\a\subeq 
\p$ a maximal abelian subspace and let $\Sigma=\Sigma(\g, \a)\subeq \a^*$ be the 
corresponding root system. Related to this root system is the root space 
decomposition according to the simultaneous eigenvalues of $\ad (H), H\in \a:$

$$\g=\a\oplus \m \oplus\bigoplus_{\alpha\in \Sigma} \g^\alpha;$$ 
here $\m=\z_\k(\a)$ and $\g^\alpha =\{ X\in \g\: (\forall H\in \a) \ 
[H, X]=\alpha(H) X\}$. For the choice of a positive system $\Sigma^+\subeq 
\Sigma$ 
one obtains the nilpotent Lie algebra $\n =\bigoplus_{\alpha\in 
\Sigma^+}\g^\alpha$. 
Then one has the Iwasawa decomposition on the Lie algebra level 
$$\g=\k\oplus \a \oplus \n.$$

Let $G_\C$ be a simply connected Lie group with Lie algebra $\g_\C$, where 
for a real Lie algebra $\l$, by $\l_\C$ we mean its complexification.  
We denote by $G, A,  A_\C, K, K_\C,  N$ and $N_\C$ the analytic subgroups of 
$G_\C$ corresponding to $\g, \a,\a_\C, \k,\k_\C, \n$ and $\n_\C$. 
If $\uu =\k\oplus i \p$ then it is a subalgebra of $\g_\C$ and the corresponding 
analytic subgroup $U=\exp(\uu)$ is a maximal compact, and in this case, simply 
connected, subgroup of $G_\C$. 

 For these choices one has for $G$ the Iwasawa decomposition, that is, 
the multiplication map 
$$K\times A\times N\to G, \ \ (k, a, n)\mapsto kan$$
is an analytic diffeomorphism. In particular, every element 
$g\in G$ can be written uniquely as $g=\kappa(g) a(g) n(g)$ with each of the 
maps $\kappa(g)\in K$,\break $a(g)\in A$, $n(g)\in N$ depending analytically on $g\in 
G$. 

We shall be concerned with finding a suitable domain in $G_\C$ on which this 
decomposition extends holomorphically. Of course, various domains having this 
property have been obtained by several individuals. What distinguishes the one 
here is its $K_\C$-$G$ double coset feature as well as a type of maximality. 
First we note the following:

\proclaim{Lemma} The multiplication mapping 
$$\Phi\: K_\C\times A_\C\times N_\C\to G_\C, \ \ (k,a,n)\mapsto kan$$
has everywhere surjective differential. 
\endproclaim

\demo{Proof}  Obviously one has $\g_\C=\k_\C\oplus\a_\C\oplus\n_\C$
and $\a_\C\oplus\n_\C$ is a subalgebra of $\g_\C$. Then following 
Harish-Chandra, since $\Phi$ is left $K_\C$ and right $N_\C$-equivariant it 
suffices to check that 
$d\Phi(\1, a, \1)$ is surjective for all $a\in A_\C$. Let $\rho_a(g)=ga$ be the 
right translation in $G_\C$ by the element $a$. Then for $X\in\k_\C$, $Y\in 
\a_\C$ 
and $Z\in \n_\C$ one has
$$d\Phi(\1,a, \1)(X,Y,Z)= d\rho_a(\1)(X+Y+\Ad(a)Z),$$
from which the surjectivity follows. 
\enddemo

 To describe the domain we extend $\a$ to a $\theta$-stable Cartan 
subalgebra $\h$ of $\g$ so that $\h=\a\oplus\t$ with 
$\t\subeq \m$. Let $\Delta=\Delta(\g_\C,\h_\C)$ be the corresponding root system 
of $\g$. Then it is known that $\Delta\res_\a\bs \{0\}=\Sigma$.   

 Let $\Pi=\{\alpha_1, \ldots, \alpha_n\}$ be the set of simple restricted 
roots corresponding to the positive roots 
$\Sigma^+$. We define elements $\omega_1, \ldots, \omega_n$ of $\a^*$ as 
follows, 
using the restriction of the Cartan-Killing form to $\a$: 
$$(\forall 1\leq i,j\leq n)\qquad \left\{ \begin{array}{ll}\la \omega_j,\alpha_i\ra =0 & \hbox{if $i\neq 
j$}\\[4pt] {2\la \omega_i,\alpha_i\ra\over \la \alpha_i,\alpha_i\ra}=1 &\hbox{if 
$\alpha_i\in\Delta$}\\[4pt] {\la \omega_i,\alpha_i\ra\over \la \alpha_i,\alpha_i\ra}=1 &\hbox{if 
$\alpha_i\not\in\Delta$ and $2\alpha_i\not\in 
\Sigma$}\\[4pt] {\la \omega_i,\alpha_i\ra\over \la \alpha_i,\alpha_i\ra}=2 & \hbox{if 
$\alpha_i\not\in\Delta$ and $2\alpha_i\in \Sigma$.} \end{array}\right.$$ 
Using standard results in structure theory relating $\Delta$ and $\Sigma$ one 
can show that $\omega_1, \ldots, \omega_n$ are algebraically integral for 
$\Delta=\Delta(\g_\C,\h_\C)$. 
The last piece of structure theory we shall recall is the little Weyl group. We 
denote by ${\cal W}_\a =N_K(\a)/ Z_K(\a)$ the {\it Weyl group} of 
$\Sigma(\a,\g)$.

We are ready to define a first approximation to the double coset domain. We 
set 
$$\a_\C^1=\{ X\in \a_\C \: (\forall 1\leq k\leq n)(\forall w\in {\cal W}_\a)\ 
|\Im \omega_k(w.X)|<{\pi\over 4}\}$$
and 
$$\a_\C^0=  2\a_\C^1.$$ On the group side we let
$A_\C^0=\exp(\a_\C^0)$ and $A_\C^1=\exp(\a_\C^1)$. 
Clearly ${\cal W}_\a$ leaves each of $\a_\C^0$, $\a_\C^1$, $A_\C^0$ and $A_\C^1$ 
invariant.

If $\alpha\in \a_\C^*$ is analytically integral for $A_\C$, then we set 
$a^\alpha=e^{\alpha(\log a)}$ for all $a\in A_\C$. Since $G_\C$ is simply 
connected, the elements $\omega_j$ are analytically integral for 
$A_\C$ and so we have $a^{\omega_k}$ well defined.

Next we introduce the domains
\begin{eqnarray*}
A_\C^{0,\leq} &\nhs=\nhs&\{a\in A_\C\: (\forall 1\leq k\leq  n) \Re 
(a^{\omega_k})>0\},\\
\noalign{\noindent and} 
 A_\C^{1,\leq}&\nhs =\nhs&(A_\C^{0,\leq})^{1\over 2}= \{a\in A_\C\: (\forall 1\leq k \leq 
 n)
|\arg(a^{\omega_k})|< {\pi\over 4}\}.\end{eqnarray*}
Note that $A_\C^0\subeq A_\C^{0,\leq}$ and $A_\C^1\subeq A_\C^{1,\leq}$.

\proclaim{Lemma}
{\rm (i)} For $\Omega\subeq A_\C$ open{\rm ,} $K_\C \Omega N_\C$ is open in $G_\C$. 
In particular{\rm ,} 
the sets $K_\C A_\C N_\C$, $K_\C A_\C^1 N_\C$, $K_\C A_\C^{1,\leq}N_\C$, 
$K_\C A_\C^0N_\C$ and $K_\C A_\C^{0,\leq} N_\C$
are open in $G_\C$. 

 {\rm (ii)} $K_\C A_\C N_\C$ is   dense in $G_\C$. 
\endproclaim

\demo{Proof} This is an immediate consequence of Lemma 1.1 as $\Phi$ is 
a morphism of affine algebraic varieties with everywhere submersive
differential.\enddemo

 \proclaim{Proposition}  Let $G_\C$ be a simply connected{\rm ,} semisimple{\rm ,} complex Lie 
group. Then the multiplication mapping 
\begin{eqnarray*}
&&\Phi\: K_\C\times A_\C^{0,\leq}\times N_\C\to G_\C,\ \ (k,a,n)\mapsto kan\\
\noalign{\noindent 
is an analytic diffeomorphism onto its open image $K_\C A_\C^{0,\leq} N_\C$. }
\end{eqnarray*}
\endproclaim

\vglue-24pt
{\it Proof}.  In view of the preceding lemmas, it suffices to show that $\Phi$ is 
injective. Suppose that 
$kan=k'a'n'$ for some $k, k'\in K_\C$, $a, a'\in A_\C^{0,\leq}$ and $n, n'\in 
N_\C$. 
Denote by $\Theta$ the holomorphic extension 
of the Cartan involution of $G$ to $G_\C$. Then we get that 
$$\Theta(kan)^{-1}kan=\Theta(k'a'n')^{-1}k'a'n'$$
or equivalently
$$\Theta(n^{-1})a^2 n=\Theta((n')^{-1}) (a')^2 n'.$$ 
Now the subgroup $\oline N_\C=\Theta(N_\C)$ corresponds to the analytic 
subgroup with Lie algebra
$\oline \n_\C=\bigoplus_{\alpha\in -\Sigma^+} \g_\C^\alpha$. 
As a consequence of the injectivity of the map 
$$\oline N_\C \times A_\C\times N_\C\to \oline N_\C A_\C N_\C, \ \ (\oline n, a, 
n )\mapsto 
\oline n a n$$
we conclude that $n=n'$ and $a^2=(a')^2$. We may assume that $a, a'\in 
\exp(i\a)$.  

To complete the proof of the proposition it remains to show that 
$a^2=(a')^2$ for $a,a'\in  A_\C^{0, \leq}$ implies that $a=a'$. Let $X_1, 
\ldots, X_n$ in $\a_\C$ be the dual basis to $\omega_1, \ldots, \omega_n$. We 
can write 
$a=\exp(\sum_{j=1}^n \phi_j X_j)$
and $a'=\exp(\sum_{j=1}^n \phi_j'X_j)$
for complex numbers $\phi_j$, $\phi_j'$ satisfying $|\Im\phi_j|<{\pi\over 2}$, 
$|\Im\phi_j'|<{\pi\over 2}$. 
Then $a^2=(a')^2$ implies that 
$$e^{2\phi_j}=a^{2\omega_j}=(a')^{2\omega_j}=e^{2\phi_j'}$$
and hence $\phi_j=\phi_j'$ for all $1\leq j\leq n$, concluding the proof of the proposition.
\phantom{running} \hfill\qed\vglue12pt

Thus every element $z\in K_\C
A_\C^{0,\leq} N_\C$ can be uniquely written as $z=\kappa(z) a(z) n(z)$ with $\kappa(z)\in K_\C$, $a(z)\in 
A_\C^{0,\leq}$ and $n(z)\in N_\C$ all depending holomorphically on $z$.   
Next we define domains using the restricted roots. We set
$$\b^0=\{ X\in\a\:(\forall \alpha\in \Sigma)\  |\alpha(X)|<\pi\}.$$
and
$$\b^1={1\over 2}\b^0.$$ Clearly 
both $\b^0$ and $\b^1$ are ${\cal W}_\a$-invariant. We set $\b_\C^j=\a + i\b^j$ 
and 
$B_\C^j=\exp(\b_\C^j)$ for $j=0,1$. 
Let $\a^0=i(\a_\C^0\cap i\a)$. Then, from the classification of restricted root 
systems and standard facts about the 
associated fundamental weights, one can verify that $\a^0\subeq \b^0$. For a 
comparison of these domains we provide below the illustrations for two rank 2 
algebras. 

 \proclaim{Lemma} Let $\omega \subeq i\b^1$ be a nonempty{\rm ,} open{\rm ,} ${\cal 
W}_\a$\/{\rm -}\/invariant{\rm ,}
 convex set. Then the set 
$$K_\C \exp(\omega)G$$
is open in $G_\C$. 
\endproclaim

\begin{center}
\BoxedEPSF{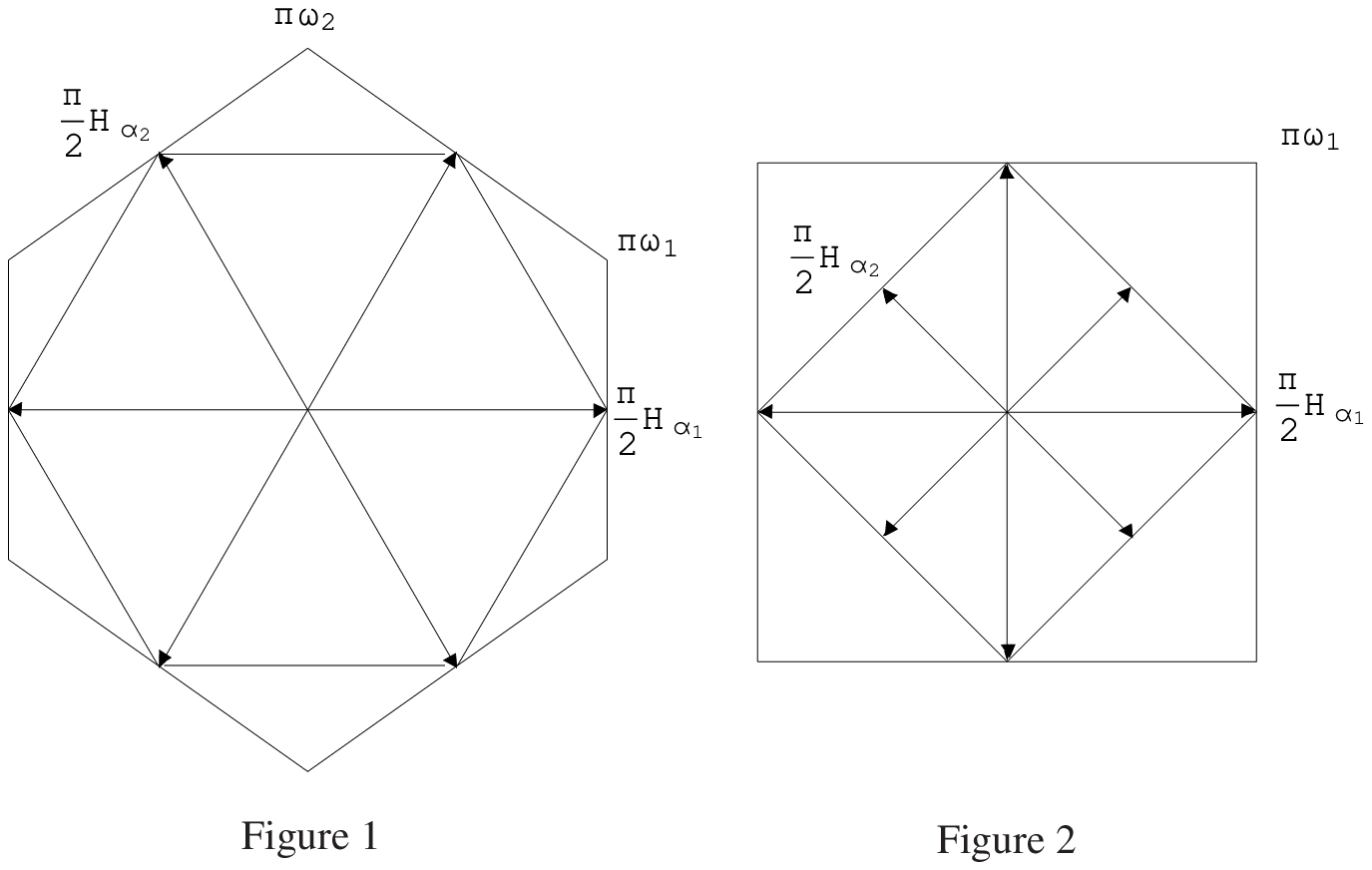 scaled 750}
 \end{center}
 
\begin{quote}
{Figure 1 corresponds to $\sL(3,\R)$ and Figure 2
to $\sp(2,\R)$. The region enclosed by an outer polygon corresponds
to $\b^0$ while that enclosed by an inner polygon corresponds 
to $\a^0$. The $H_{\alpha_i}$ denote the coroots of $\alpha_i$ and 
we identify the $\omega_i$ as elements of $\a$ via the Cartan-Killing form. }
\end{quote}

\demo{Proof} Set $W=\Ad(K)\omega$. Since $\omega$ is open, convex, and ${\cal 
W}_\a$-invariant, 
Kostant's nonlinear convexity theorem shows that $W$ is an open, convex set in~$i\p$. 
Note that $K_\C \exp(\omega )G=K_\C \exp(W) G$. 
Now [AkGi90, p.\ 4-5] shows that the  multiplication mapping 
$$m\: K_\C\times \exp(W)\times G\to G_\C, \ \ (k, a, g)\mapsto kag$$
has everywhere surjective differential. {}From that 
the assertion follows. \enddemo

For each $1\leq k\leq n$ we write $(\pi_k, V_k)$ for the real, finite-dimensional, 
highest weight representation of $G$ with highest weight $\omega_k$. 
We choose a scalar product $\la\cdot, \cdot\ra$ on $V_k$ which satisfies 
$\la \pi_k(g)v, w\ra =\la v, \pi_k(\Theta(g)^{-1})w\ra$ for all 
$v,w\in V_k$ and $g\in G_\C$. We denote by $v_k$ a normalized highest weight 
vector of $(\pi_k, V_k)$.

 \proclaim{Lemma}   For all $1\leq k\leq n${\rm ,} $a\in A_\C^1$ and $m\in \oline N${\rm ,}
$$\Re \big(\la \pi_k(\theta(m)^{-1} a^2 m)v_k, v_k\ra\big)>0.$$
\endproclaim

\demo{Proof} Fix $1\leq k\leq n$, $a$ and $m\in \oline N$, and note that $a^2\in A_\C^0$. Now,
$$\la \pi_k(\theta(m)^{-1} a^2 m)v_k, v_k\ra=\la\pi_k(a^2)\pi_k(m)v_k, 
\pi_k(m)v_k\ra.\leqno(1.1)$$
Let ${\cal P}_k\subeq \a^*$ denote the set of $\a$-weights of $(\pi_k, V_k)$. 
Then (1.1) implies that there exist  nonnegative numbers $c_\beta$, $\beta\in 
V_k$, such that 
$$\la\pi_k(\theta(m)^{-1} a^2 m)v_k, v_k\ra=\sum_{\beta\in {\cal P}_k} c_\beta 
a^{2\beta}.$$
Recall that 
$${\cal P}_k\subeq \conv ({\cal W}_\a\omega_k).$$
Since $\a_\C^0$ is convex and Weyl group invariant, to finish the proof it 
suffices to show that 
$\Re(a^{2\omega_k})>0$ for all $a\in A_\C^1$. But this is immediate from the
definition of $\a_\C^1$.  \enddemo

 \proclaim{Lemma}   Let $(b_j)_{j\in \N}$ be a convergent 
sequence in $A_\C$ and $(n_j)_{j\in \N}$ an unbounded sequence in $N_\C$. 
Then the sequence 
$$\big(\Theta(n_j)^{-1} b_j n_j\big)_{j\in \N}$$
is unbounded in $G_\C$. 
\endproclaim

\demo{Proof} Let $d(\cdot, \cdot)$ be a left invariant metric on $G_\C$. Then 
$$d(\Theta(n_j)^{-1} b_j^2 n_j, \1)=d(b_j^2 n_j, \Theta(n_j)),$$
and we see that $\lim_{j\to \infty} d(\Theta(n_j)^{-1} b_j^2 n_j, \1)=\infty $ 
(this follows for example by embedding $\Ad(G_\C)$ into $\Sl(m, \C)$, where we 
can arrange matters
so that $A_\C$ maps into the diagonal matrices and $N_\C$ in the upper 
triangular matrices). 
 \enddemo

 \proclaim{Proposition} {\rm (i)} $K_\C A_\C^1 G$ is open in $G_\C$. 

\vglue4pt {\rm (ii)} $K_\C A_\C^1 G\subeq K_\C A_\C^{1,\leq} N_\C.$

\vglue4pt {\rm (iii)} For all $\lambda\in \a_\C^*$ the mappings 
$$A_\C^1\times G\to \C, \ \ (a,g)\mapsto a(ag)^\lambda,$$
$$A_\C^1\times G\to K_\C, \ \ (a,g)\mapsto \kappa(ag)$$
are analytic{\rm ,} and holomorphic in the first variable. 
\endproclaim

\demo{Proof} (i) appears in Lemma 1.2. For (ii) take an $a\in A_\C^1$. First we show 
that $a\oline N\subeq K_\C A_\C N_\C$. 
Fix $m\in \oline N$ and let 
\begin{eqnarray*}
\Omega & =&\{ a\in A_\C^1\: am\in K_\C A_\C N_\C\}\\ &=&\{a\in A_\C^1\: \Theta(m)^{-1} a^2 m\in \oline N_\C A_\C
N_\C\}.
\end{eqnarray*}
 Then $\Omega$ is open and nonempty. We have to show that 
$\Omega=A_\C^1$. Suppose the contrary. Then there exists a 
sequence $(a_j)_{j\in \N}$ in $\Omega$ such that 
$a_0=\lim_{j\to \infty} a_j\in A_\C^1\bs \Omega$. 

Let $a\in \Omega$. Then by Proposition 1.3 we find unique
$k\in K_\C$, $b\in A_\C$ and $n\in N_\C$ such that 
$am=kbn$ or, in other words, 
$$\Theta(m)^{-1} a^2 m=\Theta(n)^{-1} b^2 n.$$
Taking matrix-coefficients with fundamental representations we thus get that 
$$b^{2\omega_k}=\la\pi_k(\Theta(n)^{-1} b^2 n)v_k, v_k\ra =\la 
\pi_k(\Theta(m)^{-1} a^2 m)
v_k, v_k\ra\leqno(1.2) $$
for all $1\leq k\leq n$. 
Applied to our sequence $(a_j)_{j\in \N}$ we get elements $k_j\in K_\C$, $b_j\in 
A_\C$ and 
$n_j\in N_\C$ with $a_jm=k_jb_jn_j$. Lemma 1.5 together with  
(1.2) imply that $(b_j)_{j\in\N}$ is bounded. If necessary, by taking a 
subsequence, we may assume that $b_0=\lim_{j\to \infty} b_j$ exists in $A_\C$. 
Since $\Theta(m)^{-1} a_0^2 m\not \in \oline N_\C A_\C N_\C$, the sequence 
$(n_j)_{j\in \N}$
is unbounded in $N_\C$. Hence $\big(\Theta(n_j)^{-1} b_jn_j\big)_{j\in \N}$ is an 
unbounded 
sequence in $G_\C$ by Lemma 1.6. 
But this contradicts the fact that $\big(\Theta(m)^{-1} a_j^2 m\big)_{j\in \N}$ 
is bounded. 
Thus we have proved that  $a\oline N\subeq K_\C A_\C N_\C$ for all 
$a\in A_\C^1$. But now (1.2) together with Lemma 1.5 actually shows that 
$b\in A_\C^{1,\leq}$, hence  $a\oline N\subeq K_\C A_\C^{1,\leq} N_\C$ for all 
$a\in A_\C^1$. 
\ssk  The Bruhat decomposition of $G$ gives
$G=\bigcup_{w\in {\cal W}_\a} \oline N wMAN$
with $M=Z_K(A)$. Since $A_\C^1$ is $N_K(A)$-invariant, we get that 
$aG\subeq K_\C A_\C^{1,\leq} N_\C$. 
Then (ii) is now clear while (iii) is a consequence of (ii) and Proposition 
1.3. \enddemo

Next we are going to prove a significant extension of Proposition 1.7. We will 
conclude the proof in the following section.

\proclaim{Theorem} Let $G$ be a classical semisimple Lie group. Then the following 
assertions 
hold\/{\rm :}\/
\begin{itemize}
\ritem{(i)} $K_\C B_\C^1 G$ is open in $G_\C;$
\ritem {(ii)} $B_\C^1 G\subeq K_\C A_\C N_\C;$
\ritem {(iii)} there exists  an analytic  function 
$$B_\C^1\times G\to \a_\C, \ \ (a,g)\mapsto H(ag),$$
holomorphic in the first  variable{\rm ,}
such that $ ag\in K_\C \exp H(ag) N_\C$ for all $a\in B_\C^1$ and $g\in G${\rm ;}
\ritem{(iv)} there exists  an analytic  function 
$$\kappa\: B_\C^1\times G\to K_\C, \ \ (a,g)\mapsto \kappa(ag),$$
holomorphic in the first  variable{\rm ,}
such that $ag\in \kappa(ag)A_\C N_\C$ for all $a\in B_\C^1$ and $g\in G$. 
\end{itemize}

\endproclaim

\demo{Proof} (i) follows from Lemma 1.2. (ii) follows from Proposition 2.5, 
Proposition 2.6 and Proposition 2.9 in the next section. 

\vglue4pt
(iii) Set $L=K_\C\cap A_\C$ and note that $L$ is a
discrete subgroup of $G_\C$. 
Then the first part of the proof of Lemma 1.3 shows that we  have a 
biholomorphic diffeomorphism 
$$(K_\C\times_L A_\C)\times N_\C \to K_\C A_\C N_\C, \ \ ([k,a],n)\mapsto kan.$$
In particular, we get a holomorphic middle projection 
$$\tilde a\: K_\C A_\C N_\C \to A_\C/ L, \ \ kan\mapsto aL, $$  
and so, by (ii), an analytic mapping 
$$\tilde \Phi\:   B_\C^1\times G\to A_\C/L, \ \ (a,g)\mapsto \tilde a(ag).$$
Now $\a_\C\to A_\C/ L$, via the map $\ X\mapsto \exp(X)L$, is the universal 
cover of $A_\C/L$. 
To complete the proof of (iii) it remains to show that $\tilde \Phi$ lifts 
to a continuous map with values in $\a_\C$. Since $\exp\: \a_\C^1\to A_\C^1$ is 
injective, Proposition 1.7
implies that $\tilde \Phi\res_{A_\C^1\times G}$ lifts to a continuous map $\Psi$ 
 with values in $\a_\C$. 
Since the exponential function restricted to $\b_\C^1$ is injective (cf. Remark 
1.9.), 
$B_\C^1$ is simply connected and so for every simply connected set $U\subeq G$ 
we get a continuous lift of $\tilde \Phi\res_{ B_\C^1\times U }$ extending 
$\Psi\res_{A_\C^1\times U}$. By the uniqueness of liftings 
we get a continuous lift of $\tilde \Phi$ completing the proof of 
(iii).

\vglue4pt
(iv) In view of (ii), we get an analytic map 
$$\tilde\kappa\: B_\C^1\times G\to K_\C/ L,\ \  (a,g)\mapsto \tilde \kappa(ag)$$
even holomorphic in the first variable and such that $ag\in \tilde\kappa(ag)A_\C 
N_\C$. 
Thus in order to prove the assertion in (iv), it suffices that $\tilde \kappa$ 
lifts to a continuous 
map $\kappa\: B_\C^1\times G\to K_\C$. But this is proved as in (iii).
 \enddemo

\numbereddemo{{R}emark}    The simply connected hypothesis on $G_\C$  that has been made is 
not necessary. More generally, if $G$ is classical, semisimple and contained in 
its complexification, then  Theorem 1.8 is valid. Indeed, let $\g$ be a 
semisimple Lie algebra with Cartan 
decomposition $\g=\k\oplus \a\oplus\n$, $\g_\C$ its complexification and let
$G_\C$ be a simply connected Lie group with Lie algebra $\g_\C$. As before, let 
$G$ be the analytic subgroup of $G_\C$ with Lie algebra $\g$.

Let now $G_1$ be another connected Lie group with Lie algebra $\g$ and 
suppose that $G_1$ sits in its complexification $G_{1,\C}$. 
Write $G_1=K_1 A_1 N_1$ for the Iwasawa decomposition of $G_1$ corresponding 
to $\g=\k\oplus \a\oplus\n$. Set $B_{1,\C}^1=A_1 \exp_{G_{1,\C}} (i\b^1)$. 
Since $G_\C$ is simply connected, we have a covering homomorphism 
$$\pi\: G_\C \to G_{1,\C}.$$
Hence Theorem 1.8 (ii) implies that 
$$B_{1,\C}^1 G_1\subeq K_{1,\C} A_{1,\C} N_{1,\C}.$$
To see that Theorem 1.8 (iii), (iv) remains true for $G_1$ contained in 
$G_{1,\C}$ one needs that $B_{1,\C}^1$ is simply connected. But this will follow 
from
the fact that $\exp_{G_{1,\C}}\: \b_\C^1\to B_{1,\C}^1$ is injective. 
To see this, note that this map is injective if and only if the map 
$$f\: \b^1\to A_{1,\C}, \ X\mapsto \exp_{G_{1,\C}}(X)$$
is injective. If $f$ were not injective, then there would exist  an element 
$X\in \b^0$,\break $X\neq 0$, such that $\exp_{G_{1,\C}}(X)=\1$. Hence $\alpha(X)\in 
i{2\pi\Z}$
for all $\alpha\in \Sigma$\break (cf.\ [Hel78, Ch.\ VII,\S4, Prop.\ 4.1]), a
contradiction to $X\in \b^0\backslash \{0\}$.  \enddemo

The next proposition will be used in a later section. It has independent 
interest as it can be considered as a principle of convex inclusions 
and as such is related to Kostant's nonlinear convexity theorem. 

 Suppose that $E$ is a subset in a complex vector space $V$. We 
denote by $\conv E$ the convex hull of $E$ and by
$\cone E =\R^+ E$ the cone generated by $E$.   

 \proclaim{Proposition} Let $0\in \omega\subeq \b^0$ be a connected subset. Set 
$\b_\C^\omega =\a+ i \omega$ 
and $B_\C^\omega =\exp(\b_\C^\omega)$. Then{\rm ,}
$$ B_\C^\omega G\subeq K_\C A_\C N_\C \Rarrow   B_\C^{\conv\omega} G\subeq K_\C 
A_\C N_\C.$$ 
\endproclaim

\demo{Proof} Fix $g\in G$. It suffices to show the existence of a holomorphic 
function 
$$f_g\: B_\C^{\conv\omega} \to\a_\C, \ \ a\mapsto f_g(a)$$
such that $ag\in K_\C \exp(f_g(a)) N_\C$ for $a\in B_\C^{\conv\omega}$ holds. 
We already know from Theorem 1.8(iii) 
that a holomorphic function $\tilde f_g\: B_\C^\omega \to\a_\C$ with $ag\in K_\C 
\exp(\tilde f_g(a)) N_\C$
for $ a\in B_\C^\omega$ exists. Now $B_\C^{\conv\omega}$ is the holomorphic hull 
of $B_\C^\omega$
and so $\tilde f_g$ extends to a holomorphic mapping $f_g\: B_\C^{\conv\omega} 
\to\a_\C$. 

It remains to show that $ag\in K_\C \exp(f_g(a)) N_\C$ for $a\in 
B_\C^{\conv\omega}$. If not, then 
we find a convergent sequence $(a_n)_{n\in\N}$ with $\lim_{n\to\infty} 
a_n=a_0\in B_\C^{\conv\omega}$, 
$a_ng\in K_\C A_\C N_\C$ but $a_0g\not\in K_\C A_\C N_\C$. 
Hence we find a sequence $m_n\in N_\C$ such that 
$$\Theta(g)^{-1} a_n^2 g= \Theta(m_n)^{-1} f_g(a_n)^2 m_n$$
but $\Theta(g)^{-1}a_0^2 g\not\in \oline N_\C A_\C N_\C$. As 
$\big(f_g(a_n)\big)_{n\in\N}$ is bounded, we conclude (cf.\ Lemma 1.6) that 
$(m_n)_{n\in\N}$ is unbounded, a contradiction.  \enddemo

\section{Matrix calculations}

We shall prove (ii) of Theorem 1.8 by various results about matrices. First we 
shall treat the group $G=\Sl(m,\R)$, $m\geq 2$. Then we shall give a class of 
subgroups of $\Sl(m,\R)$ whose roots have a hereditary property similar to one 
held by Levi factors of parabolic subgroups. This will allow us to take care of 
most of the classical groups. The remaining cases are treated at the end of this 
section.

Here we obviously have $G\subeq G_\C=\Sl(m,\C)$ with $G_\C$ simply connected. We 
let $\k=\so(m,\R)$ and choose   

\centerline{${\displaystyle\a =\{ \diag (x_1, x_2, \ldots, x_m)\in M(m,\R)\: \sum_{i=1}^m x_i=0\}}$}
\vglue12pt\noindent 
as a maximal abelian subalgebra in $\p=\Symm(m,\R)\cap\sL(m,\R)$. 
Define elements $\eps_j\in \a^*$ by 
setting 
$$\eps_j(\diag (x_1, \ldots, x_m))=x_j.$$
Then $\Sigma=\{ \eps_i-\eps_j\: 1\leq i\neq j\leq m\}$ and we 
take $\Sigma^+ =\{ \eps_i-\eps_j\: i<j\}$ as a positive system. The associated 
system of simple restricted roots is given by 
$$\Pi=\{ \eps_1-\eps_2, \ldots, \eps_{m-1}-\eps_m\}.$$
As $\g$ is {\it split} we have $\Sigma=\Delta$. In particular, the $\omega_j$, 
$1\leq 
j\leq m-1$ are the usual fundamental weights and are given by 
$$\omega_j=\eps_1+\ldots+\eps_j, \qquad (1\leq j\leq m-1).$$ 
The Weyl group ${\cal W}_\a$ of $\Sigma(\a,\g)$ is 
the group of permutations on the  $m$ elements $\eps_1, \ldots, 
\eps_m$.

 In matrix notation the nilpotent groups $N$ and $\oline N$ are given by:
$$
N=\{\left( \begin{array}{ccccc}1& x_{12}&\ldots&  x_{1m}\\ & 1 &x_{23}  \ldots &x_{2m}\\ &   & \ddots &\vdots       \\ &   &        & 1    &
\end{array}\right)\: x_{ij}\in\R\}$$
and
$$
\oline N=\{\left( \begin{array}{cccc}
1      &        &             &\\ x_{21} & 1      &             &\\ \vdots & \ddots & \ddots      &\\ x_{m1} &\ldots  & x_{m,m-1}& 1\end{array}\right)\:
x_{ij}\in\R\}.$$

 For each $1\leq j\leq m$ we set $e_j=(\delta_{k-j, 
l-j})_{k,l}\in\diag(m, \R)$. 
Further we associate to each $\omega_j$ the element $X_{\omega_j}=\sum_{k=1}^j 
e_j 
-{j\over m} \sum_{j=1}^m e_j$.

 \proclaim{Lemma} {\rm (i)} 
$$\b^0=\Int \Big(\conv\big(\{ \pm \pi w.X_{\omega_j}\: w\in {\cal W}_\a, \ 1\leq 
j\leq m-1\}\big)\Big);$$

{\rm (ii)}  
$$\b^0\subeq \a\cap \Big({m-1\over m}\bigoplus_{j=1}^m ]-\pi, \pi[ e_j\Big).$$
\endproclaim

\demo{Proof} (i) Set 
$$\b' =\Int\Big( \conv\big(\{ \pm \pi w X_{\omega_j}\: w\in {\cal W}_\a, \ 1\leq 
j\leq m-1\}\big)\Big)
.$$
Both $\oline {\b^0}$ and $\oline {\b'}$ are closed, convex and ${\cal 
W}_\a$-invariant.
Thus by the convexity of 
$\b^0$ and $\b'$ we have to show only that  $\oline {\b^0}= \oline {\b'}$. 
Now 
${\cal W}_\a$ rotates the extreme points of both  $\oline {\b^0}$ and $\oline 
{\b'}$,
and the extreme points of $\oline {\b'}$ are given by 
$\pm w. \pi X_{\omega_j}$. We shall prove the result by double containment.

``$\supeq$'': By the Krein-Milman Theorem it suffices to show that $\pm \pi X_{\omega_j}\in \oline {\b^0}$ for all $1\leq j\leq m-1$. Every 
$\alpha\in \Sigma^+$ can be written as 
$\alpha=\sum_{j=1}^{m-1} \delta_j (\eps_j -\eps_{j+1})$ with coefficients 
$\delta_j\in \{0,1\}$. Thus $\alpha(X_{\omega_j}) \in \{0,1\}$ and the 
inclusion ``$\supeq$'' follows from the definition of $\b^0$. 

``$\subeq$'': Notice that $\omega_1, \ldots, \omega_{m-1}$ 
constitute  a basis of $\a^*$. Hence every $X\in\oline {\b^0}$ can be written 
as $X=\sum_{j=1}^{m-1} \lambda_j X_{\omega_j}$ with coefficients 
$\lambda_j\in\R$. 
{}From the definition of $\b'$ we may assume that $\lambda_j\geq 0$ for all 
$1\leq j\leq m-1$. In particular, we see that  
$$(\eps_1-\eps_m)(X) =\sum_{j=1}^{m-1}(\eps_j-\eps_{j+1})(X) =
\sum_{j=1}^{m-1}\lambda_j\in [0, \pi[, $$
concluding the proof of ``$\subeq$''.

(ii) For $X=\diag(x_1, \ldots, x_m)=\sum_{j=1}^m x_j e_j\in \b^0$,
$$-\pi< 2x_1+x_2+\ldots+x_m =x_1-x_m<\pi ,$$
and 
$$-\pi<x_1-x_j<\pi \qquad \hbox{for all $2\leq j\leq m-1$}.$$  
By summing these inequalities we obtain 
$$-(m-1)\pi <mx_1 <(m-1)\pi ,$$
or equivalently $ |x_1|< {m-1\over m}\pi$. Similarly, 
 $|x_j|<{m-1\over m}\pi$ for all $1\leq j\leq m$.  \enddemo

\numbereddemo{{R}emark}  Notice that $\a^0$ is strictly smaller than $\b^0$,
although  they have common boundary points (cf.\ Figure 1). In particular, Lemma 2.1 shows 
that 
$$\partial \a^0\cap \partial \b^0\supeq \{ {\pi\over 2} (e_i-e_j)\: 
1\leq i\neq j\leq m-1\}. 
$$
\enddemo

For every $1\leq k\leq m$ we denote by $\Delta_k(A)$ the $k^{\rm th}$ principal
minor of a matrix $A\in M(m,\C)$. 
For every $g=(g_{ij})_{1\leq i, j\leq m}\in M(m,\C)$ and $1\leq k\leq m$
we define $g_{(k)}\in  M(k,\C)$ by $g_{(k)}=(g_{ij})_{1\leq i, j\leq k}$.

 \proclaim{Proposition}  Let $G=\Sl(m,\R)$ with $G_\C=\Sl(m,\C)$. Then for all  $1\leq 
k\leq m${\rm ,}  
$a\in B_\C^0$ and $g\in \Sl(m,\R)$ with $g_{(k)}\in \Gl(k,\R)${\rm ,}   \begin{itemize}
\ritem{(i)} $\Delta_k(gag^t)\neq 0;$
\ritem{(ii)} $\Spec \big((gag^t)_{(k)}\big)\subeq
\cone\big ( \conv (\Spec(a))\big)$.
\end{itemize}

\endproclaim

\demo{Proof} (i) Fixing $1\leq k\leq m$,  $a\in B_\C^0$ and 
$g\in\Sl(m,\C)$ with\break $g_{(k)}\in \Gl(k,\C)$, we write
$a=\diag(r_1 e^{i\phi_1}, \ldots, r_me^{i\phi_m})$ with 
$r_i>0$, $-{m-1\over m}\pi < \phi_i < {m-1\over m}\pi$ (cf.\ Lemma 2.1(ii)).  
Set  

$$g=\left( \begin{array}{cc} g_{(k)} & B \\ \ast &\ast\end{array}\right)$$
with $g_{(k)}\in \Gl(k, \R)$ and $B\in M(k\times (m-k), \R)$. 
Then,
\begin{eqnarray*}
\Delta_k(gag^t)&=&
\Delta_k\Big(\left( \begin{array}{cc} g_{(k)} & B \\ \ast&\ast
 \end{array}\right)\diag(r_1e^{i\phi_1}, \ldots r_me^{i\phi_m})\left( \begin{array}{cc}
g_{(k)}^t & *\\ B^t & *\end{array}\right)\Big)\\
&=&\det_k\Big(g_{(k)}\diag(r_1e^{i\phi_1},\ldots, r_ke^{i\phi_k})g_{(k)}^t\\ &&\hskip.5in+
B\diag(r_{k+1}e^{i\phi_{k+1}},\ldots, r_me^{i\phi_m})B^t\Big).
\end{eqnarray*}
In order to show that $\Delta_k(gag^t)\neq 0$ we have to 
show that the $k\times k$-matrix 
$$X_{(k)} =g_{(k)}\diag(r_1e^{i\phi_1},\ldots, r_ke^{i\phi_k})g_{(k)}^t+
B\diag(r_{k+1}e^{i\phi_{k+1}},\ldots, r_me^{i\phi_m})B^t$$
is invertible.

 Assume first that $k\leq m-k$. Then we can write $B=(B_1,B_2)$ with 
$B_1\in M(k,\R)$ and $B_2\in M(k\times (m-2k),\R)$. Hence we obtain that 
\begin{eqnarray*}
B\diag(r_{k+1}e^{i\phi_{k+1}},\ldots, r_me^{i\phi_m})B^t&=&
B_1\diag(r_{k+1}e^{i\phi_{k+1}},\ldots, r_{2k}e^{i\phi_{2k}})B_1^t\\ && + B_2\diag(r_{2k+1}e^{i\phi_{2k+1}},\ldots,
r_me^{i\phi_m})B_2^t.
\end{eqnarray*}

Let $\la\cdot, \cdot\ra$ be the 
usual hermitian inner product on $\C^k$.  In particular, if $v\in \C^k$, $v\neq 
0$, then we get 
\begin{eqnarray*}
\la X_{(k)}v,v\ra&=&\la\diag(r_1e^{i\phi_1},\ldots, 
r_ke^{i\phi_k})g_{(k)}^tv, g_{(k)}^tv\ra\\[5pt]
&&
+\la \diag(r_{k+1}e^{i\phi_{k+1}},\ldots, r_{2k}e^{i\phi_{2k}})B_1^tv, B_1^tv\ra 
\\[5pt] && +\la \diag(r_{2k+1}e^{i\phi_{2k+1}},\ldots, r_me^{i\phi_m})B_2^tv, 
B_2^tv\ra.\end{eqnarray*}
So there exist numbers $c_1, \ldots, c_m\geq 0$,  not all zero,  such that 
$$\la X_{(k)}v,v\ra=\sum_{j=1}^m c_j e^{i\phi_j}.\leqno(2.1)$$
Similarly one shows that (2.1) holds for the case $k\geq m-k$. 
Now (i) follows from (2.1) and Lemma 2.4 below.

(ii) Since $X_{(k)}= (gag^t)_{(k)}$, (ii) follows from (2.1). 
 \enddemo

We denote by $\C^+ =\{ z\in \C\: z\not\in ]-\infty, 0]\}$ the split plane  in 
$\C$.

 \proclaim{Lemma}  Let $\phi_1, \ldots, \phi_m\in\R$ be such that 
$\diag(\phi_1,\ldots, \phi_m)\in \b^0$. Then for all 
sequences of nonnegative numbers $c_1, \ldots, c_m${\rm ,} not all 
zero{\rm ,}  
$$\sum_{j=1}^m c_j e^{i\phi_j}\in\C^+.$$
In particular $\sum_{j=1}^m c_j e^{i\phi_j}\neq 0$. 
\endproclaim

\demo{Proof} As $\b^0$ is ${\cal W}_\a$-invariant 
there is no loss of generality to assume that 
$\phi_1\leq \ldots\leq \phi_m$. Then 
$0\leq \phi_j-\phi_1< \pi$ for all $1\leq j\leq m$. Since $\sum_{j=1}^m \phi_j 
=0$ we have $\phi_m\geq 0$. 
Thus$\sum_{j=1}^m c_j e^{i\phi_j}$
is a sum of vectors not all zero in the real convex cone 
$$C =\{z\in \C\:  \phi_m-\pi <\arg(z)\leq \phi_m\}$$
in $\C$. 
In particular $\sum_{j=1}^m c_j e^{i\phi_j}$ is nonzero since the 
convex cone $C$ is pointed (i.e. contains no affine lines).  
Since $0\leq \phi_m< {m-1\over m}\pi$ (cf.\ Lemma 2.1(ii)) we also have 
$C\bs\{0\}\subeq \C^+$, concluding the proof of the lemma. \enddemo 

\proclaim{Proposition}   For $G=\Sl(m, \R)${\rm ,} 
$$K_\C B_\C^1 G\subeq K_\C A_\C N_\C.$$ 
\endproclaim

\demo{Proof} Take $a\in B_\C^1$ and recall that $a^2\in B_\C^0$. First we show that 
$a\oline N\subeq K_\C A_\C N_\C$. 
Let $\oline n\in \oline N$. Then Proposition 2.3(i) says that all 
principal minors of the complex symmetric matrix $\oline n^t a^2 \oline n$ are 
nonzero. 
Hence a theorem of Jacobi (cf.\ [Koe83, p.\ 124]) implies that 
there exist unique elements $b_0\in A_\C$ and $m\in N_\C$ such that 
$$\oline n^t a^2 \oline n=m^tb_0 m.$$  Let 
$a_0\in A_\C$ be such that $a_0^2=b_0$. Then we have 
$$a\oline n =k a_0 m$$
with $k\in K_\C$ given by $k=a\oline n m^{-1} a_0^{-1}$.

Using, as before, the Bruhat decomposition  
$G=\bigcup_{w\in {\cal W}_\a} \oline N wMAN,$
together with the $N_K(A)$-invariance of $B_\C^1$, we get that 
$aG\subeq K_\C A_\C N_\C$ for all $g\in G$, completing the proof. \enddemo

With $G=\Sl(m,\R)$ out of the way we want to use an observation that will allow 
us to obtain a proof of Theorem 1.8(ii) for appropriate 
subgroups.  The groups that will be covered in this way are: $\Sp(n,\R)$, 
$\Sp(p,q)$, $\Sp(n,\C)$,
$ \SU(p,q)$, $\SO^*(2n)$,  $\Sl(n, \C)$ and  $\Sl(n, \H)$.

Recall that a Levi subalgebra $\m$ of a standard parabolic subalgebra must be of 
the form $\m = \m (\Theta)$ for $\Theta \subeq\Pi$. If, moreover, $\m$ is 
$\theta$-stable, then the Iwasawa decomposition for $\m$ is compatible with that 
of $\g$. More generally, for $\g=\sL(m,\R)$ we consider  $\theta$-stable 
subalgebras $\g_1\subeq \g$ with a property that will give them Iwasawa 
decompositions compatible with that of $\g$. 
Set $\k_1=\k\cap\g_1$ and $\p_1 =\p\cap\g_1$ so that $\g_1=\k_1\oplus\p_1$ is 
a Cartan decomposition of~$\g$. Let $\a_1\subeq \p_1$ be a maximal abelian 
subspace. Since we can extend $\a_1$ to a maximal abelian subspace of $\p$ and 
since all maximal abelian subspaces of $\p$ are conjugate under $\Ad(K)$, we may 
assume that 
$\a_1\subeq \a$. Choose a positive system $\Sigma_1^+$ of $\Sigma_1=\Sigma(\g_1, 
\a_1)$. Then we can find 
a positive system $\Sigma^+$ of $\Sigma$ such that $\Sigma_1^+\subeq 
\Sigma^+\res_{\a_1}$. 
Write $\n_1 =\bigoplus_{\alpha\in\Sigma_1^+} \g_1^\alpha$
and note that $\n_1\subeq \n$. 

  We now impose the following condition on the restricted roots: 
$$\Sigma\res_{\a_1}\bs\{0\}=\Sigma_1.\leqno({\rm I})$$
It can be checked that (I) holds for example for the standard imbeddings of the 
subalgebras
$\g_1=\sp(n,\R)$ (with $2n=m$), $\su(p,q)$ (with $2p+2q=m$), $\sp(p,q)$ (with  
$2p+2q=m$) or 
$\so^*(2n)$ (with $2n=m$)(in all cases the fact that makes things work is that 
the restricted root system of $\g_1$ is either of type $C_n$ or $BC_n$). Further 
examples are $\g_1=\sL(n,\C)$ (with $2n=m$), $\sp(n,\C)$ (with $2n=m$)
or $\sL(n,\H)$ (with $4n=m$) (here the explanation is that the root system 
$\Sigma_1$ is 
of type $A$).    
Set 
$$\b_1^0 =\{X\in \a_1\: (\forall \alpha\in \Sigma_1)\ |\alpha(X)|<\pi\}$$
 and 
 $$\b_1^1 ={1\over 2}\b_1^0.$$ 
Then condition (I) guarantees that 
$$\b_1^1\subeq \b^1.\leqno(2.2)$$

  We denote by $G_1$ the analytic subgroup of  $G$  which is associated to 
$\g_1$. 
We assume that $G_1$ is closed. 
Further we denote by $K_1$, $A_1$, $N_1$ and $\oline N_1$ the 
analytic subgroups of $G_1$ corresponding to $\k_1$, $\a_1$, $\n_1$ and $\oline 
\n_1$. 
Finally we set $B_{1,\C}^1 =\exp(\a_1+i\b_1^1)$. 
In order to prove Theorem 1.8(ii) for the group $G_1$ we have to show that 

$$ B_{1,\C}^1 G_1\subeq K_{1,\C} A_{1,\C} N_{1,\C}$$
or equivalently 
$$(\forall b\in B_{1,\C}^1)(\forall g\in G_1) (\exists a\in A_{1,\C}, m\in 
N_{1,\C}),
\qquad g^t b g = m^t a m.\leqno(2.3)$$
In view of (2.2) and the validity of (2.3) for $G$ we deduce that for all 
$b\in B_{1,\C}^1$, $g\in G_1$ there exist  unique elements $m=m(b,g)\in N_\C$, 
$a=a(b,g)\in A_\C$
such that $g^t b^2 g = m^t a m$. Moreover $a=a(b,g)$ and $m=m(b,g)$ are 
analytic functions in the variables $b\in B_{1,\C}^1$, $g\in G_1$. Since we 
already 
know that $ a(A_{1,\C}^1,G_1)\subeq A_{1,\C}$ and  $m(A_{1,\C}^1,G_1)\subeq 
N_{1,\C}$ (cf.\ Proposition 1.7),
the analyticity of both $a$ and $m$ implies that   
$a(B_{1,\C}^1,G_1)\subeq A_{1,\C}$ and  $m(B_{1,\C}^1,G_1)\subeq N_{1,\C}$ 
proving \pagebreak (2.2). 

We summarize the above discussion with

 \proclaim{Proposition}  Assume that $G$ is one of the groups $\Sl(n,\R)${\rm ,} 
$\Sp(n,\R)${\rm ,} $\Sp(p,q)${\rm ,} $\SU(p,q)${\rm ,} 
$\SO^*(2n)${\rm ,} $\Sl(n,\C)${\rm ,} $\Sl(n,\H)$ or $\Sp(n,\C)$. Then 
$$K_\C B_\C^1 G\subeq K_\C A_\C N_\C.
$$
\endproclaim

There remain the restricted root systems for $\g=\so(p,q)$ and $\g=\so(n,\C)$. So 
first we recall some facts concerning these root systems of type $B_n$ and 
$D_n$. 
\vglue4pt
${\bf B}_{\bf n}$:   The root system $B_n$ is given by 
$$\Sigma=\{ \pm \eps_i\pm \eps_j\: 1\leq i\neq j\leq n\}\cup\{\pm\eps_i\: 
1\leq i\leq n\}.$$
A basis of $\Sigma$ is 
$$\Pi =\{\eps_1-\eps_2, \eps_2-\eps_3,\ldots, \eps_{n-1}-\eps_n, \eps_n\}.$$
If $\g$ is a split real Lie algebra with restricted root system $\Sigma$, then 
the $\omega_i$ are  the fundamental weights associated to 
$\Pi$ and given by 
$$\omega_1=\eps_1, \ \omega_2=\eps_1+\eps_2,\ \ldots,\   \omega_{n-1}=
\eps_1+\ldots+\eps_{n-1}, \ \omega_n={1\over 2}(\eps_1+\ldots+\eps_n).$$

\vglue4pt ${\bf D}_{\bf n}$:   The root system $D_n$ is given by 
$$\Sigma=\{ \pm \eps_i\pm \eps_j\: 1\leq i\neq j\leq n\}$$
and a basis of $\Sigma$ is given by 
$$\Pi=\{\eps_1-\eps_2, \eps_2-\eps_3,\ldots, \eps_{n-1}-\eps_n, 
\eps_{n-1}+\eps_n\}.$$
If $\g$ is a split real Lie algebra with restricted root system $\Sigma$, then 
the $\omega_i$ are the fundamental weights:
$$\omega_1=\eps_1, \ \omega_2=\eps_1+\eps_2,\ \ldots, \   \omega_{n-2}=
\eps_1+\ldots+\eps_{n-2}$$
and 
$$\omega_{n-1}={1\over 2}(\eps_1+\ldots+\eps_{n-1}-\eps_n), \ \omega_n={1\over 
2}
(\eps_1+\ldots+\eps_n).$$
\vglue4pt

To indicate the dependence of $\b^0$ on the root system, we shall write 
$\b^0(\Sigma)$ for $\b^0$. 

 \proclaim{Lemma} There exists $\b^0(D_n)=\b^0(B_n)$. 
\endproclaim

\demo{Proof} This is immediate from the equality $\conv(B_n)=\conv(D_n)$. \enddemo

The final goal of this section is to prove the inclusion 
$$B_\C^1 G\subeq K_\C A_\C N_\C\leqno(2.4)$$
for $G=\SO(p,q)$ or $G=\SO(n,\C)$. 

 We shall repeat the strategy used for the target group of type $A_n$. So 
assume for the moment that 
(2.4) holds for $G=\SO(n,n)$. Assume that $p\geq q$ and embed 
$G_1=\SO(p,q)$ into $\SO(p,p)$ in the natural way (upper left corner block). 
Then if we restrict $\Sigma$ to $\a_1$ we get a root system of type 
$B_q$ or $D_q$. Hence by our restriction procedure from the 
preceding section and Lemma~2.7,  we get (2.4) also for the subgroup $G_1$. 
Thus it suffices to prove (2.4) for $G=\SO(n,n)$ and $G=\SO(n,\C)$, with both
$\so(n,n)$ and $\so(n,\C)$ split.

 In what follows $\g$ denotes either $\so(n,n)$ or $\so(n,\C)$. We embed 
$\g$ into $\sL(2n,\R)$ 
as in the previous section. Then if we restrict the weights of $\sL(2n,\R)$ 
to $\g$, we obtain a root system of type $C_n$ or $BC_n$.  
We set 
$$\b_{\rm res}^0=\b^0(C_n)=\b^0(BC_n)$$
and 
$$\b_{\rm res}^1 ={1\over 2}\b_{\rm res}^0.$$ 
On the group side we define  $B_{{\rm res}, \C}^j =\exp(\a +i \b_{\rm res}^j)$ 
for $j=0,1$. 
In particular we get that 
$$B_{{\rm res}, \C}^1 G\subeq K_\C A_\C N_\C.\leqno(2.5)$$

We write $(\pi_n, V_n)$ for the $n^{\rm th}$ fundamental representation of $\tilde G$ 
with highest weight $\omega_n={1\over 2} (\eps_1+\ldots+\eps_n)$. 
We write ${\cal P}_n$ for the set of $\a$-weights of $(\pi_n, V_n)$
and set 
$$\b(\pi_n)^0=\{X\in\a\: (\forall \alpha\in {\cal P}_n)\  |\alpha(X)|<{\pi\over 
2}\}.$$
As usual we put $\b(\pi_n)^1={1\over 2}\b(\pi_n)^0$.

 \proclaim{Lemma} The following holds\/{\rm :}\/
$$\conv(\b_{\rm res}^0\cup \b(\pi_n)^0)\supeq \b^0.$$
\endproclaim

\demo{Proof} We claim that the extreme points of $\oline {\b^0}$ are given by 
$$\Ext(\oline{\b^0})=\left\{ \begin{array}{ll}\{ \pm \pi e_i, {\pi\over 2}
(\pm e_1\pm \ldots\pm e_n)\} & \hbox{for $n\geq 3$},\\ \{ \pm \pi e_i\} & \hbox{for $n=2$}.\end{array}\right.$$ 
In fact we have $\b^0=\b^0(B_n)$ by Lemma 2.7 and so $\Ext(\oline{\b^0})$ is 
invariant under the Weyl group ${\cal W}(B_n)=(\Z_2)^n\rtimes S_n$. {}From that 
the claim follows. 

 Now we have ${\pi\over 2}(\pm e_1\pm \ldots\pm e_n)\in \oline{\b_{\rm 
res}^0}$
and $\pm \pi e_i\in \oline {\b(\pi_n)^0}$. 
Hence the assertion of the lemma follows from the Krein-Milman theorem. \enddemo 

\proclaim{Proposition} Assume that $G=\SO(p,q)$ or $G=\SO(n,\C)$. Then 
for all $a\in B_\C^1${\rm ,}
$$K_\C B_\C^1 G\subeq K_\C A_\C N_\C.$$ 
\endproclaim

\demo{Proof} {}From what we have already done it is enough to prove the inclusion for 
$G=\SO(n,n)$ or $G=\SO(n,\C)$. By passing to a covering group if necessary we 
can also replace $G$ by $\tilde G$. 
Set $B_\C(\pi_n)=\exp(\a+i\b(\pi_n)^1)$. 
In view of Proposition 1.9, Lemma 2.7, (2.5) and Lemma 2.8 it remains to check 
that 
$$ B_\C(\pi_n)G\subeq K_\C A_\C N_\C.\leqno(2.6)$$
Suppose that (2.6) is false. Then we can find a $g\in G$ and a 
convergent sequence  $(a_j)_{j\in \N}$ in 
$B_\C(\pi_n)$ with $\lim_{j\to\infty} a_j=a_0\in B_\C(\pi_n)$, 
$\Theta(g)^{-1} a_j^2 g\in \oline N_\C A_\C N_\C$ for all 
$j\in \N$ but $\Theta(g)^{-1} a_0^2 g\not\in \oline N_\C A_\C N_\C$. 
In particular we find elements $m_j\in N_\C$ and $b_j\in A_\C$ with 
$$\Theta(g)^{-1} a_j^2 g=\Theta(m_j)^{-1} b_j m_j.$$
To arrive at a contradiction we have to show that $(b_j)_{j\in \N}$ is bounded 
(cf.\ Lemma 2.4).  
Let $\la\cdot, \cdot\ra$ denote an hermitian inner product on $V_n$ with 
$\la \pi_n(g)v, w\ra\break=\la v, \pi_n(\Theta(\oline g)^{-1})w\ra $ for all 
$g\in G_\C$, $v,w\in V$.  
Let $Q\subeq V_n\bs\{0\}$ be a compact subset. Then the definition of 
$\b(\pi_n)^0$ 
shows that 
$$\inf_{v\in Q} \Re \la \pi_n(\Theta(g)^{-1} a_j^2 g)v, v\ra >0.\leqno(2.7)$$

 If $v={\pi_n(m_j)^{-1}v_\alpha\over 
\|\pi_n(m_j)^{-1}v_\alpha\|}$ for a normalized weight vector $v_\alpha$
with weight $\alpha$, we get 
$$ {b_j^{\alpha}\la v_\alpha, \pi_n(\oline {m_j} m_j^{-1})v_\alpha\ra\over 
\|\pi_n(m_j^{-1})v_\alpha\|^2} =\la \pi_n(\Theta(m_j)^{-1} b_j m_j)v, v\ra 
=\la \pi_n(\Theta(g)^{-1} a_j^2 g)v, v\ra$$
for all $j\in \N$. 
In particular,  (2.7) implies that there are constants $C_1, C_2>0$ such that 
$$(\forall\alpha\in {\cal P}_n)\qquad C_1 > {|b_j^{\alpha}|\cdot |\la v_\alpha, 
\pi_n(\oline {m_j} m_j^{-1})v_\alpha\ra|\over 
\|\pi_n(m_j^{-1})v_\alpha\|^2}> C_2. $$

Recall that the weight spaces of the spin-representation $(\pi_n, V_n)$ are 
one-dimensional. Hence it follows that $\la \pi(n)v_\alpha, v_\alpha\ra=\la 
v_\alpha,v_\alpha\ra$ 
for all $n\in N_\C$ and all weight vectors $v_\alpha\in V_n$.  For the same 
reason we get 
$\|\pi_n(m_j)(v_\alpha)\|^2\geq 1$ for all $m_j$. In particular we obtain that 
$$(\forall\alpha\in {\cal P}_n)\qquad |b_j^{\alpha}|>C\leqno(2.8)$$
for some constant $C>0$. Now we have ${\cal P}_n=-{\cal P}_n$ and so (2.8) 
actually implies 
that $(b_j)_{j\in \N}$ is bounded. \enddemo

It seems reasonable to expect that a better technique would show  Theorem~1.8 to 
be 
valid also for the exceptional groups. Thus we formulate

\demo{Conjecture {\rm A}} Let $G$ be a semisimple Lie group with $G\subeq G_\C$ and $G_\C$ 
simply connected. Then 
\vglue12pt \hfill $B_\C^1 G\subeq K_\C A_\C N_\C .$ \hfill\qed
\enddemo

\numbereddemo{{R}emark} In [KrStII] we clarify the geometry of the domain, thereby giving 
more evidence for 
its  naturality. We show that the domain $K_\C \bs K_\C B_\C^1 G$ 
is bi-holomorphic 
to a maximal Grauert tube of $K\bs G$ having complex structure the adapted one. 
We also show the existence of a subdomain of $K_\C \bs K_\C B_\C^1 G$ 
bi-holomorphic to a Hermitian symmetric space but not isometric.  \enddemo

\vglue-9pt
\section{Holomorphic extension of irreducible
representations}
 
\vglue-6pt
We now come to our first application of the preceding construction, the 
holomorphic extension of representations. Additional applications of this will 
be given in subsequent sections for specific situations, such as principal 
series of representations, specific groups, or eigenfunctions on (locally) 
symmetric spaces.  The notation from representation theory needed for this 
section is standard and may be found explained in, say, [Kn86]. 

\vglue6pt {\it Notation.} As per Conjecture A we shall write $\Omega$ for 
$B_\C^1$ if $G$ 
is classical, and $A_\C^1$ otherwise.

\proclaim{Theorem}  Let $G$ be a linear{\rm ,} simple Lie group and let $(\pi, E)$ be 
an irreducible  Banach representation of $G$. 
Then for any $K$\/{\rm -}\/finite vector $v\in E_K${\rm ,} the orbit map
$$G\to E, \ g\mapsto \pi(g)v$$
extends to a $G$\/{\rm -}\/equivariant holomorphic map on $G\Omega K_\C$. 
\endproclaim

\demo{Proof} Set $V=E_K$, the collection of $K$-finite vectors of $(\pi, E)$. 

 Casselman's subrepresentation theorem (cf.\ [Wal88, 3.8]) gives the 
existence of a 
$(\g,K)$-embedding of $V$ into a principal series representation 
$$V\to (\Ind_{P_{\rm min}}^G(\sigma \otimes \lambda\otimes \1), {\cal 
H}_{\sigma, \lambda})\leqno(3.1)$$
where $P_{\rm min}=MAN$ is a minimal parabolic subgroup. In the next section we 
recall the standard terminology for principal series; in summary we 
set $\pi_{\sigma, \lambda}=\Ind_{P_{\rm min}}^G 
(\sigma\otimes\lambda\otimes\1)$; we write $(W_\sigma,\la\cdot,\cdot\ra_\sigma)$ 
for the representation Hilbert space of $\sigma$; we realize ${\cal 
H}_{\sigma,\lambda}$ as a Hilbert subspace of $L^2(K/M, W_\sigma)$, and we use 
induction from the right. 

 Write ${\cal H}$ for the completion of $V$ in ${\cal H}_{\sigma,\lambda}$. 
Let us first assume that $E={\cal H}$. 
Fix $v\in V$ and write $f_v\: G\to {\cal H}\subeq {\cal H}_{\sigma,\lambda}, \ 
g\mapsto \pi_{\sigma,\lambda}(g)v$ 
for the corresponding orbit map. 
Then we have for all $g\in G$ that  
$$(f_v(g))(kM)=a(g^{-1}k)^{\lambda-\rho} v(\kappa(g^{-1}k))\qquad (k\in 
K).\leqno(3.2)$$
Hence it follows from either Theorem 1.8 (for $\Omega=B_\C^1$) or Proposition 
1.7 
(for $\Omega=A_\C^1$) that analytic continuation of (3.2) gives rise to a map
\begin{eqnarray}
&&\tilde f_v\: G\Omega K_\C\to C^\infty(K/M,W_\sigma), \ \ 
\big(g\mapsto(kM\mapsto a(g^{-1}k)^{\lambda-\rho} 
v(\kappa(g^{-1}k)))\big).\speqnu{3.3}\\
\noalign{\noindent 
Note that $\tilde f_v\res_G=f_v$. }\nonumber
\end{eqnarray}

\vglue-19pt
 We claim that $\im\tilde f_v\subeq {\cal H}$. Write ${\cal H}^\bot$ for the 
orthogonal complement of ${\cal H}$ in the Hilbert space
$L^2(K/M, W_\sigma)$. Choose $w\in {\cal H}^\bot$. 
In order to show that $\la w,\im \tilde f_v\ra =\{0\}$, we may assume that $w$ 
is a $K$-finite, continuous function on $K/M$. Consider the function 
$$F\: G\Omega K_\C\to \C, \ \ g\mapsto \la \tilde f_v(g), w\ra=\int_K 
a(g^{-1}k)^{\lambda-\rho} \la v(\kappa(g^{-1}k)), 
w(k)\ra_\sigma \ dk,$$
with the equality on the right-hand side following from (3.3). 
Since $w$ is a bounded function, it is easy to see that $F$ is holomorphic. 
Since $F\res_G=0$ and $F$ is holomorphic we have 
$F=0$. This concludes  the proof of the claim. 

 Next we show that $\tilde f_v$ is holomorphic. 
Since $V\subeq C^{\infty}(K/M, W_\sigma)$ is dense in ${\cal H}$ 
and because weak holomorphicity implies holomorphicity, it is enough to show  
that for all $w\in V$ the analytically continued matrix coefficients
$$\pi_{v, w}\: G\Omega K_\C \to \C, \ \ g\mapsto \la \tilde f_v(g), w\ra  $$
are holomorphic. Again (3.3) gives that 
$$\pi_{v, w}(g)=\int_K a(g^{-1}k)^{\lambda-\rho} \la v(\kappa(g^{-1}k)), 
w(k)\ra_\sigma \ dk\qquad (g\in G\Omega K_\C)$$
and the holomorphicity of $\tilde f_v$ follows. 
\ssk Before we can deduce the general case from the case $E={\cal H}$ we need a 
little more refined information on the orbit maps. Note that $\im \tilde 
f_v\subeq {\cal H}^\infty$, 
${\cal H}^\infty$ the $G$-module of smooth vectors (indeed $\im \tilde f_v\subeq 
{\cal H}^\omega$, ${\cal H}^\omega$
the analytic vectors). Thus 
$\tilde f_v$ also induces a map $\hat f_v\: G\Omega K_\C\break\to {\cal H}^\infty$. 
Recall that the (Fr\'echet) topology on  ${\cal H}^\infty$ is induced from 
the seminorms
$${\cal H}^\infty\ni v\mapsto \|d\pi(u)v\|\qquad (u\in {\cal 
U}(\g_\C)).$$
{}From the explicit formula (3.2) of the induced action, one then deduces that 
$\hat f_v$ is continuous.
In particular $\hat f_v$ is 
holomorphic, since it is continuous and since for all $w$ in the 
dense subspace $V\subeq ({\cal H}^\infty)'$ the function 
$\la \hat f_v, w\ra=\pi_{v,w}$ is holomorphic.

 Finally we have to show how the general case follows from the case where 
$E={\cal H}$. We use the Casselman-Wallach globalization theorem (cf.\ [Wal92, 
11.6.7(2)]) which 
implies that the embedding (3.1) extends to a  $G$-equivariant topological 
embedding on the 
level of smooth vectors:
$$(\pi, E^\infty)\to (\pi_{\sigma,\lambda}, {\cal H}_{\sigma,\lambda}^\infty).$$
Hence the Fr\'echet representations $(\pi, E^\infty)$ and 
$(\pi_{\sigma,\lambda}, {\cal H}^\infty)$
are equivalent. As $\hat f_v$ was shown to be holomorphic for every 
$v\in V$, the proof of Theorem 3.1 is now \pagebreak complete.
 \enddemo

The holomorphic extension of the orbit map, $ g\mapsto \pi(g)v$, raises the 
question of the dependence of $\|\pi(g)v\|$ on $g$. This we address in a 
subsequent section. The holomorphic extension of a representation also gives 
rise to a holomorphic extension of its $K$-finite matrix coefficients. In the 
next section, we obtain estimates for the holomorphically extended matrix 
coefficients. 

\vglue-3pt
\section{Principal series representations}
 \vglue-6pt

{\it Integral formulas}.
We shall look in more detail at the growth properties of the holomorphic 
extension of matrix coefficients of principal series representations induced off 
a minimal parabolic subgroup. For now, we shall focus on the case of spherical 
principal series for two reasons: we shall use these results to obtain 
estimates on automorphic functions for locally symmetric spaces; the extension 
to the general case 
requires considering Eisenstein integrals and, albeit with many technicalities, 
given the holomorphic properties of the decompositions in Theorem 1.8, this 
presents no fundamentally new difficulties.

Set $\rho ={1\over 2} \sum_{\alpha\in \Sigma^+} m_\alpha \alpha\in \a^*$ with 
$m_\alpha =\dim \g^\alpha$. For $\lambda\in\a_\C^*$ 
we define a vector space
$${\cal D}_\lambda =\{ f\in C^\infty(G)\: 
(\forall man\in MAN)(\forall g\in G)\  f(gman)=a^{\lambda-\rho} f(g)\}.$$ 
The group $G$ acts on ${\cal D}_\lambda$ by left translation in the arguments, 
i.e., 
we obtain a representation $(\pi_\lambda, {\cal D}_\lambda)$ of $G$ 
given by $(\pi_\lambda(g)f)(x) =f(g^{-1}x)$ for $g,x\in G$, $f\in {\cal 
D}_\lambda$. Besides this realization we shall need the standard realizations of 
these representations that are called the compact (resp. noncompact) picture.    
The {\it compact realization} has for Hilbert space
$${\cal K}_\lambda=\oline {{\cal D}_\lambda\res_K}^{L^2(K)}\subeq L^2(K),$$
while the {\it noncompact realization} has
$${\cal N}_\lambda=\oline {{\cal D}_\lambda\res_{\oline N}}^{L^2(\oline N,\ 
a(\oline n)^{-2\Re (\lambda)}\ d\oline n)}\subeq 
L^2(\oline N,\ a(\oline n)^{-2\Re (\lambda)}\ d\oline n).$$
The representations $(\pi_\lambda, {\cal K}_\lambda)$ and $(\pi_\lambda, {\cal 
N}_\lambda)$ 
are continuous representations of~$G$. Moreover, the mapping 
$$\ f\res_K\to f\res_{\oline N} \qquad 
(f\in {\cal D}_\lambda)$$
extends to a unitary equivalence $(\pi_\lambda, {\cal K}_\lambda)\to 
(\pi_\lambda, 
{\cal N}_\lambda)$, provided $L^2(K)$ is obtained from a normalized Haar measure 
on $K$ 
and $L^2(\oline N)$ is obtained from a Haar measure $d\oline n$ which satisfies 
$\int_{\oline N} a(\oline n)^{-2\rho}\ d\oline n=1$. 
For $\lambda\in i\a^*$ the representations  
$(\pi_\lambda, {\cal K}_\lambda)$ and $(\pi_\lambda, {\cal N}_\lambda)$ are 
unitary. 
We will write $(\pi_\lambda, {\cal H}_\lambda)$ if we do not want to emphasize  
a particular realization.

We recall that for $(\pi, {\cal H})$ a continuous 
representation of a Lie group $G$ on some Hilbert space ${\cal H}$, a vector 
$v\in {\cal H}$ is called {\it analytic} if the orbit map 
$f_v\: G\to {\cal H}, \ g\mapsto \pi(g) v$ is analytic. Suppose that $G$ is 
contained in its universal complexification 
$G_\C$, and denote by $g\mapsto \oline g$ the complex conjugation 
in $G_\C$ with respect to the real form $G$.
Then for every analytic vector $v\in {\cal H}$ there 
exists a left $G$-invariant open neighborhood $U$ of $\1\in G_\C$ with $U=\oline 
U$ such that $f_v$ extends to a holomorphic map $\tilde f_v\: U\to {\cal H}, \ \ 
g\mapsto \pi(g) v$. With $\pi^*$ denoting the contragradient representation one 
has
$$\la \pi(g)v, v\ra =\la v, \pi(\oline g)^* v\ra\leqno(4.1)$$
for all $g\in U$.

For $G$ a Lie-group, $K<G$ a compact subgroup, 
and $(\pi, V)$ a continuous representation of $G$ on some topological vector 
space $V$, the representation $(\pi, V)$ is called {\it $K$-spherical} if 
$V^K\neq \{0\}$, $V^K=\{ v\in V\: (\forall k\in K) \pi(k) v=v\}$. 
For all $\lambda\in \a_\C^*$ the induced representation 
$(\pi_\lambda, {\cal D}_\lambda)$ is $K$-spherical, 
$\dim {\cal D}_\lambda^K=1$ and the function 
$$f_0\: G\to \C, \ \ x\mapsto a(x)^{\lambda-\rho}$$
is a generator of ${\cal D}_\lambda^K$. Moreover we have 
$$(\forall g,x\in G)\qquad (\pi_\lambda(g) f_0)(x)=a(g^{-1}x)^{\lambda-\rho}.$$
In the other realizations one has $v_0=f_0\res_K=\1_K\in {\cal K}_\lambda^K$, 
and $w_0=f_0\res_{\oline N}\in {\cal N}_\lambda^K$ given by 
$w_0(\oline n)=a(\oline n)^{\lambda-\rho}$.

 \proclaim{Proposition}   Let $(\pi_\lambda, {\cal K}_\lambda)$ be the 
compact realization of a spherical principal series representation with 
parameter $\lambda\in\a_\C^*$ and let $v_0=\1_K\in {\cal K}_\lambda^K$. 
Then the orbit map 
$$F\: G\to {\cal K}_\lambda, \ \ g\mapsto \pi_\lambda(g) v_0$$
extends to a holomorphic map 
$$\tilde F\: G\Omega K_\C\to {\cal K}_\lambda$$ 
on the open domain $G\Omega K_\C\subeq G_\C$. 
 \endproclaim

{\it Remark.} Note that a slight modification 
of Theorem 3.1 to representations of finite length  implies the proposition. 
But we shall give here a more direct proof avoiding the heavy machinery of 
representation theory. 

\demo{Proof} We consider the map 
$$\Phi\: G\times K\to \C,\ \  (g,k)\mapsto a(g^{-1}k)^{\lambda-\rho}$$
and note that $\Phi_\lambda(g, \cdot)=F(g)$. By Proposition 1.7 and Theorem 1.8 
the function 
$\Phi$ extends to an analytic map 
$$\tilde \Phi\: G \Omega K_\C\times K\to \C,\ \  (z,k)\mapsto 
a(z^{-1}k)^{\lambda-\rho}$$
which is holomorphic in the first argument.  
It is obvious that $\tilde \Phi(z,\cdot)\in L^2(K)$ for all $z\in G\Omega K_\C$. 
Let 
$P\:L^2(K)\to {\cal K}_\lambda$ denote the orthogonal projection and define 
$$\tilde F\:   G\Omega K_\C\to {\cal K}_\lambda, \ \ z\mapsto P(\tilde 
\Phi(z,\cdot)).$$
Then $\tilde F\res_G=F$ and it remains to show that $\tilde F$ is holomorphic. 
For that however it suffices to show that 
$$ G\Omega K_\C\to \C,\ \  z\mapsto \la \tilde F(z), f\ra$$
is holomorphic for all $f\in {\cal D}_\lambda\res_K\subeq C^\infty (K)$. But 
this in 
turn follows from 
$$\la \tilde F(z), f\ra= \int_K a(z^{-1} k)^{\lambda- \rho} \oline{f(k)} \ dk$$
by the compactness of $K$, the continuity of $\tilde \Phi$ and the 
holomorphy of $\tilde\Phi(\cdot, k)$.  \enddemo

If $\lambda\in \a_\C^*$ and $(\pi_\lambda, {\cal K}_\lambda)$ is
the induced  representation realized in the compact picture, then the matrix 
coefficient of the $K$-fixed vector with itself is the familiar zonal spherical 
function,
$$\phi_\lambda(g)=\la \pi_\lambda(g^{-1}) v_0, v_0\ra.\leqno(4.2)$$

The holomorphic extension to $G \Omega K_\C$ of $\pi_\lambda(g^{-1})v_0$ gives a 
holomorphic extension of the matrix coefficient $\phi_\lambda(g)$. However, this 
is not the largest domain of analyticity for $\phi_\lambda(g)$. Since we will 
estimate the norm of $\pi_\lambda(g^{-1}) v_0$ by means of $\phi_\lambda(g)$, in 
order to obtain optimal estimates on the norm it will be important to have an 
expression that represents $\phi_\lambda(g)$ in its entire domain of holomorphy. 
In terms of the pairing in the compact realization, $\phi_\lambda(g)$ is given 
by the well-known integral formula 
$$\phi_\lambda(g) = \int_K a(gk)^{\lambda-\rho}\ dk.$$

By the $K$-bi-invariance of $\phi_\lambda$ and in light of Proposition 1.7, this 
defining integral formula for the spherical function can be extended to $K_\C 
\Omega K_\C$. But in general the integral formula need not extend
to any larger domain (cf. Example 4.3). There are a couple of reasons for this. 
First, the integrand $k\mapsto a(a^{-1}k)^{\lambda-\rho}$
becomes singular if  $a$ leaves $A_\C^1$, and secondly, it is no longer possible 
to take holomorphic square roots
(the $\rho$-exponent frequently involves a square root). We shall present an 
alternative integral formula valid on a domain about twice as large and this 
will be crucial for the estimates on the norm of $\pi_\lambda(g^{-1}) v_0$.

To state the result we recall the notation $\Omega$, viz.  if $G$ is classical, 
then $\Omega=B_\C^1$ and otherwise $\Omega=A_\C^1$. Consistent with this and the 
notation $B_\C^0$ (resp. $A_\C^0$), we use $\Omega^2=B_\C^0$ if $G$ is classical 
and otherwise set $\Omega^2=A_\C^0$.

\proclaim{Theorem} Let $\lambda\in \a_\C^*$ and $\phi_\lambda$ be the 
spherical function with parameter $\lambda$ associated to 
$G/K$.    
\begin{itemize}
\ritem{(i)} The spherical function $\phi_\lambda$ extends to a $K_\C$\/{\rm -}\/bi\/{\rm -}\/invariant 
function on\break
$K_\C \Omega^2K_\C\subeq G_\C$ which is holomorphic when restricted to 
$\Omega^2$. 
\ritem{(ii)}  $(\forall b\in A)(\forall a\in \exp(i\a)\cap\Omega)${\rm ,}
$$\phi_\lambda(ba^2)=\int_K a(bak)^{\lambda -\rho}\cdot  
\oline{a(ak)^{\lambda -\rho}} \cdot
 \oline{a(ak)^{-2\Re\lambda}} \ dk.$$
In particular{\rm ,} for $\lambda\in i\a^*$ we get for all  $a\in \exp(i\a)\cap\Omega$
$$\phi_\lambda(a^2)=\int_K |a(ak)^{2(\lambda -\rho)}| \ dk.$$
\ritem{(iii)} $(\forall b\in A)(\forall a\in \exp(i\a)\cap\Omega)${\rm ,}
$$\phi_\lambda(ba^2)=\int_{\oline N} a(ba\oline n)^{\lambda -\rho}\cdot  
\oline{a(a\oline n)^{\lambda -\rho}} \cdot
 \oline{a(a\oline n)^{-2\Re\lambda}} \ d\oline n .$$
In particular{\rm ,} for $\lambda\in i\a^*${\rm ,} for all $a\in \exp(i\a)\cap\Omega,$ 
$$\phi_\lambda(a^2)=\int_{\oline N}|a(a\oline n)^{2(\lambda -\rho)}| \ 
d\oline n.$$
\end{itemize}

\endproclaim

\demo{Proof} (i)  It suffices to show that $\phi_\lambda\res_A$ extends to 
a holomorphic function on $\Omega^2$. 
We will work with the compact realization $(\pi_\lambda, {\cal K}_\lambda)$. Let 
$a\in A$. Then (2.1) implies that 
$$(\forall a\in A) \qquad \phi_\lambda(a^2)=\la \pi_\lambda(a^{-1}) v_0, 
\pi_\lambda(a^{-1})^* v_0\ra. \leqno(4.3)$$
We now analytically continue the right-hand side of (4.3). 
Recall from [Kn86, p.\ 170] that for $f\in {\cal D}_\lambda$ and $x,g\in G$ one 
has 
$$(\pi_\lambda(g)^* f)(x)=a(gx)^{-2\Re \lambda} f(gx).\leqno(4.4)$$
Hence $\pi_\lambda(g)^* f=a(g\cdot)^{-2\Re \lambda} \pi_\lambda(g^{-1}) f$
for all $g\in G$. Similarly as in Proposition~4.1 one shows that $g\mapsto 
\pi(g)^* v_0$
extends to a holomorphic ${\cal K}_\lambda$-valued map on $G\Omega K_\C$. 
Thus Proposition 4.1 implies that 
the function 
$$A\to \C,\ \  a\mapsto \la \pi_\lambda(a^{-1}) v_0, 
 \pi_\lambda(a^{-1})^* v_0\ra$$
extends to a holomorphic function on $\Omega$. Since 
we have a unique holomorphic square root on $\Omega^2$, 
namely 
$$\Omega^2\to \Omega, \ a=\exp(X)\mapsto\sqrt{a} =\exp({1\over 2} X),$$
the assertion of (i) now follows from (4.3).
\vglue4pt
(ii) In view of the proof of (i), (ii) is immediate from the analytic 
extensions of (4.3) and (4.4) to $\exp(i\a)\cap\Omega$.  
\vglue4pt
(iii) This is proved as (ii) is by use of the noncompact realization 
$(\pi_\lambda, {\cal N}_\lambda)$ instead of $(\pi_\lambda, {\cal 
K}_\lambda)$. \enddemo

\numbereddemo{{E}xample} We explicate the theorem for the group 
$G=\Sl(2,\R)$. Clearly $G\subeq G_\C=\Sl(2,\C)$ and $G_\C$ 
is simply connected. We let $\k=\so(2)$,  
$$\a =\{ \left(\begin{array}{cc} x & 0\\ 0 & -x\end{array}\right)\: x\in \R\}\quad \hbox{and}\quad \n =\{ 
\left( \begin{array}{cccc}0 & n\\ 0 & 0\end{array}\right)\: n\in \R\}.$$ 
For $z\in \C^*$, $x\in \C$ and  $\theta\in\C$  we set 
$$a_z =\left( \begin{array}{cccc}z & 0\\ 0 & z^{-1}\end{array}\right)\in A_\C, \quad 
n_x =\left( \begin{array}{cccc} 1& x\\ 0 & 1\end{array}\right)\in N_\C $$
and
$$
k_\theta =\left( \begin{array}{cccc}\cos \theta & \sin \theta\\ -\sin \theta & \cos\theta\end{array}\right)\in 
K_\C.$$
Then 
$$A_\C^0=B_\C^0=\{ a_z\: \Re(z)>0\}\quad\hbox{and}\quad A_\C^1=\{ 
a_z\:|\arg(z)|<{\pi\over 4}\}.$$
Since $\omega =\omega_1=\rho$, we may identify $\a_\C^*$ with $\C$ by means of 
the isomorphism
$$\R\to\a^*, \ \ \lambda\mapsto \lambda\omega.$$

Let us consider the spherical function with parameter $\lambda\in i\a^*$. Then 
Proposition A.1(i)
in the appendix shows that (4.2), the defining integral formula for 
$\phi_\lambda$, 
extends to $A_\C^1$ and we have 

$$(\forall a_z\in A_\C^1)\qquad \phi_\lambda(a_z)={1\over 2\pi}\int_0^{2\pi} 
{d\theta\over (z^2+
\sin^2\theta({1\over z^2} -z^2))^{{1\over 2}-\lambda}}.\leqno (4.5)$$
 
It is easy to see that it is no longer possible to take consistently 
an analytic square root of $(z,\theta)\mapsto  z^2+
\sin^2\theta({1\over z^2} -z^2)$ if $|\arg(z)|$ becomes larger than ${\pi\over 
4}$. 
Also as this function has zeros, the integrand of the integral expression 
above becomes singular (although in this case the singularity is integrable). 
On the other hand Theorem 4.2(ii), (iii), together with Proposition A.1 imply, 
for 
all $ -{\pi\over 2}< \phi<{\pi\over 2}$, $r>0$, that 
\begin{eqnarray*}
&&\hskip-.25in\phi_\lambda(a_{re^{i\phi}})\\[3pt]
&&\nhs={1\over 2\pi}\int_{-\pi}^\pi
{d\theta\over (r^2e^{i\phi} +\sin^2\theta({1\over r^2} e^{-i\phi} -r^2 
e^{i\phi}))^{{1\over 2}-\lambda}
(e^{-i\phi} +\sin^2\theta(e^{i\phi} -e^{-i\phi}))^{{1\over 2}+\lambda}}\end{eqnarray*}
and 
$$\phi_\lambda(a_{re^{i\phi}})=\int_{-\infty}^\infty
{d x\over (r^2e^{i\phi}+ {1\over r^2} e^{-i\phi} x^2)^{{1\over 2}-\lambda}
 (e^{-i\phi} + e^{i\phi} x^2)^{{1\over 2}+\lambda}}.$$
If one examines the integral over $K$, one sees that the second factor is 
identically 1 when evaluated on the real group $A$, so that it comes into play 
only on the complex domain.  \enddemo

{\it An upper estimate}.
We can give a soft upper estimate along the convex hull of extreme points of the 
domain.

 \proclaim{Proposition}  Let $(\pi_\lambda, {\cal H}_\lambda)$ be a unitary principal 
series representation of $G$. Let $\omega$ be $\b^0$ if $G$ is classical and 
$\a^0$ otherwise. 
If $a,b\in \omega${\rm ,} then 
$$\sup_{0\leq t\leq 1} \phi_\lambda\Big(\exp(i(ta+(1-t)b))\Big)\leq 
\max\{\phi_\lambda(\exp(ia)), 
\phi_\lambda(\exp(ib))\}.$$
\endproclaim

\demo{Proof} Set $S_{[0,1]}=\{ z\in \C\: \Re z\in [0,1]\}$ and $X=\{za+(1-z)b\: 
z\in S_{[0,1]}\}$. Then 
$$S_{[0,1]}\to X, \ \ z\mapsto g(z)=za+(1-z)b$$
defines a bi-holomorphism of complex manifolds with boundary. We set 
$$f\: S_{[0,1]}\to\C, \ \ z\mapsto \phi_\lambda(\exp(ig(z))).$$
Then $f$ is holomorphic on $\Int  S_{[0,1]}$ and we claim that $f$ is bounded. 
In fact, we have 
$$f(z)=\la \pi_\lambda(\exp(\Im g(z)))\pi_\lambda(\exp(-i{1\over 2}\Re 
g(z)))v_0,
\pi_\lambda(\exp(-i{1\over 2}\Re g(z)))v_0\ra, $$
and so by the unitarity of $\pi_\lambda$
$$|f(z)|\leq \la \pi_\lambda(\exp(-i{1\over 2}\Re g(z)))v_0,
\pi_\lambda(\exp(-i{1\over 2}\Re g(z)))v_0\ra=f(\Re z).$$
This implies our claim and so the assertion of the proposition follows from the 
Phragmen-Lindel\"of principle.  \enddemo

{\it A radial lower estimate}.
A precise estimate of the nature of the singularity along the entire boundary of 
$\Omega$  appears difficult. However, for an approach to the boundary along the 
direction of roots (or co-roots), we can obtain estimates. In this regard see 
Remark 5.5. We recall our 
standing hypothesis that $G$ is a semisimple Lie group contained in its 
complexification $G_\C$.

Let $H\in \a$ and assume that 
there is a $\theta$-invariant  $\sL(2,\R)$-triple\break $\{ H, X, \theta(X)\}\subeq \g 
$, i.e., 
$[H,X]=2X, \ [H, \theta(X)]=-2\theta(X), \ [X, \theta(X)]= H$.

 We shall want to give estimates for the radial behaviour  
$$\phi_\lambda(\exp(i{\pi\over 2} (1-\eps)H))$$ for $\eps\to 0$ provided that 
${\pi\over 2}H\in \partial \b^0$ and $(\pi_\lambda, {\cal H}_\lambda)$ is 
unitarizable. We do this by restriction of the representation to a subgroup 
isomorphic to $\Sl(2,\R)$ or $\PSl(2,\R)$. 
The triple gives 
$$\g_0 =\span_\R \{H, X, \theta (X)\},$$
a $\theta$-stable subalgebra of $\g$ isomorphic to $\sL(2,\R)$. Its 
Cartan \pagebreak decomposition is given by $\g_0=\k_0\oplus\p_0$ with 
$\k_0=\k\cap \g_0$ and $\p_0=\p\cap \g_0$. Further, we set 
$\a_0=\R H$. We denote the analytic subgroups of $G$ 
corresponding to $\g_0$, $\a_0$ and $\k_0$ by 
$G_0$, $A_0$ and $K_0$.

 We write $\hat G_0$ for the set of equivalence classes of 
unitary irreducible representations of $G_0$ (the {\it unitary dual of $G_0$}). 
Since $G_0$ is semisimple and thus of type I, 
there is a natural $T_1$-topology on $\hat G_0$,  the {\it hull-kernel 
topology}, 
which we denote by $\tau_{\rm hk}$,   
and a Borel measure $\mu$ on  $\hat G_0$ such that 
$$(\pi_\lambda\res_{G_0}, {\cal H}_\lambda)=
\Big(\int_{\hat G_0}^\oplus \pi_\sigma\otimes I 
\ d \mu(\sigma), \int_{\hat G_0}^\oplus {\cal H}_\sigma\hat \otimes V_\sigma 
\ d \mu(\sigma)\Big).$$
Here $(\pi_\sigma, {\cal H}_\sigma)$ denotes a representative of $\sigma\in \hat 
G_0$. In particular, $v_0\in {\cal H}_\lambda^K$ disintegrates 
as 
$$v_0=\int_{\hat G_0} v_0^\sigma \ d\mu(\sigma)$$
with $1=\|v_0\|^2=\int_{\hat G_0} \|v_0^\sigma\|^2 \ d\mu(\sigma)$. 
As each $v_0^\sigma$ is $K_0$-fixed, for all $\sigma$ with $v_0^\sigma\neq 0$,  
$$(\forall g\in G_0)\qquad 
\phi_\sigma^0(g)= {1\over \|v_0^\sigma\|^2}\la \pi_\sigma(g^{-1})v_0^\sigma, 
v_0^\sigma\ra$$
defines  a spherical function on $G_0$. 

In particular we get for all $a\in A_0$
$$\phi_\lambda(a)=\la\pi_\lambda(a^{-1})v_0,v_0\ra =\int_{\hat G_0} 
\la\pi_\sigma(a^{-1})v_0^\sigma, v_0^\sigma\ra \ d\mu(\sigma) =\int_{\hat G_0} 
\phi_\sigma^\alpha(a)\|v_0^\sigma\|^2 d\mu(\sigma). \leqno(4.6)$$

 \proclaim{Proposition}  Let $G$ be a semisimple Lie group with Lie algebra $\g$ and 
assume that $G\subeq G_\C$. Suppose that $\{H, X, \theta(X)\}$ with $H\in\a$ 
forms an 
$\sL(2,\R)$\/{\rm -}\/triple in $\g$. Assume that ${\pi\over 2}H\in \partial \b^0$ and 
that 
$(\pi_\lambda, {\cal H}_\lambda)$ is unitarizable.  
Then we have 
$$\phi_\lambda(\exp(i{\pi\over 2}(1-\eps) H))\geq C|\log \eps| $$
for $0<\eps\leq 1$ and a constant $C>0$.

\demo{Proof} Write $Y_0=\{ \sigma\in \hat G_0\: \pi_\sigma \ \hbox{is 
$K$-spherical}\}$. {}From the well-known details of the unitary 
dual of $G_0$ we know  that 
there is  a natural parametrization $i\R\,\cup\,]0,1[\,\to Y_0$. 
Moreover if we equip $i\R\,\cup\, ]0,1[$ with its Euclidean topology, then this 
parametrization becomes continuous (this essentially follows from the fact that 
the assignments 
$\sigma\mapsto \phi_\sigma^0(a)$ are continuous with respect to the Euclidean 
topology; cf.\ the technique of [Wal92, 14.12.3]). 
In particular this parametrization induces the (possibly) stronger Euclidean 
topology $\tau_e$ on 
$Y_0$. Hence if $\bigcup_{n\in\N} Q_n$ is an exhaustion of 
$(Y_0,\tau_e)$ by compact sets,  then $\bigcup_{n\in\N} Q_n$ defines an 
exhaustion 
of quasicompact Borel sets of $(Y_0,\tau_{\rm hk})$. So  we can find 
a $Q_n$ with $\int_{Q_n} v_0^\sigma\ d \mu(\sigma)\neq 0$. 
It follows from (4.6) that we have 
$$\phi_\lambda(\exp(i{\pi\over 2}(1-\eps) H))=\int_{Y_0} 
\phi_\sigma^0 (\exp(i{\pi\over 2}(1-\eps) H)) \|v_0^\sigma\|^2  \ d\mu(\sigma)$$
for $\eps>0$. In particular we get that  
$$\phi_\lambda(\exp(i{\pi\over 2}(1-\eps) H))\geq \int_{Q_n} 
\phi_\sigma^0 (\exp(i{\pi\over 2}(1-\eps) H)) \|v_0^\sigma\|^2  \ d\mu(\sigma)$$
for $\eps>0$. Since $Q_n$ is compact, Theorem 5.1 (to follow) implies that 
there is a constant $C\geq 0$ such that $\phi_\sigma^0 (\exp(i{\pi\over 
2}(1-\eps) H))\geq C|\log \eps|$
for all $\sigma\in Q_n$. 
Hence we get that 
$$\phi_\lambda(\exp(i{\pi\over 2}(1-\eps) H))\geq C |\log\eps| \int_{Q_n} 
\|v_0^\sigma\|^2  \ d\mu(\sigma), $$
proving the theorem. \enddemo

Set $\Sigma_0=\{\alpha\in \Sigma\: 2\alpha\not\in \Sigma\}$. 
Let $\alpha\in \Sigma_0$ and $H_\alpha\in \a$ be the corresponding {\it co\/{\rm -}\/root}, 
i.e., 
$H_\alpha\in [\g^\alpha, \g^{-\alpha}]\cap \a$ such that $\alpha(H_\alpha)=2$.

\proclaim{{C}orollary}  Suppose that $G$ is one of the groups $\Sl(n,\R)${\rm ,} $\Sl(n,\C)${\rm ,} 
$\Sl(n,\H)${\rm ,} $\Sp(n,\R)${\rm ,} $\SO^*(2n)$ or $\SU(p,q)$. Let $\alpha\in \Sigma_0$ 
and $H_\alpha$ be its 
co\/{\rm -}\/root.\break 
Assume that $(\pi_\lambda, {\cal H}_\lambda)$ is unitarizable. 
Then there exists a constant $C${\rm ,} depending only on $\lambda${\rm ,} such that 
$$\phi_\lambda(\exp(i{\pi\over 2}(1-\eps) H_\alpha))\geq C|\log \eps| $$
for $0<\eps\leq 1$.
\endproclaim

\demo{Proof} Since all 
restricted root systems are either of type $A_n$, $C_n$ or $BC_n$ we have 
${\pi\over 2}H_\alpha\in \partial \b^0$. 
Then the assertion follows from Proposition 4.5. \enddemo

\section{Real rank one}

{\it Singularity of spherical functions}.
We consider Lie algebras of real rank one. For these  we are able to obtain 
sharp  asymptotic behaviour of the holomorphically 
extended spherical functions. 

As $\g$ has real rank one, $\dim \a=1$. As is the custom, we set $p =\dim 
\g^\alpha$, $q =\dim \g^{2\alpha}$, and $c ={1\over 4(p+2q)}$. Here 
$\oline \n=\g^{-\alpha}\oplus\g^{-2\alpha}$ and $[\g^{-\alpha}, 
\g^{-2\alpha}]=\{0\}$.

We have the familiar formula for all $\oline n=\exp(X+Y)=\exp(X)\exp(Y)$, $X\in  
\g^{-\alpha}$, $Y\in \g^{-2\alpha}$:
$$a(\oline n)^\rho=[(1+ c\|X\|^2)^2+ 4c\|Y\|^2]^{{p+2q\over 4}}.\leqno(5.1)$$ 
Here $\|Z\|^2 =-\kappa(Z,\theta Z)$ for all $Z\in\g$ with $\kappa$ denoting the 
Cartan-Killing form of $\g$. 
Computations involving 
$a(\oline n)^\rho$ have appeared many times. We include the following 
computations only because they  involve the holomorphic extension of $a(\oline 
n)^\rho$ and 
for this we have no convenient reference. 

 Let $A_\alpha\in\a$ be defined by $\alpha(A_\alpha)=1$. For convenience we 
shall identify 
$\a_\C$ and $\a_\C^*$ with $\C$ by means of the isomorphisms
\begin{eqnarray*}
\C\to\a_\C, && z\mapsto z A_\alpha,\\
 \C\to\a_\C^*, && \lambda\mapsto \lambda\alpha.
\end{eqnarray*}

Thus, $\phi_\lambda(e^z) := \phi_\lambda$( exp $ zA_\alpha)$.

Let $\Omega_\g\subeq \a$ denote $\b^0$ if $\g$ is classical and $\a^0$ 
otherwise.

Here and henceforth we   use the notation  $f(\eps)\asymp g(\eps)$ for 
two positive valued 
functions $f(\eps)$, $g(\eps)$ if there exist constants $c_1, c_2>0$ such that 
$c_1 f(\eps)\leq g(\eps)\leq c_2 f(\eps)$ for all $\eps$. 
 
\proclaim{Theorem} Let $G$ be a connected Lie group of real rank one
contained in its universal complexification $G_\C$.  
\begin{itemize}
\ritem{(i)} For all $\lambda\in i\a^*$ the maximal tube domain of 
definition of $\phi_\lambda\circ \exp_A$ is given by 
$$T_{\lambda, {\rm max}}=\b_\C^0.$$

\ritem{(ii)} For $X\in \partial \b^0${\rm ,} and a fixed $\lambda\in \a_\C^*${\rm ,} there 
exist a $C>0$ such that  for 
$\eps\to 0$, $\eps>0$
$$|\phi_\lambda(\exp(\pm i(1-\eps)X)|\leq C \left\{ \begin{array}{ll}|\log \eps| & \hbox{for 
$p=1${\rm ,} $q=0$},\\ \eps^{-p+1} & \hbox{for $p>1${\rm ,} $q=0$},\\ |\log \eps| & \hbox{for $q=1$},\\ \eps^{-q+1} & \hbox{for $q>1$}.
\end{array}\right.$$  If in addition $\lambda\in i\a^*${\rm ,} then 
$$|\phi_\lambda(\exp(\pm i(1-\eps)X)|\asymp\left\{ \begin{array}{ll}|\log \eps| &\hbox{for  
$p=1${\rm ,} $q=0$},\\ \eps^{-p+1} & \hbox{for $p>1${\rm ,} $q=0$},\\ |\log \eps| & \hbox{for $q=1$},\\ \eps^{-q+1} & \hbox{for
$q>1$}.\end{array}\right.$$ 
\end{itemize}

 \endproclaim

\numbereddemo{{R}emark}  For some $\lambda\in \a_\C^*\bs i\a^*$ it can happen that 
$T_{\lambda, {\rm max}}= \a_\C$; i.e., the spherical function 
$\phi_\lambda\res_A$ extends 
holomorphically to $A_\C$. Simply consider $G=\Sl(2,\C)$. 
Then $\g\cong\so(3,1)$, which in our previous notation corresponds to the case 
$p=2$ and $q=0$. The explicit formula for spherical functions on complex groups 
specialized to $\Sl(2,\C)$ then reads 
$$\phi_\lambda(e^z)={1\over \lambda} {e^{(\lambda)z} 
-e^{-(\lambda)z}\over e^z -e^{-z}}.$$
Hence we see that $T_{\lambda, {\rm max}}= \b_\C^0$ for $\lambda\not\in \Z$ 
while for 
$\lambda\in \Z$ one has $T_{\lambda, {\rm max}}= \a_\C$.  \enddemo

The proof of Theorem 5.1 is computational and is presented in lemmas for the 
various cases.

With our parametrization it follows from (5.1) and Theorem 4.2(iii) that 
for all $\lambda\in i\a^*$, $\phi\in \Omega_\g$,
\begin{eqnarray*}
&&\hskip-.25in \phi_\lambda(e^{-i\phi})=e^{-i\lambda\phi} \\[4pt] &&
\times\int_{\R^p}\int_{\R^q}   {\textstyle{dX \ dY\over 
[(1+c e^{i\phi} \|X\|^2)^2 +4 c e^{2i\phi} \|Y\|^2]^{{p+2q\over 4}-{\lambda\over 
2}}
[(1+c e^{-i\phi}\|X\|^2)^2 +4 c e^{-2i\phi} \|Y\|^2]^{{p+2q\over 
4}+{\lambda\over 2}}}}.\end{eqnarray*} 
Using polar coordinates  we thus obtain 
\begin{eqnarray*}
&&\hskip-.25in \phi_\lambda(e^{-i\phi})=\tilde c e^{-i\lambda\phi}\\
&&\times \int_0^\infty \int_0^\infty 
{u^{p-1} v^{q-1}du \ dv\over 
[(1+e^{i\phi} u^2)^2 +e^{2i\phi} v^2]^{{p+2q\over 4}-{\lambda\over 2}}
[(1+e^{-i\phi} u^2)^2 +e^{-2i\phi} v^2]^{{p+2q\over 4}+{\lambda\over 2}}}\end{eqnarray*}
for a constant $\tilde c$  depending on  only $p$ and $q$. 
Finally with the substitution $r=u^2$, $s=v^2$ we arrive at 
\begin{eqnarray}
&& \phi_\lambda(e^{-i\phi})=Ce^{-i\lambda\phi}\speqnu{5.2}\\
&&\quad\times \int_0^\infty \int_0^\infty 
{r^{{p-2\over 2}} s^{{q-2\over 2}}dr \ ds\over 
[(1+e^{i\phi} r)^2 +e^{2i\phi} s]^{{p+2q\over 4}-{\lambda\over 2}}
[(1+e^{-i\phi} r)^2 +e^{-2i\phi} s]^{{p+2q\over 4}+{\lambda\over 
2}}} \nonumber
\end{eqnarray}
for all $\lambda\in i\a^*$, $\phi\in \Omega_\g$ and a constant 
$C$ which is independent of $\lambda$. 
We distinguish three cases.

\demo{Case 1: ${p=1}$, ${q=0}$}  In this case we have $G=\Sl(2,\R)$, 
the root system is split 
(i.e. $\Delta=\Sigma$), and so $\omega =\omega_1={1\over 2}\alpha$. Hence 
$$\a_\C^0=\b_\C^0=\{ z\in\C\: |\Im z|<\pi\}$$ and 
(5.2) boils down to 
$(\forall \lambda\in i\R)(\forall -\pi<\phi<\pi)$,
$$ \phi_\lambda(e^{-i\phi})=C e^{-i\lambda\phi}\int_0^\infty {dr \over 
 \sqrt{r} (1+e^{i\phi}r)^{{1\over 2}-\lambda} (1+e^{-i\phi}r)^{{1\over 
2}+\lambda}}.\leqno(5.3)$$
\enddemo
 
 \proclaim{Lemma}  For $p=1${\rm ,} $q=0$ and $\lambda\in i\a^*${\rm ,}
$T_{\lambda, {\rm max}}=\b_\C^0=\a_\C^0$.  
Moreover{\rm ,}  for a fixed $\lambda$ the asymptotics at the boundary are given by 
$$\phi_\lambda(e^{-i(\pi - \eps)}) \asymp |\log\eps|$$
for $\eps\to 0${\rm ,} $\eps>0$.
\endproclaim

\demo{Proof} This is immediate from (5.3). \enddemo

\demo{Case 2: ${p>1}$, ${q=0}$} Here $\g=\so(p+1, 1)$ is 
a classical Lie algebra and so $\Omega_\g=\b^0$. We have $\omega 
=\omega_1=\alpha$ and so 
$\a_\C^0=\{ z\in\C\: |\Im z|<{\pi\over 2}\}$ and $\b_\C^0=2\a_\C^0$. Formula  
(5.2) simplifies to 
$$\phi_\lambda(e^{-i\phi})=e^{-i\lambda\phi}\int_0^\infty {r^{{p-2\over 2}} 
dr\over 
(1+e^{i\phi} r)^{{p\over 2}-\lambda}(1+e^{-i\phi} r)^{{p\over 
2}+\lambda}} \leqno(5.4)$$
for all $\lambda\in i\R$ and $ -\pi <\phi<\pi$. 
\enddemo

 \proclaim{Lemma}  For $p>1${\rm ,} $q=0$ and for all $\lambda\in i\a^*${\rm ,}   
$T_{\lambda,{\rm max}}=\b_\C^0=2\a_\C^0=\{ z\in\C\: |\Im z|<\pi\}$. 
Moreover{\rm ,} for a fixed $\lambda$ the asymptotics at the boundary are given by 
$$\phi_\lambda(e^{-i(\pi - \eps)}) \asymp {1\over \eps^{p-1}}$$
for $\eps\to 0$, $\eps>0$. 
\endproclaim

\demo{Proof} We will estimate $\phi_\lambda(e^{-i(\pi - \eps)})$ 
for $\eps\to 0$   $(\eps>0)$. Since the unbounded 
contribution to the integral is local (at $r=1$) we may henceforth assume 
that $\lambda=0$. Then (5.4) gives 
\begin{eqnarray*}
\phi_0(e^{-i(\pi-\eps)})&=&\int_0^\infty {r^{{p-2\over 2}} dr\over 
|(1+e^{i(\pi-\eps)} r)|^p}\\[5pt] &\asymp&  \int_0^2 {r^{{p-2\over 2}} dr\over 
|1+(-1+i \eps)r|^p}\\[5pt] &\asymp&  \int_0^2 {r^{{p-2\over 2}} dr\over 
(|1-r|+r\eps )^p}\\[5pt] &\asymp&  \int_{-1}^1 {(r+1)^{{p-2\over 2}} dr\over 
(|r|+(r+1)\eps)^p}\\[5pt] &\asymp & \int_{{-1\over 2}}^{1\over 2}  {dr\over 
(|r|+\eps)^p}\\[5pt] &\asymp&  \eps^{-(p-1)}.\end{eqnarray*}
In the calculation above we used the first order approximation 
$e^{i(\pi-\eps)}\approx -1+i\eps$ for $\eps>0$, 
$\eps\to 0$ which, as one easily convinces oneself, is justified.   
 \enddemo

\demo{Case 3: ${p>1}$, ${q>0}$} In this case we have $\omega 
=\omega_1=2\alpha$ and so 
$\a_\C^0=\{ z\in\C\: |\Im z|<{\pi\over 4}\}$. Formula  (5.2) and Theorem 
4.2(iii) 
then 
imply for all $ \lambda\in i\R$, $-{\pi\over 4}<\phi<{\pi\over 4}$ and $t>0$ 
that 
\begin{eqnarray}
&&\phi_\lambda(t^{-1}e^{-i\phi})=Ct^{p+2q\over 
2}e^{-i\lambda\phi}\speqnu{5.5}\\ &&\times \int_0^\infty \int_0^\infty {r^{{p-2\over 2}} s^{{q-2\over 2}}dr \ ds\over 
[(1+t^2e^{i\phi} r)^2 +t^4e^{2i\phi} s]^{{p+2q\over 4}-{\lambda\over 2}}
[(1+e^{-i\phi} r)^2 +e^{-2i\phi} s]^{{p+2q\over 4}+{\lambda\over 2 
}}}.\nonumber
\end{eqnarray}

In particular (5.5) implies  that $\phi_\lambda\circ\log_A$
extends to a holomorphic function on $\b_\C^0=2\a_\C^0$. Again, this turns out 
to be the maximal 
domain, as we will show below. 
\enddemo

 \proclaim{Lemma}   For $p>1${\rm ,} $q>0$ and for all $\lambda\in i\a^*${\rm ,} $T_{\lambda,{\rm 
max}}=\b_\C^0=2\a_\C^0=\{ z\in\C\: |\Im z|<{\pi\over 2}\}$. 
Moreover{\rm ,} for a fixed $\lambda$ the asymptotics at the boundary are given by 
$$\phi_\lambda(e^{-i(\pi - \eps)}) \asymp \left\{ \begin{array}{ll}{1\over \eps^{q-1}} & \hbox{if $q>1$}, 
\\ |\log\eps| & \hbox{if $q=1$},\end{array}\right.$$
for $\eps\to 0$, $\eps>0$. 
\endproclaim

\demo{Proof} We will estimate $\phi_\lambda(e^{i({\pi\over 2} - \eps)})$ 
for $\eps\to 0$,  $\eps>0$. Since the  unbounded 
contribution of the integral is local (near $r=0$ and $s=1$), we may henceforth 
assume 
that $\lambda=0$. Then (5.5) gives that  
\begin{eqnarray*}
\phi_0(e^{-i({\pi\over 2}-\eps)})
&\asymp&\int_0^\infty \int_0^\infty {r^{{p-2\over 2}} s^{{q-2\over 2}}\ ds \ 
dr\over 
|(1+e^{i({\pi\over 2}-\eps)} r)^2 +e^{i(\pi-2\eps)} s|^{p+2q\over 2}}\\ &\asymp&\int_0^\infty \int_0^\infty {r^{{p-2\over 2}} s^{{q-2\over 2}}\
ds \  dr\over 
|1+2re^{i({\pi\over 2}-\eps)}+ e^{i(\pi -2\eps)} r^2 +e^{i(\pi-2\eps)} 
s|^{p+2q\over 2}}\\ &\asymp&\int_0^{1\over 2} \int_0^2 {r^{{p-2\over 2}} s^{{q-2\over 2}}\ ds \ 
dr\over 
|1+2r(i+\eps) + (-1 +i2\eps)r^2 + (-1+i2\eps) s|^{p+2q\over 
2}}\\ &\asymp&\int_0^{1\over 2}\int_0^2 {r^{{p-2\over 2}} s^{{q-2\over 2}}\ ds \ 
dr\over 
|(1 + 2r\eps -r^2-s) +i2(r+ \eps (r^2 + s))|^{p+2q\over 
2}}\\ &\asymp&\int_0^{1\over 2} \int_0^2 {r^{{p-2\over 2}} s^{{q-2\over 2}}\ ds \ 
dr\over 
\big| |1 + 2\eps r -r^2-s| +2 |\eps(r^2+s) +r|\big|^{p+2q\over 
2}}\\ &\asymp&\int_0^{1\over 2} \int_{-1}^1 {r^{{p-2\over 2}} (s+1)^{{q-2\over 2}}\ ds 
\ 
dr\over 
\big(|2\eps r -r^2-s| +2 \eps(r^2+s+1) +2r\big)^{p+2q\over 
2}}\\ &\asymp&\int_0^{1\over 2} \int_{-1}^1 {r^{{p-2\over 2}} \ ds \ 
dr\over 
\big(|2\eps r -r^2-s| +2 \eps(r^2+s+1) +2r\big)^{p+2q\over 
2}}.\end{eqnarray*}

Elimination of the absolute value in the integrand gives 
\begin{eqnarray*}
\phi_0 (e^{-i({\pi\over 2}-\eps)})&\asymp&\int_0^{1\over 2} \int_{-1}^{ 
2\eps r -r^2}
{r^{{p-2\over 2}} \ ds \ dr\over 
\big(2\eps r -r^2-s +2 \eps(r^2+s+1) +2r\big)^{p+2q\over 
2}}\\ && +\int_0^{1\over 2} \int_{2\eps r -r^2}^1 
{r^{{p-2\over 2}} \ ds \ dr\over 
\big(s-2\eps r +r^2 +2 \eps(r^2+s+1) +2r\big)^{p+2q\over 
2}}\\ & \asymp&\int_0^{1\over 2} \int_{-1}^{ 2\eps r -r^2}
{r^{{p-2\over 2}} \ ds \ dr\over 
\big(2\eps r -r^2 +2 \eps(r^2+1) +2r + s(-1+2\eps)\big)^{p+2q\over 
2}}\\ &&+\int_0^{1\over 2} \int_{2\eps r -r^2}^1 
{r^{{p-2\over 2}} \ ds \ dr\over 
\big(-2\eps r +r^2 +2 \eps(r^2+1) +2r + s(1+2\eps)\big)^{p+2q\over 
2}}\\ & \asymp&\int_0^{1\over 2}{r^{{p-2\over 2}} \ dr\over 
(2\eps +2r + 4\eps r^2)^{{p+2q\over 2}-1}} \\ & \asymp&\int_0^{1\over 2}{r^{{p-2\over 2}} \ dr\over 
(r+\eps)^{{p+2q\over 2}-1}}\\ & \asymp& \eps^{-{p+2q\over2}+1} \int_0^{1\over 2}{r^{{p-2\over 2}} \ dr\over 
({r\over \eps}+1)^{{p+2q\over 2}-1}}\\ & \asymp& \eps^{-{p+2q\over2}+1} \eps^{p-2\over 2}\int_0^{1\over 
2}{\big({r\over \eps}\big)^{{p-2\over 2}} \ dr\over 
({r\over \eps}+1)^{{p+2q\over 2}-1}}\\ & \asymp& \eps^{-{p+2q\over2}+1} \eps^{p-2\over 2}\eps\int_0^{1\over 
2\eps}{r^{{p-2\over 2}} \ dr\over 
(r+1)^{{p+2q\over 2}-1}}\\ & \asymp &\eps^{-q+1}\int_0^{1\over 2\eps}{r^{{p-2\over 2}} \ dr\over 
(r+1)^{{p+2q\over 2}-1}}\\ & \asymp &\eps^{-q+1}\int_1^{1\over 2\eps}r^{-q} dr\\ & \asymp&\left\{ \begin{array}{ll}{1\over
\eps^{q-1}} & \hbox{if $q>1$}, \\ |\log\eps| & \hbox{if $q=1$}.\end{array}\right.\\
\noalign{\vskip-36pt}
\end{eqnarray*}
 \enddemo
\vglue12pt

We remark that in order to obtain upper estimates only,  the assumption that 
$\lambda\in i\a^*$ was not used in view of 
the degree of generality of the formula in Theorem 4.2(iii). Collecting the 
preceding results we have proved   Theorem~5.1.

Everything that we will have proved about radial limits, namely Theorem~5.1 and 
Theorem  4.5, is consistent with the following conjecture.

\demo{Conjecture {\rm B}} Let $\alpha\in \Sigma_0$ and $H_\alpha\in \a$ be the corresponding 
{\it co-root}; i.e., 
$H_\alpha\in [\g^\alpha, \g^{-\alpha}]\cap \a$ such that $\alpha(H_\alpha)=2$.  
Let $c_\alpha\in \R$ such that $c_\alpha H_\alpha\in \partial \b^0$. Further,
set $m_\alpha=\dim \g^\alpha$ for all $\alpha\in \Sigma$. Then for all 
$\alpha\in \Sigma_0$ and $\lambda\in\a_\C^*$ we have 
$$|\phi_\lambda(\exp(i(1-\eps)c_\alpha H_\alpha))|\asymp\left\{ \begin{array}{ll} {1\over 
\eps^{m_\alpha-1}} & \hbox{if $m_\alpha>1$}, \\ |\log\eps| & \hbox{if $m_\alpha=1$}. \end{array}\right.
$$
\enddemo

\numbereddemo{{R}emark}  Correspondence with G. Heckman and E. Opdam suggests that the 
nature of the singularity of the holomorphically extended spherical function in 
co-root directions might be obtained from properties of the monodromy associated 
to solutions of the system of invariant differential operators. \enddemo

{\it Lower estimates}.
In a later application to automorphic functions we will also  need  lower 
estimates for the norm of the $K$-fixed vector in the holomorphically continued 
region, for all $\eps >0$, not only at the singularity. 
The result is obtained in a way similar to the preceding.

 \proclaim{Proposition} Let $(\pi_\lambda, {\cal H}_\lambda)$ be a unitary spherical 
principal
series representation of a group $G$ of real rank one. Let $X\in \partial \b^1$ 
and set 
$$(v_0)_\eps =\pi_\lambda(\exp(i(1-\eps)X))v_0.$$ 
Then there exists a constant $C$ independent of $\lambda$ such that 
$$\|(v_0)_\eps\|^2=|\phi_\lambda(\exp(-2i(1-\eps)X))|\geq C  \left\{ \begin{array}{ll} e^{(\pi- 
7\eps)|\lambda|} & \hbox{for $q=0$}, \\ \eps e^{({\pi\over 2}-21\eps)|\lambda|} & \hbox{for $ q>0$}, \end{array}\right.$$
for all $0<\eps\leq 1$.
\endproclaim

\demo{Proof} As usual we restrict ourselves to the case of $\lambda$ imaginary.

\vglue9pt {\it Case} 1: $q=0$. Here we have that 
$$ \|(v_0)_\eps\|^2 =\phi_\lambda(e^{i(\pi-\eps)})=e^{-i\lambda(\pi-\eps)}
\int_0^\infty {r^{p-2\over 2}\over |(1+e^{i(\pi-\eps)}r)^{{p\over 
2}-\lambda}|^2}.$$

By the Weyl group invariance of $\phi_\lambda$ we have 
$\phi_\lambda=\phi_{-\lambda}$ and so we may assume that 
$\lambda\in i\R^+$, i.e., $\lambda=i|\lambda|$. Then we get 
$$ \|(v_0)_\eps\|^2 \geq e^{(\pi-\eps)|\lambda|}
\int_0^{1\over 2} {r^{p-2\over 2}\ dr \over |(1+e^{i(\pi-\eps)}r)^{{p\over 
2}-\lambda}|^2}.$$
If $z$ is a complex number, then we write $-\pi\leq \arg(z)<\pi$ for the 
argument 
of $z$ and $m(z)$ for the modulus of $z$. Then for $0\leq r\leq {1\over 2}$ we 
have 
$0\leq \arg(1+e^{i(\pi-\eps)}r)<3\eps$ and $m(1+e^{i(\pi-\eps)}r)\leq 2$. Hence 
$${1\over |(1+e^{i(\pi-\eps)}r)^{{p\over 2}-\lambda}|^2}\geq 2^{-{p\over 2}} 
e^{-6\eps|\lambda|} $$
and the assertion of the proposition for $q=0$ follows.

\vglue9pt {\it Case} 2: $q>0$. Here we have that
\begin{eqnarray*}
 \|(v_0)_\eps\|^2& =&\phi_\lambda(e^{i({\pi\over 2}-\eps)})\\
&=&
e^{-i\lambda({\pi\over 2}-\eps)}
\int_0^\infty \int_0^\infty
{r^{p-2\over 2} s^{{p-2\over 2}}\ dr \ ds \over 
|\big((1+e^{i({\pi\over 2}-\eps)}r)^2+ se^{i(\pi-2\eps)}\big)^{{p+2q\over 4}-
{\lambda\over 2}}|^2}.
\end{eqnarray*}
By the Weyl group invariance of $\phi_\lambda$ we may assume that 
$\lambda\in i\R^+$ and  hence get 
$$ \|(v_0)_\eps\|^2\geq 
e^{({\pi\over 2}-\eps)|\lambda|}
\int_0^{1\over 2} \int_0^\eps
{r^{p-2\over 2} s^{{p-2\over 2}}\ dr \ ds \over 
\big|((1+e^{i({\pi\over 2}-\eps)}r)^2+ se^{i(\pi-2\eps)})^{{p+2q\over 4}-
{\lambda\over 2}}\big|^2}.$$

Now for $0\leq r\leq \eps$ and $0\leq s\leq {1\over 2}$ we have 
$0\leq \arg((1+e^{i({\pi\over 2}-\eps)}r)^2+ se^{i(\pi-2\eps)})\leq 10\eps$
and  $m((1+e^{i({\pi\over 2}-\eps)}r)^2+ se^{i(\pi-2\eps)})\leq 2$. Hence the 
assertion follows as in Case 1.  \enddemo

\section{Invariant seminorms}

Bernstein and Reznikov, in [BeRe99], introduced the notion of a maximal 
invariant seminorm associated to Sobolev 
norms of vectors in representations. For the $K$-fixed vector of spherical 
principal series representations for 
$G=\Sl(2,\R)$ they coupled this with some estimates on the holomorphically 
extended spherical functions into 
a beautiful technique to get estimates on Rankin-Selberg integrals for Maa{\ss} 
forms. 

We shall extend their technique in several directions. First, by using a more 
representation theoretic viewpoint we will be able to treat the case of real rank one groups. When 
specialized to $G=\Sl(2,\R)$ this will allow us to get a small improvement over 
the corresponding results in [BeRe99]. Secondly, in Section~9 we are able to consider 
some higher rank groups for which we obtain estimates on triple 
products of Maa{\ss} forms. These higher rank results are likely new, but should be 
viewed as a sample of the technique rather than as sharp results.

\advance\theoremcount by 1

\demo{Definition {\rm  6.1.  (cf.\ [BeRe99, App.\ A])}} 

(a) Let $V$ be a real or complex 
vector space and $(N_i)_{i\in I}$  a family of seminorms on it. Then   
 
$$(\inf_{i\in I} N_i)(v) \:=\inf_{\sum_{i\in I}  v_i=v}\sum_{i\in I} N_i(v_i)$$
also defines a seminorm on $V$ and satisfies $\inf_{i\in I} N_i\leq N_j$ for 
every $j\in I$. 

\vglue4pt 
(b) Let $G$ be a semigroup acting on $V$ and $N\: V\to [0,\infty[$ a 
single 
seminorm. 
Then for $g\in G$ define a seminorm $N_g$ by $N_g(v) =N(g\cdot v)$. 
As in (a) one obtains a seminorm $N^G$ by setting
$$N^G =\inf_{g\in G} N_g.
$$

\numbereddemo{Definition} Let $(\pi, {\cal H})$ be a unitary  representation of a Lie 
group $G$
on some Hilbert space ${\cal H}$.  Let $\{X_1, \ldots X_n\}$ be a basis of $\g$. 
Then the {\it $k^{\rm th}$ Sobolev norm} on ${\cal H}^\infty$ is defined by 
$$S_k(v) =\sum_{0\leq m_1+\ldots +m_n\leq k}
 \| d\pi(X_1^{m_1}\ldots X_n^{m_n})v\|\qquad (v\in {\cal H}^\infty).$$
It is easy to see that a different choice of basis leads to an equivalent 
seminorm. We remark that 
$S_k$, $k>1$,  is usually {\it not} $G$-invariant. As in (b) above, we set
$$S_k^G(v) =\inf_{g\in G}S_k(\pi (g)v).$$
\enddemo

Then it is a natural problem to estimate $S_k^G(\cdot)$ for the various 
representations of $G$. Fix an irreducible unitary  
representation of a semisimple Lie group $G$ having a nonzero $K$-fixed vector 
$v_0$. Let $v\in {\cal H}_{\lambda, K}$ be a $K$-finite 
vector. Recall from Proposition 4.1 that the orbit map 
$G\to {\cal H}_\lambda, \ g\mapsto \pi_\lambda(g)v$
extends to a holomorphic map on $G\Omega K_\C$. 
Write $\Omega=A\Omega_i$ with $\Omega_i\subeq \exp(i\a)$, and notice that 
$\Omega_i$ has compact closure.
We shall show for real rank one groups that $S_k^G(\pi_\lambda(a)v)$ 
is comparable to $\|\pi_\lambda(a)v\|$ uniformly in $a\in \Omega$ for all
$K$-finite vectors. Similar results will be obtained for holomorphic discrete 
series in Section~8. But first we explain how for spherical principal series the 
case of an arbitrary $K$-finite vector $v$ can be reduced 
to the spherical vector $v_0$.

\vglue9pt {\it Reduction to a spherical vector}.

 \proclaim{Lemma} Let $G = KAN$ be any Iwasawa decomposition and set $L=AN$.
Suppose that $(\pi_\lambda, {\cal H}_\lambda)$ is an irreducible unitary  
representation of a semisimple Lie group $G$ having a nonzero $K$-fixed vector 
$v_0$.   
\begin{itemize}
\ritem{(i)} The $K$\/{\rm -}\/spherical vector 
$v_0$ is $L$-cyclic{\rm ,} i.e. ${\cal H}_\lambda=\oline 
{\span_\C\{\pi_\lambda(L)v_0\} }.$
\ritem{(ii)} If ${\cal H}_{\lambda, K}$ denotes the $K$\/{\rm -}\/finite vectors of 
$(\pi_\lambda, {\cal H}_\lambda)${\rm ,} 
then 
$$  {\cal H}_{\lambda, K}= d\pi_\lambda({\cal U}(\l_\C))v_0, $$
where $\l$ denotes the Lie algebra of $L$. 
\end{itemize}

\endproclaim

\demo{Proof} (i) This follows from $\pi_\lambda(L)v_0=\pi_\lambda(G)v_0$ and the 
irreducibility of 
$(\pi_\lambda, {\cal H}_\lambda)$. 
\vglue4pt
(ii) This is immediate from (i). \enddemo

Let $(\pi,{\cal H})$ be a  Hilbert representation of $G$. For a closed subgroup 
$L<G$ write ${\cal H}_L^\infty$ for the smooth vectors for $\pi\res_L$. If $\pi$ 
is irreducible, then from the Casselman-Wallach theory of 
smooth globalizations of Harish-Chandra modules (cf.\ [Wal92, Ch.\ 11])
one has that ${\cal H}^\infty ={\cal H}_K^\infty.$

If $H<G$ is a subgroup, denote by $S_{k, H}$ the $k^{\rm th}$ Sobolev norm for the 
representation $\pi\res_H$. In particular, the Fr\'echet topology on ${\cal 
H}^\infty$ is also induced by the Sobolev norms $(S_{k,K})_{k\in \N}$.

 \proclaim{Lemma} Let $G = KAN$ be any Iwasawa decomposition and set $L=AN$. Suppose 
that $(\pi_\lambda, {\cal H}_\lambda)$ is an irreducible unitary  
representation of a semisimple Lie group $G$ having a nonzero $K$\/{\rm -}\/fixed vector 
$v_0$. 
\begin{itemize}
\ritem{(i)} For every $k\in \N$ there exist  an $l\in \N$ and a constant $C>0$ 
 such that 
$$(\forall a\in \Omega_i) \qquad 
S_k(\pi_\lambda(a)v_0)\leq C S_{l, L}(\pi_\lambda(a)v_0).$$

\ritem{(ii)} For every $v\in {\cal H}_{\lambda, K}$ and $k\in \N$ there exist  an 
$l\geq k$ and a constant $C>0$ such that 
$$(\forall a\in \Omega_i) \qquad 
S_k(\pi_\lambda(a)v)\leq C S_{l, L}(\pi_\lambda(a)v_0).$$
\end{itemize}
\endproclaim

\demo{Proof} (i) We identify ${\cal U}(\g_\C)$ with ${\cal S}(\g_\C)$. Then the 
natural grading of 
${\cal S}(\g_\C)$ yields a direct sum decomposition ${\cal 
U}(\g_\C)=\bigoplus_{k\in \N} {\cal U}(\g_\C)^k$. 
Fix a norm $\|\cdot\|$ on $\g_\C$ and take its natural extension to ${\cal 
S}(\g_\C)$. For any $g\in G_{\Bbb C}$, $\Ad(g)$ maps ${\cal U}(\g_\C)^k$ to 
itself boundedly, so has a norm, say, $\|\Ad(g)\|_k$.   
If $X\in {\cal U}(\g_\C)^k$ with $\|X\|=1$, 
then  
$$\| X\pi_\lambda(a)v_0\|=\|\pi_\lambda(a) (\Ad(a)^{-1}X)v_0\|\leq 
\|\Ad(a^{-1})\|_k \sup_{Y\in 
{\cal U}(\g_\C)^k\atop \|Y\|\leq 1} \|\pi_\lambda(a)Yv_0\|.$$
Here $C =\sup_{a\in \Omega_i}\|\Ad(a^{\pm 1})\|_k$ is finite by the relative 
compactness of $\Omega_i$. 
Hence from Lemma 6.3 there exist an $l\in \N$ and an $r>0$ such 
that 
$$\| X\pi_\lambda(a)v_0\|\leq C \sup_{Y\in 
{\cal U}(\g_\C)^k\atop \|Y\|\leq 1} \|\pi_\lambda(a)Yv_0\|\leq C \sup_{Z\in 
{\cal U}(\l_\C)^l\atop \|Z\|\leq r} \|\pi_\lambda(a)Zv_0\|.$$
Now as $\l$ is normalized by $\a$, we get that 
$$\| X\pi_\lambda(a)v_0\|\leq C^2 \sup_{Z\in 
{\cal U}(\l_\C)^l\atop \|Z\|\leq r} \|Z\pi_\lambda(a)v_0\|\leq C' 
S_{l,L}(\pi_\lambda(a)v_0)$$
for some constant $C'$ independent of $X$. 
\vglue8pt

(ii) By Lemma 6.3(ii) there exists an $X\in {\cal U}(\l_\C)$ such 
that $v=Xv_0$. Since $\a$ normalizes $\l$ the assertion follows 
now from (i). \enddemo

Throughout this section we shall follow the custom that a constant `C' 
depends on any quantifiers preceding it in the statement. Thus in the previous 
result (ii), `C' depends on $\pi_\lambda$, $k$, and $v$ but not on $a$.

\vglue12pt {\it Compressing Sobolev norms}.
For any choice of positive roots $\Sigma^+$ we set $\a^+ =\{ X\in \a\: (\forall 
\alpha\in \Sigma^+)\ \alpha(X)>0\}$ and $\a^- 
=-\a^+$, and, on the group side, let $A^{\pm} =\exp(\a^\pm)$. 

 \proclaim{Lemma} Let $(\pi, {\cal H})$ be a unitary representation of $G$ and 
$v\in {\cal H}^\infty$. Then for $k\in \N_0${\rm ,} 
\begin{itemize}
\ritem{(i)} $S_{k,N}^{A^+}(v)=\|v\|;$
\ritem{(ii)} $S_{k,AN}^G(v)=S_{k,A}^G(v)$. 
\end{itemize}
\endproclaim

\demo{Proof} (i) Let $\{X_1, \ldots,X_s\}$ be a basis of root vectors  
of $\n$ corresponding to  roots $\alpha_1, \ldots, \alpha_s\in \Sigma^+$. 
Then for any $v\in {\cal H}^\infty$  
$$S_{l,N} (v)=\|v\|+\sum_{1\leq m_1+\ldots+m_s\leq l} \|d\pi(X_1^{m_1}\cdot
\cdot X_s^{m_s})v\|.$$
For $a=\exp(X)\in A$,  \pagebreak $X\in\a$,
\begin{eqnarray*}
&&\hskip-.5in S_{l,N} (\pi(a)v)\\
&  
=&\|\pi(a)v\|+\sum_{1\leq m_1+\ldots+m_s\leq l} \|d\pi(X_1^{m_1}\cdot
\cdot X_s^{m_s})\pi(a)v\|\\ &=&\|v\|+\sum_{1\leq m_1+\ldots+m_s\leq l} 
\|\pi(a)d\pi((\Ad(a^{-1})X_1)^{m_1}\cdot
\cdot (\Ad(a^{-1})X_s)^{m_s})v\|\\ &=&\|v\|+\sum_{1\leq m_1+\ldots+m_s\leq l}e^{-\sum_{j=1}^s m_j\alpha_j(X)} 
\|d\pi(X_1^{m_1}\cdot
\cdot X_s^{m_s})v\|.\end{eqnarray*}
If we choose $X\in \a^+$, 
\begin{eqnarray*}
S_{l,N}^{A^+} (v)&\leq& \inf_{a\in A^+}S_{l,N}(\pi(a)v)\leq 
\inf_{t>0}S_{l,N}(\pi(\exp(tX))v)\\ &=&\inf_{t>0} \Big(\|v\|+\sum_{1\leq m_1+\ldots+m_s\leq l}e^{-t\sum_{j=1}^s 
m_j\alpha_j(X)} \|d\pi(X_1^{m_1}\cdot
\cdot X_s^{m_s})v\|\Big)\\
&=&\|v\|.
\end{eqnarray*}
On the other hand, clearly $\|v\|\leq S_k^{A^+}(v)$. Thus 
$\|v\|=S_{l,N}^{A^+}(v)$ completing 
the proof of (i).
\vglue4pt
(ii) One has the obvious inequality $S_{k,AN}(v)\geq S_{k,A}(v)$, so that 
$S_{k,AN}^G(v)\geq S_{k,A}^G(v)$.  On the other hand,
\begin{eqnarray*}
S_{k, AN}^{G} (v)&\leq &\inf_{g\in G}S_{k, AN}(\pi(g)v)\\ &\leq &\inf_{h\in A^+}S_{k, AN}(\pi(h)v)\\ &=& {S_{k, A}}(v),\end{eqnarray*}
so that $S_{k,AN}^G(v)\leq S_{k,A}^G(v)$. 
 \enddemo

{\it The case of $G=\Sl(2,\R)$}.
Our goal is to estimate $S_k^G (\pi (a)v_0)$ for all
$a\in \Omega_i$. In this section we shall present extensive details for 
$G=\Sl(2,\R)$ as this will be the model for the proof later for rank one 
groups. Here we will consider an irreducible unitary  
spherical principal series representation $(\pi_\lambda, {\cal H}_\lambda)$. The 
complementary series and nonspherical principal series representations can be 
shown similarly. Discrete series however will be obtained rather differently in Section~8.  

We identify $N$ with $\R$ via the mapping $n_x\mapsto x$ (see Appendix A for 
notation). We are going to work in the noncompact realization of $\pi_\lambda$
on $L^2(N)=L^2(\R)$. With $g^{-1}=\left(\begin{array}{cc} a & b\\ c & d\end{array}\right)$ the action 
of $\pi_\lambda(g)$ is given by 
$$(\pi_\lambda(g)f)(x)=|cx+d|^{\lambda -1} f\Big({ax +b \over cx 
+d}\Big)\leqno(6.1)$$
for all $f\in L^2(\R)$ and $x\in\R$. 
For this module one has 
$${\cal H}_\lambda^\infty=\{ f\in C^\infty(\R)\: |x|^{\lambda -1} f({1\over 
x})\in C^\infty(\R)\}.$$

We use a usual basis for the Lie algebra of $\g$:  
$$H=\left(\begin{array}{cc} 1 & 0 \\ 0 & -1\end{array}\right), \qquad E=\left(\begin{array}{cc}  0 & 1\\ 0 & 0\end{array}\right), 
\qquad F=\left(\begin{array}{cc} 0 & 0 \\ 1& 0\end{array}\right).$$
Then $\a=\R H$, $\n=\R E$ and $\oline \n=\R F$. With $U=E-F$ we have $\k=\R U$. 
Differentiating (6.1) one obtains the formulas  \advance\eqcount by 1
\begin{eqnarray}
d\pi_\lambda(H) &=& (\lambda -1 ) -2 x{d\over dx},  \\ d\pi_\lambda(E)&= &- {d\over dx},  \\ 
d\pi_\lambda(F)&=&
(1-\lambda)x + x^2 {d\over dx},  \\ d\pi_\lambda(U) &=& (\lambda-1) - (1+x^2) {d\over dx},   \\ d\pi_\lambda(E+F) &=& (1-\lambda)x -
(1-x^2) {d\over dx}.
\end{eqnarray}

We also define the {\it radial  operators} 
by 
$$(R_jf)(x)=(x^j {d^j\over dx^j}f)(x)$$
and define the {\it radial Sobolev norms} by 
$$S_{k, {\rm rad}}(f)=\sum_{j=0}^k \|R_j f\|.$$
From the action of $d\pi_\lambda(H)$ and $R^j$ it is clear that 
there exists  a constant $C>0$, depending on $k$ and $\lambda$, such that for 
all $f\in {\cal S}(\R)$
$${1\over C} S_{k,{\rm rad}} (f)\leq 
S_{k,A}(f)\leq C S_{k, {\rm rad}}(f).\leqno(6.7)$$

As remarked by the referee, in (6.2) and (6.4) the coefficient of the derivative 
term has a zero; consequently $S_k(v)$ cannot be majorized by $S_{k,A\oline 
N}(v)$ or by $S_{k,A}(v)$ in general. However, we shall show in the next proposition that there is such a relationship for the $G-$invariant
Sobolev  norms.  
 
 \proclaim{Proposition}  Let $G=\Sl(2,\R)$ and $(\pi_\lambda,{\cal H}_\lambda)${\rm ,} $\lambda\in i\a^*${\rm ,} be an 
irreducible unitary spherical principal series representation. Then for every 
$k\in \N_0$ there exists a $C>0$ such that for $v\in {\cal H}_\lambda^\infty,$ 
$$ S_k^G (v)\leq C S_{k,A}^G(v).$$
\endproclaim

\demo{Proof} The $A$ action on $K/M\cong S^1$ has two fixed points, corresponding to 
the two Bruhat cells. In the noncompact realization $N$ they become the origin 
and the point at infinity. We shall estimate $ S_k^G (f)$ by using first a 
cutoff function at infinity, $\oline \n$, and an elementary estimate there. Near 
the origin a dilated cutoff localizes sufficiently high derivatives of $f$ to 
get an estimate. Away from the fixed points, motivated by an argument in [BeRe99] and classical Littlewood-Paley theory, we use a family of suitably dilated cutoff functions which compress the $\n$ derivatives in the definition of $G$-invariant norm to {\it radial} derivatives thereby obtaining the desired 
estimate. 

For $j\in \Z$ we denote by $I_j$ the set $\{ x\in\R \: 2^{-j-1}\leq |x|\leq 
2^{-j+1}\}$. 
For a function $\psi$ on $\R$ we write $\psi_j(x)=\psi(2^j x)$. Notice that if 
$\psi$ is supported in $I_0$ then $\psi_j$ is supported in $I_j$, and 
$${\rm supp}(\psi_j)\ \cap \ {\rm supp} (\psi_{j+1}) \subseteq \pm\lbrack {1\over 
2^{j+1}}, {1\over 2^j}\rbrack.$$  
We take a smooth, nonnegative function $\phi$ supported in $I_0$ and such that 
for every $m\in \N_0$, 
$$\sum_{j=0}^m \phi_j(x)=\left\{ \begin{array}{ll}0 & \hbox{if $|x|\leq  2^{-m-1}$}, \\ 1 & \hbox{if $2^{-m}\leq |x|\leq 1$},\\ 0 & \hbox{if
$2\leq |x|$}.\end{array}\right.$$

Choose a nonnegative function $\tau\in C^\infty(\R)$ with support in $\{ 
x\in \R\: \break 1\leq |x|\}$ 
such that $(\tau+\phi)(x)=1$ for $|x|\geq 1$. 
Finally for each $m\in \N$ define the function $\tau_m\in C_c^\infty(\R)$ by 
$\tau_m =\1 -\tau-\sum_{j=0}^m \phi_j$. Notice that $\supp 
\tau_m\subeq\{x\in\R\:  |x|\leq 2^{-m}\}$
and $\tau_m(x)=1$ for $|x|\leq 2^{-m-1}$. From the properties of the $\phi_j$ 
and $\tau$ it is easy to see that for any $l\geq 1$, $\tau^{(l)}_m (x) = 
-2^{lm}\phi^{(l)}(2^mx)$.

Let $f\in {\cal H}_\lambda^\infty$. Since
\begin{eqnarray*}
\1 & = & \tau + \1 - \tau\\ & =& \tau + \tau_m + \sum_{j=0}^m 
\phi_j\\ & = &\tau + \phi + \tau_m + \sum_{j=1}^m \phi_j,\\
\noalign{\noindent then}
 f &=& (\tau + \phi)f + \tau_mf + \sum_{j=1}^m {\phi_jf}.
\end{eqnarray*}
For any choices of $g, g_1,\ldots, g_m\in G$, using the definition of $S_k^G$, 
we get
$$S_k^G(f)\leq S_k((\tau +\phi) f) + S_k(\pi_\lambda(g)(\tau_m f)) + \sum_{j=1}^m 
S_k(\pi_\lambda(g_j) (\phi_j f)).
\leqno(6.8)$$

First we consider the term $S_k((\tau +\phi)f)$. From an examination of 
formulas (6.2)--(6.4) one sees that $S_k((\tau +\phi) f)\leq C S_{k,\oline 
N}((\tau +\phi)f)$
for all $f\in {\cal H}_\lambda^\infty$. (Throughout this proof $C$ will denote a constant depending only on $k$, $\tau$, $\phi$ and $\lambda$.)
Hence we have 
$$S_{k}((\tau +\phi)f)\leq CS_{k,\oline N}((\tau +\phi) f)\leq C S_{k,\oline 
N}(f)$$
for all $f\in {\cal H}_\lambda^\infty$.  Majorizing this term in (6.8) we get 
$$S_k^G(f)\leq C S_{k,\oline N}(f) + S_k((\pi_\lambda(g)\tau_m f))+ 
\sum_{j=1}^m S_k(\pi_\lambda(g_j) (\phi_j f))\leqno(6.9)$$
for all $f\in {\cal H}_\lambda^\infty$.

Next we specify a good choice of the elements 
$g, g_1,\ldots, g_m\in G$. For every $t>0$ denote by $b_t$ the element 
$$b_t =\left(\begin{array}{cc} {1\over \sqrt t} & 0 \\ 0 & \sqrt t\end{array}\right)\in A.$$
From (6.1) it follows that 
$$(\pi_\lambda(b_t)f)(x)= t^{{1\over 2}(1-\lambda)} f(tx)$$
for all $t>0$ and $x\in \R$. 
Take $g_j=b_{2^{-j}}$ for all $1\leq j\leq m$ and $g = b_{2^{-(m+1)}}$. 
Notice that for every $m$ all the $\pi_\lambda(g_j) (\phi_j f)$ are  
supported in $[-2,2]$, as is $\pi_\lambda(g) (\tau_m f)$. For any smooth 
function $h$ supported in $[-2,2]$ we can conclude from the formulas (6.2)--(6.5) that
$S_k(h)\leq CS_{k,N}(h)$. Using this in (6.9) we get
$$S_k^G(f)\leq C S_{k,\oline N}(f) +C S_{k,N}(\pi_\lambda(g) (\tau_mf)) + 
C\sum_{j=1}^m S_{k,N}(\pi_\lambda(g_j) (\phi_j f))\leqno(6.10)$$
for all $f\in {\cal H}_\lambda^\infty$. 

Estimating  $S_{k,N}(\pi_\lambda(g) (\tau_mf))$, we use Leibniz 
on $\tau_m f$ and $L^\infty$ estimates on $\tau_m^{(j)} = 
-2^{jm}\phi^{(j)}(2^mx)$. From (6.3) one sees that $S_{k,N}(h)=\sum_{l=0}^k 
\|h^{(l)}\|$. Then
\begin{eqnarray}
&& \speqnu{6.11}\\
&&\hskip-.25in S_{k,N}(\pi_\lambda (g) (\tau_m f))  = \sum_{l=0}^k  \| 
{d^l\over dx^l} 2^{-{(m+1)\over 2}(1-\lambda)} (\tau_m f)(2^{-(m+1)}\cdot)\|\nonumber \\ & &=\sum_{l=0}^k \vert 2^{-{(m+1)\over
2}(1-\lambda)}\vert\nonumber\\
&&\quad \times \left [\int\Big|\sum_{n=0}^l 2^{-(m+1)l}{l \choose l-n} \tau_m^{(l-n)}(2^{-(m+1)}x) f^{(n)} 
(2^{-(m+1)}x)\Big|^2 \ dx\right]^{1\over 2} \nonumber \\ && \leq\sum_{l=0}^k \vert 2^{-{(m+1)\over 2}(1-\lambda)}\vert
\nonumber\\
&&\quad\times \sum_{n=0}^l
\Big[\int_{|x|\leq 2}\Big| 2^{-(m+1)l}{l \choose l-n} \tau_m^{(l-n)}(2^{-(m+1)}x) 
f^{(n)}(2^{-(m+1)}x)\Big|^2 \ dx \Big]^{1\over 2}\nonumber \\ && =\sum_{l=0}^k \Big|2^{{(m+1)\over 2}\lambda}\Big| \sum_{n=0}^l
\Big[\int_{|y|\leq {1\over 2^m}} \Big|2^{-(m+1)l} {l \choose l-n} \tau_m^{(l-n)}(y) f^n(y) \Big|^2\  dy 
\Big]^{1\over 2}\nonumber \\ && \leq\sum_{l=0}^k | 2^{{(m+1)\over 2}\lambda}| \sum_{n=0}^l{l \choose l-n}{ \|2^{(l-n)m}\phi^{(l-n)}\|_\infty\over
2^{(m+1)l}}
\Big[\int_{|y|\leq {1\over 2^m}} |f^{(n)} (y) |^2\  dy \Big]^{1\over 2}\nonumber \\ && =\sum_{n=0}^k | 2^{{(m+1)\over 2}\lambda}| {1\over 2^{mn}}
\sum_{l=n}^k{l \choose l-n}{\|\phi^{(l-n)}\|_\infty\over 2^l} \Big[\int_{|y|\leq {1\over 2^m}} 
|f^{(n)}(y)|^2\ dy\Big]^{1\over 2}\nonumber \\ && =\sum_{n=0}^k | 2^{{(m+1)\over 2}\lambda}| {1\over 2^{(m+1)n}}
\sum_{j=0}^{k-n}{j+n \choose n} {\|\phi^{j}\|_\infty\over 2^j} \Big[\int_{|y|\leq {1\over 2^m}} 
|f^{(n)}(y)|^2\ dy\Big]^{1\over 2}\nonumber \\ && \leq\big(\sum_{j=0}^k{\|\phi^{(j)}\|_\infty\over {j!2^j}}\big)\sum_{n=0}^k 
{k!\over 
n!2^{(m+1)n}}\Big[\int_{|y|\leq {1\over 2^m}}| f^{(n)}(y)|^2 \ dy\Big]^{1\over 
2}.\nonumber\end{eqnarray}
Now $k$ is fixed and each of the at most $k$ derivatives $f^{(n)}$ is in $L^2$, 
hence the integrals can be made uniformly small. So for each $f$ we can choose 
an $m$ so that the last line above is at most $\|f\|$. Then we have
 
$$S_k^G(f)\leq C S_{k,\oline N}(f) + C\|f\| + 
C \sum_{j=1}^m S_{k,N}(\pi_\lambda(g_j) (\phi_j f))$$
for any $f\in {\cal H}_\lambda^\infty$.   Thus from (6.10) we obtain   
$$S_k^G(f)\leq C S_{k,\oline N }(f) + C\|f\| +
C\sum_{l=0}^k \sum_{j=1}^m \| {d^l\over dx^l} (2^{-{j\over 2}(1-\lambda)}\phi 
f(2^{-j}\cdot))\|.\leqno(6.12)$$
As in (6.11), using Leibniz on $\phi f$, $L^\infty$ estimates on $\phi^{(j)}$, and majorizing the binomial coefficients, we get
\begin{eqnarray}
\qquad\quad \sum_{l=0}^k \sum_{j=1}^m \| {d^l\over dx^l} (2^{-{j\over 2}}\phi 
f(2^{-j}\cdot))\|
&\nhs  \leq \nhs&C 
\sum_{l=0}^k \sum_{j=1}^m \Big(\int_{I_0}  2^{-j-2l}|f^{(l)} (2^{-j}x)|^2 \ 
dx\Big)^{1\over 2}\speqnu{6.13}
\\ & \nhs=\nhs & C \sum_{l=0}^k \sum_{j=1}^m \Big(\int_{I_j} 2^{-2l} |f^{(l)} (x)|^2 \ 
dx\Big)^{1\over 2}\nonumber\\ &\nhs\leq\nhs& 4  C \sum_{l=0}^k \sum_{j=1}^m \Big(\int_{I_j}  | x^lf^{(l)} (x)|^2 \ 
dx\Big)^{1\over 2}\nonumber\\ &\nhs\leq \nhs&4  C S_{k,{\rm rad}} (f) \leq 4 C S_{k,A} (f), \nonumber
\end{eqnarray}
where the last inequality follows from (6.7) and again $ C$ depends only on $\tau$, $\phi$, $k$ and $\lambda$. 
Thus we get from (6.12) and (6.13) that 
$$S_k^G (f) \leq C S_{k,\oline N}(f) +  C\|f\|  + C S_{k,A}(f)\leq  C\|f\| + C 
S_{k,A\oline N}(f)$$
for all $f\in {\cal H}_\lambda^\infty$. Now, 
$$S_k^G \leq C S_{k,A\oline N}^G$$ and, by Lemma 6.5(ii),  $S_k^G \leq C 
S_{k,A}^G$ as was to be shown.
 \enddemo

In $G=\Sl(2,\R)$ the element 
$$k_0={1\over \sqrt 2}\left( \begin{array}{cccc}1 & 1 \\ -1 & 1\end{array}\right)$$
is in $K$ and is a square root of the
Weyl group element. It will turn out that $k_0$ provides a uniform minimizer for 
Sobolev norms.

\proclaim{Theorem} Let $G=\Sl(2,\R)$ and $(\pi_\lambda,{\cal H}_\lambda)${\rm ,} $\lambda\in i\a^*${\rm ,} be an 
irreducible unitary spherical principal series representation.
Then for every $k\in \N_0$ there exists a $C>0${\rm ,} depending on k and $\lambda${\rm ,} 
such that for all $a\in \Omega_i$
$$ S_{k,A}(\pi_\lambda(k_0)\pi_\lambda(a)v_0)
\leq C \|\pi_\lambda(a)v_0\|.$$
In particular{\rm ,}   for all $a \in \Omega_i$
$$ S_k^G(\pi_\lambda(a)v_0)\leq  C \|\pi_\lambda(a)v_0\|.$$
\endproclaim

\demo{Proof} In view of Proposition 6.6 the second assertion follows 
from the first one. 
To prove the first assertion notice that 
$$(\pi_\lambda(k_0)f)(x)=|x+1|^{\lambda-1} f\Big({x-1\over x+1}\Big)
\leqno(6.14)$$
for all $f\in L^2(\R)$.  

We parametrize $\Omega_i$ with $a_\eps$ and
$a_\eps^{-1}$, where  
$$a_\eps =\left( \begin{array}{cccc}e^{i{\pi\over 4}(1-\eps)} & 0 \\ 0 & 
e^{-i{\pi\over 4}(1-\eps)}\end{array}\right)$$
for $0<\eps\leq 1$. 
Then, in the noncompact realization, $\pi_\lambda(a_\eps)v_0$ is of the form 
$c(\lambda,\eps) f_\eps$ where 
$$f_\eps(x)= {1 \over (1+ e^{i\pi(1 -\eps)} x^2)^{{1\over 2}-\lambda}}$$
and $c(\lambda,\eps)$ is a constant depending on $\lambda$ and 
$\eps$, and is uniformly bounded in $\eps$ (as can be seen from \S 5). 
Notice that the poles of $f_\eps$, as $\eps\to 0$, are 
at $x=\pm 1$. Thus if we take a smooth cut-off function  $\tau\in 
C_c^{\infty}(\R)$
with, say, $\tau\res_{[-2,2]}=1$, then 
$$ S_{k,A}(\pi_\lambda(k_0)f_\eps)
\leq S_{k,A}(\pi_\lambda(k_0) \tau f_\eps)+ S_k(\pi_\lambda(k_0) (1-\tau) 
f_\eps)
\leq S_{k,A}(\pi_\lambda(k_0) \tau f_\eps) +  C. \leqno(6.15)$$
Here $C$ is a positive constant independent of $\eps$ because, on the support of 
$(1-\tau)$, one has $\|(1-\tau)f_\eps\|\asymp\|(1-\tau)x^{-1}\|$, with similar 
results on the norms of derivatives.

With $g_\eps =\tau f_\eps$, in view of (6.7), (6.14) and (6.15), it 
suffices to show that 
$$S_{k,{\rm rad}}(\pi_\lambda(k_0)g_\eps)\leq C \|\pi_\lambda(a_\eps)v_0\|$$
for all $\eps$ and some constant $C>0$.  
By the radial Sobolev norms and the estimate in Lemma 5.3,
$ \|g_\eps\|\asymp\|\pi_\lambda(a_\eps)v_0\|\asymp \sqrt{|\log\eps|}$, 
it is enough to show that 
$$\|R_j \pi_\lambda(k_0)g_\eps\|\leq C_j \sqrt{|\log\eps|}$$
for all $j\geq 1$.

{} For $f\in C_c^\infty(\R)$ and  from (6.14),
\begin{eqnarray*}
(R_1 \pi_\lambda(k_0) f)(x)&= &2x |x+1|^{\lambda-1} 
{1\over (x+1)^2} f'\Big({x-1\over x+1}\Big)\\
&& + 
 \eps(x) (\lambda-1) x|x+1|^{\lambda-2}
f\Big({x-1\over x+1}\Big)\end{eqnarray*}
with $\eps(x)=1$ for $x>-1$ and $\eps(x)=-1$ for $x<-1$. 
Disregarding the sign function $\eps(x)$, and using induction we have 
$$(R_j\pi_\lambda(k_0)f)(x)= x^j\sum_{m=0}^j c_m f_m^j(x)$$
for some constants $c_m$ independent of $f$, and where
$$f_m^j(x)= |x+1|^{\lambda-1-j-m} f^{(m)}\Big({x-1\over x+1}\Big).$$
Thus to estimate $S_{k, {\rm rad}}(\pi_\lambda(k_0) g_\eps)=\sum_{j=0}^k \|R_j 
\pi_\lambda(k_0) g_\eps\|$ we must show that 
$$|\la x^j g_{\eps,m}^j, x^j g_{\eps, n}^j\ra|\leq C |\log\eps|\leqno(6.16)$$
for all $m,n\leq j$. 

Now, consider an expression of the form $|\la x^jf_m^j, x^jf_n^j\ra|$ where
\begin{eqnarray}
&\nhs\nhs&\speqnu{6.17}\\
|\la x^jf_m^j, x^jf_n^j\ra| &\nhs\leq \nhs& \int_\R  x^{2j}\
|x+1|^{-2-m-n-2j}
\ |f^{(m)}\Big({x-1\over x+1}\Big)|\ |f^{(n)}\Big({x-1\over x+1}\Big)|
\ dx\nonumber
\\ &\nhs=\nhs&2 \int_\R  \Big|{x+1\over 1-x}\Big|^{2j}\ \Big|{x+1\over 1-x} +1
\Big|^{-m-n-2j}\ |f^{(m)}(x)f^{(n)}(x)|\ dx\nonumber \\ &\nhs=\nhs& 2^{-(m+n)-2j +1} 
\int_\R  |x+1|^{2j}\ |1-x|^{m+n}\ |f^{(m)}(x)f^{(n)}(x)|\ dx.\nonumber \end{eqnarray}

Next, as $\eps\to 0$, the functions $g_\eps$ have poles at
$x=1$ and $x=-1$. Similarly, as $\eps\to 0$, $g_\eps^{(m)}(x)$ has poles 
only at $x=\pm 1$ and of 
order at most $m+{1\over 2}$. Examining (6.17) with $f=g_\eps$
we see that the factor $|x+1|^{2j} |1-x|^{(m+n)}|$
cancels poles. In particular
$$|x+1|^{2j} |1-x|^{m+n}|g_\eps^{(m)}(x)g_\eps^{(n)}(x)|$$
has poles at $x=\pm 1$ for $\eps\to 0$ of order no more than that of $g_\eps$. 
This establishes (6.16) and concludes the proof of the 
theorem. \enddemo

\numbereddemo{{R}emark}   (a) The second estimate, $ S_k^G(\pi_\lambda(a)v_0)\leq  C 
\|\pi_\lambda(a)v_0\|$, in Theorem 6.7 is optimal
in the sense that by $G$-invariance one has $\|v\|\leq S_k^G(v)$
for any smooth vector $v$ in a unitary representation 
$(\pi,{\cal H})$ of $G$.
\vglue4pt
(b) One can modify the  proof of Theorem 6.7 to give the more general result 
$$(\forall a\in \Omega_i)\qquad S_{k,A}(\pi_\lambda(k_0)\pi_\lambda(a)v)
\leq C \|\pi_\lambda(a)v\|$$
for an arbitrary $K$-finite vector $v$.
\vglue4pt
(c) The estimate $ S_k^G(\pi_\lambda(a_\eps)v_0)\leq  C 
\sqrt{|\log\eps|}$ in Theorem 6.7 is a little sharper 
than the estimate (0.5) in [BeRe99], viz. $ S_k^G(\pi_\lambda(a_\eps)v_0)\leq  C 
|\log\eps|$.
\vglue4pt
(d) Theorem 6.7 can be easily generalized to complementary 
series using the results on spherical 
functions in Theorem 4.2. \enddemo 

Part of the method for $G=\Sl(2,\R)$ generalizes to all groups of real rank one. For example, the element $k_0\in K$ can
 be found in these groups and gives a uniform minimizer for $S_{k,A}$.

 \proclaim{Lemma} Let $\g$ be a semisimple Lie algebra
with Iwasawa decomposition $\g=\k\oplus\a\oplus\n$. 
Suppose that the restricted root system $\Sigma$ satisfies 
one of the following assumptions\/{\rm :}\/
\begin{itemize} 
\ritem{(1)} $\Sigma$ is of type $A_1$ or $BC_1${\rm ,} i.e.{\rm ,} $\g$ is 
of real rank one{\rm ;}  
\ritem{(2)} $\Sigma$ is of type $C_n$ or $BC_n$ for $n\geq 2$. 
\end{itemize}
Then there exists a $k_0\in K$ such that 
$$\Ad(k_0)\a\subeq \k\oplus\oline \n =\k\oplus \n.$$
\endproclaim

\demo{Proof} First recall that all maximal abelian 
subspaces in $\p$ are conjugate under $\Ad(K)$.

Suppose that (1) is satisfied. Then $\a$ is one-dimensional. 
Pick a nonzero root vector $X_\alpha\in\g^\alpha$. Then 
$\e=\R(X_\alpha-\theta(X_\alpha))$ is a maximal abelian 
subspace in $\p$ which lies in $\k\oplus\oline\n$. Hence there 
exists a $k_0\in K$ such that $\Ad(k_0)\a=\e$. 

Suppose then that (2) is satisfied. Since $\Sigma$ is of type 
$C_n$ or $BC_n$ we can find a maximal set $\gamma_1,\ldots,\gamma_n$
of long strongly orthogonal roots. But via $\sL(2,\R)$-reduction, 
the assertion follows from the already established rank one 
case above. \enddemo

We shall make the standing assumption, for the rest of this subsection, that $G$ 
has real rank one. We need to make the element $k_0$ more explicit. Let $\beta$ denote the long positive root. Then we have $\beta=\alpha$ if 
$q=0$, otherwise $\beta=2\alpha$. Choose an $\sL(2,\R)$-triple
$\{ E,F, H\}$ in $\g$ such that $E$ lies in the root space
$\g^\beta$, $F= -\theta E$ and such that 
$$H=[E,F] \qquad [H,E]=2E\qquad [H,F]=-2F.$$

With $U=E-F$ we choose 
$$k_0 =\exp({\pi\over 4} U).$$
Then 
$$\Ad(k_0)^{-1} H= E+F.$$
Notice that 
$$\Omega_i=\{ \exp(i\phi H)| \phi\in ]-{\pi\over 4}, {\pi\over 4}[\}$$
and introduce  elements $a_\eps$ by
$$a_\eps=\exp(i{\pi\over 4}(1-\eps)H).$$

 \proclaim{Proposition}  Suppose that $G$ is of real rank 
one and that $(\pi, {\cal H})$ is an irreducible 
unitary representation with $K$\/{\rm -}\/spherical vector $v_0$. 
Then for all $k\in\N_0$ there exists a constant 
$C>0$ such that 
$$(\forall a\in \Omega_i) \qquad S_{k,A}(\pi(k_0)\pi(a)v_0)
\leq C  \sum_{j=0}^k |a^\beta +a^{-\beta}|^j  S_j(\pi(a)v_0).$$
\endproclaim

\demo{Proof} $S_{k,A}$ is given by  
$$S_{k,A}(v)=\sum_{j=0}^k \|H^j v\|.$$
  
We are going to prove the proposition 
by induction on $k$. As the case $k=0$ is obvious, we start with
the case $k=1$, so that  
$$S_{1,A}(v)=\|v\|+\| H v\|.$$

For $a\in\Omega_i$ we obtain  
\begin{eqnarray*}
H\pi(k_0)\pi(a)v_0&=&\pi(k_0)(E+F)\pi(a)v_0
=\pi(k_0)\pi(a)(a^{-\beta}E +a^\beta F)v_0\\ &=&\pi(k_0)\pi(a)\big((a^\beta+a^{-\beta}) F + a^{-\beta}
\underbrace{(E-F)}_{\in \k}\big)v_0\\ &=&\pi(k_0)\pi(a)\big((a^\beta+a^{-\beta}) F \big)v_0\\ &=&\pi(k_0)\big(a^{-\beta}(a^\beta+a^{-\beta})
F
\big)\pi(a)v_0.\end{eqnarray*}
Using the unitarity of $\pi$ we get 
\begin{eqnarray*}
S_{1,A}(\pi(k_0)\pi(a)v_0)&=&\|\pi(a)v_0\| +\|H\pi(k_0)\pi(a)v_0\|\\
&=&\|\pi(a)v_0\| +|a^\beta +a^{-\beta}|\ \|F\pi(a)v_0\|.
\end{eqnarray*}
Since $\|F\pi(a)v_0\|\leq C S_1(\pi(a)v_0)$, 
the proof of the $k=1$ case is complete.

Suppose that the statement holds for $k-1$. We must show that 
$$\|H^k \pi(k_0) \pi(a) v_0\|\leq C \sum_{j=0}^k |a^\beta +a^{-\beta}|^j 
S_j(\pi(a)v_0).$$ 
As before we have 
$$ H^k\pi(k_0)\pi(a)v_0 =\pi(k_0)(E+F)^k\pi(a)v_0
=\pi(k_0)\pi(a)(a^{-\beta}E +a^\beta F)^kv_0. $$
Now we must arrange the expressions 
$(a^{-\beta}E +a^\beta F)^k$ in an appropriate way. With 
$U=E-F$,
$$(a^{-\beta}E +a^\beta F)^k= \big ((a^\beta+a^{-\beta})F +a^{-\beta}U)^k.$$
Now using repeatedly the fact 
$$(\forall X\in {\cal U}(\g_\C)^j)\qquad 
UXv_0=([U,X]+ XU)v_0=[U,X]v_0,$$
with $[U,X]\in {\cal U}(\g_\C)^{j-1}$, we obtain elements $Z_{j,a}\in {\cal 
U}(\g_\C)^j$,  uniformly bounded 
depending on $a$,  such that  
$$(a^{-\beta}E +a^\beta F)^kv_0 =
\sum_{j=0}^k (a^\beta+a^{-\beta})^j Z_{j,a}v_0$$
just as in the $k=1$ case. \enddemo

In (6.6) one can see that the derivative term of $E+F$ has coefficient vanishing precisely at $x=\pm 1$, the eventual poles of $f_\epsilon$.
This is the key to the proof of Theorem 6.7. For the other rank one groups the singularity of $f_\epsilon$ lies on a hypersurface, such as $\|x\|
= \pm 1$, and this cannot be dominated by a single operator $E+F$, nevertheless an argument of this type will be developed.

\demo{$\oline N$ vector fields from principal series representations}
Before we can complete our discussions for the rank 
one case, we first have to provide some simple facts
related to the Bruhat decomposition. The facts collected below 
hold for an arbitrary semisimple Lie group. 

For every $k\in K$ we write $\lambda_k\: K/M\to K/M,\  xM\mapsto k xM$ for the 
left 
translation on $K/M$. There is the standard action of $G$ on $K/M$ by  
$$G\times K/M\to K/M, \ \ (g, kM)\mapsto \kappa(gk)M.$$
Then every $X\in \g$ defines a vector field $\tilde X$ on $K/M$ via 
$$\tilde X_m ={d\over dt}\Big|_{t=0} \kappa(\exp(tX)m)$$
for all $m\in M$. Write $p_\k\:\ \ \g\to \k$ for the projection 
along $\a+\n$ and $p_{\k/\m}\: \g\to \k/\m$ for the composition 
of $p_\k$ and the quotient mapping $\k \to \k/\m$. 
\enddemo

 \proclaim{Lemma}  For all $X\in \g$ and $m=kM\in K/M${\rm ,}
$$\tilde X_m=d\lambda_k(\1)p_{\k/\m}(Ad(k^{-1})X).$$
\endproclaim

\demo{Proof} By definition we have that 
$$\tilde X_m={d\over dt}\Big|_{t=0} \kappa(\exp(tX)k)M=
{d\over dt}\Big|_{t=0} \kappa(k\exp(t\Ad(k^{-1})X))M.$$
Set $Y =\Ad(k^{-1})X$. Then there exist  smooth curves 
$Y_\k(t)\in\k$, $Y_\a(t)\in\a$, $Y_\n(t)\in \n$ such that 
$$\exp(tY)=\exp(Y_\k(t))\exp(Y_\a(t))\exp(Y_\n(t))$$
for all $t\in\R$. Differentiation at $t=0$ yields   
$$Y=Y_\k'(0)+Y_\a'(0)+Y_\n'(0)$$
and so $Y_\k'(0)=p_\k(Y)$. Hence we get that 
\begin{eqnarray*}
\tilde X_m&=&{d\over dt}\Big|_{t=0} \kappa(k\exp(Y_\k(t))
\exp(Y_\a(t))\exp(Y_\n(t)))M\\ &=&d\lambda_k(\1){d\over dt}\Big|_{t=0} \kappa(\exp(Y_\k(t)))M\\
&=&d\lambda_k(\1)p_{\k/\m}(\Ad(k^{-1})X).\\
\noalign{\vskip-36pt}
\end{eqnarray*}
 \enddemo
\vglue12pt

The Bruhat decomposition of $G$
$$G=\bigcup_{w\in {\cal W}} \oline N w MAN$$
gives the familiar decomposition of the flag manifold into Schubert cells 
$$K/M=\bigcup_{w\in {\cal W}} \kappa(\oline N w M)/M.$$

 \proclaim{Lemma} For all $m\in \kappa(\oline N)M/M\subeq K/M${\rm ,}
$$T_m K/M=\{ \tilde X_m\: X\in \oline \n\}.$$
\endproclaim

\demo{Proof} This result is known but we include it for completeness. In view of 
Lemma  6.11, it suffices to show that for 
$k=\oline n man$, $\oline n\in \oline N, m\in M, a\in A,\break n\in N$,
$$p_\k(\Ad(k^{-1})\oline \n) +\m =\k$$
or equivalently 
$$\Ad(k^{-1})\oline \n+ \m +\a +\n =\g.$$
By the special choice of $k$,

\begin{eqnarray*}
\noalign{\vskip-12pt}
\Ad(k^{-1})\oline \n+ \m +\a +\n&=&\Ad(man)^{-1}\big(\Ad(\oline 
n^{-1})\oline\n +\Ad(man)(\m+\a+\n)\big)\\ &=&\Ad(man)^{-1}(\oline\n +\m+\a+\n)=\g.\\
\noalign{\vskip-36pt}
\end{eqnarray*}
  \enddemo

\vglue9pt

 \proclaim{Lemma}   Let $(\pi_\lambda, {\cal H}_\lambda)$ be a spherical principal 
series representation 
of $G$ realized in $L^2(K/M)$. Then for all $m\in \kappa(\oline N) M/M\subeq 
K/M${\rm ,}
$$\span_\C\{ \{d\pi_\lambda(X)_m\: X\in \oline \n\}\cup\{\1\}\}\supeq T_m 
(K/M).$$
\endproclaim

\demo{Proof} Recall that the $G$-action on ${\cal H}_\lambda$ is given by 
$$(\pi_\lambda(g)f)(kM)=f(\kappa(g^{-1}k)M) a(g^{-1}k)^{\lambda-\rho}.$$ 
Hence we get for all $X\in\g$ that 
\begin{eqnarray*}
(d\pi_\lambda(X)f)(kM)&=&{d\over dt}\Big|_{t=0} f(\kappa(\exp(-tX)k)M) 
a(\exp(-tX)k)^{\lambda-\rho}\\ &=&-(\tilde Xf)(kM)+ f(k) (\lambda-\rho)\big(p_\a(\Ad(k^{-1})X)\big)\end{eqnarray*}
with $p_\a\:\ \ \g\to \a$ the projection along $\k+\n$. In view of Lemma~6.12, this concludes the proof of the lemma. 
\enddemo

{\it Truncation at infinity}.
Recall that $G$ is a semisimple Lie group of real rank one and 
$(\pi_\lambda, {\cal H}_\lambda)$, $\lambda\in i\a^*$, a 
unitary principal series representation of $G$ which is realized on 
$L^2(\oline N)$, with $\oline N \cong \R^p\oplus \R^q$.

For $0<\eps\leq 1$ let $a_\eps\in \Omega_i$ be as before. 
Define a function on $\R^p\oplus \R^q$ by  
$$f_\eps(X,Y)={1\over [(1+ e^{i({\pi\over 2}-\eps)} \|X\|^2)^2+
e^{i(\pi-2\eps)}\|Y\|^2]^{{p+2q\over 4}-{\lambda\over 2}}}\.$$

Then, in the noncompact realization of $\pi_\lambda$ on $L^2(\oline N)$, the 
vector
$\pi_\lambda(a_\eps)v_0$ is of the form
$c(\lambda,\eps)f_\eps$ where 
$c(\lambda,\eps)$ is a constant depending 
only on $\lambda$ and $\eps$ and uniformly bounded in $\eps$ (cf. \S 5).

In order to estimate $S_k(\pi_\lambda(a_\eps)v_0)$ first we show that 
the behaviour at infinity does not contribute to the singularity as $\eps\to 0$.

 \proclaim{Lemma} Let $(\pi_\lambda, {\cal H}_\lambda)${\rm ,} $\lambda\in i\a^*${\rm ,} 
be a unitary spherical principal series representation of a semisimple 
Lie group $G$ of real rank one. Let $k\in \N$. 
Then there exists  a constant $C>0$ depending on $\lambda$ and $k$ such that
$$(\forall a\in \Omega_i)\qquad S_k (\pi_\lambda(a)v_0)\leq C S_{k,\oline 
N}(\tau \pi_\lambda(a)v_0)$$
for a function $\tau\in C_c^\infty(\oline N)$ with support in 
$B\:=\{Z\in \oline N\: \|Z\|\leq 2\}$. 
\endproclaim

\demo{Proof} To begin we use the compact realization of ${\cal H}_\lambda$ as 
$L^2(K/M)$. Then 
${\cal H}_\lambda^\infty=C^\infty(K/M)$ and the topology 
on ${\cal H}_\lambda^\infty$ is induced from the usual Sobolev norms 
$(S_{k,K})_{k\geq 0}$ and, in addition, each $S_k$ is equivalent 
to $S_{k,K}$.  

The Bruhat decomposition gives 
$$K/M=\kappa(\oline N)M/ M \cup\{ w\}$$
where ${\cal W}_{\a}=\{\1, w\}$. Recall 
that 
$$(\pi_\lambda(a)v_0)(k)=a(a^{-1}k)^{\lambda-\rho}.$$
Let $\phi\in C^\infty(K/M)$ with 
$\phi\equiv 1$ in a small neighborhood of $w$. 
We claim that we can make $\supp \phi$ small enough so that 
$$S_{k, K}(\phi \pi_\lambda(a)v_0)\leq C$$
for all $a\in \Omega_i$. In fact, we have 
$(\pi_\lambda(a)v_0)(w)=e^{(\lambda-\rho)(\log wa^{-1}w)}$
from which our claim easily follows.

In order to estimate $S_k(\pi_\lambda(a)v_0)$ we have 
just seen that we can truncate the function $\pi_\lambda(a)v_0$ away 
from infinity. In particular, we claim that there exists a constant $C>0$
such that 
$$S_k(\pi_\lambda(a)v_0)\leq C S_{k,\oline N}((1-\phi)\pi_\lambda(a)v_0)$$
for all $a\in \Omega_i$. 
In fact it follows from Lemma 6.13 that, uniformly 
for all $a\in \Omega_i$,
\begin{eqnarray*}
S_{k,K}(\pi_\lambda(a)v_0)&\leq &S_{k,K}((1-\phi)\pi_\lambda(a)v_0)
+S_{k,K}(\phi\pi_\lambda(a)v_0)\\ & \leq &C S_{k,\oline N}((1-\phi)\pi_\lambda(a)v_0) +C .\end{eqnarray*}

Now set $\tau= 1-\phi$ so that $\tau$ is in $C_c^\infty(\oline N)$ and we can arrange the support so that $\supp\tau\subeq B$. 
\enddemo

{\it Local estimates for invariant Sobolev norms}.
As before $G$ denotes a semisimple Lie group of real rank one and 
$(\pi_\lambda, {\cal H}_\lambda)$, $\lambda\in i\a^*$, a 
unitary principal series representation of $G$ realized 
on $L^2(\oline N)$.

Recall that $\a=\R A_\alpha $ and define for $t>0$ the elements 
$b_t=\exp(\log t A_\alpha)\break\in A$. 
Then for $f\in L^2(\oline N)$,
$$(\pi_\lambda(b_t)f)(X,Y)= t^{{p+2q\over 2}-\lambda} f(tX, t^2Y) $$
for all $t>0$ and $(X,Y)\in \R^p\oplus \R^q$. 
Also note the action of $\oline N$ on $L^2(\oline N)$: 
$$(\pi_\lambda(Y)f)(X)=f(X-Y)$$
for all $X,Y\in \oline N$.
We also will use the notation $u=\|X\|$ and $v=\|Y\|$ for $(X,Y)\in \R^p\oplus 
\R^q$. 

Our goal is now to obtain estimates for $S_k^G(f)$ for functions with support 
in the ball $B=\{ Z\in\oline N\: \|Z\|\leq 2\}$. Especially we are interested 
in estimating $S_k^G (\tau f_\eps)$. 
Now for $\eps\to 0$ the singularity of $f_\eps$ lies on the sphere 
$u=1$ for $q=0$  while for $q>0$ the singularity lies on 
$u=0$ and $v=1$. This makes it necessary to distinguish the cases 
$q=0$ and $q>0$.

For a multi-index $\gamma=(\gamma_1, \ldots, \gamma_{p+q})\in \N_0^{p+q}$ we set 
$|\gamma|=\gamma_1+\ldots+\gamma_{p+q}$ and define the differential operator 
$$\partial^\gamma={\partial^{\gamma_1} \over\partial 
X_1^{\gamma_1}}\cdot \cdot {\partial^{\gamma_p} \over\partial X_p^{\gamma_p}}\cdot 
{\partial^{\gamma_{p+1}} \over\partial Y_1^{\gamma_{p+1}}}\cdot \cdot {\partial^{\gamma_{p+q}} 
\over\partial Y_q^{\gamma_{p+q}}}.$$

\demo{The case of $q=0$} $\oline N=\R^p$ in this case
and the singularities of $f_\eps$ for $\eps\to 0$ lie on the 
sphere $u=1$. 
\enddemo

 \proclaim{Lemma} Suppose that $q=0$ and let $B=\{ X\in \R^p\: \|X\|\leq 2\}$. Then 
for 
every $k\in \N_0$ there exists a constant $C>0$ depending on $k$ and $\lambda$ 
such that 
$$S_k^G(f)\leq C \sum_{|\gamma|\leq k} \|\  |u-1|^{|\gamma|} \partial^\gamma 
f\|$$
for all smooth functions $f$ with support in $B$. 
\endproclaim

\demo{Proof} The proof is very similar to the proof of Proposition 6.6, just more technical and with more notation. We present the details for the case 
$p=2$, the general case simply having more spherical coordinate variables.

Let $\tau_1, \ldots , \tau_n$ be smooth nonnegative functions with 
$\sum_{i=1}^n \tau_i=\1$ on $B$. Then for any smooth function $f$ with support 
in 
$B$ the definition of the invariant Sobolev norms implies that 
$$S_k^G (f)\leq \sum_{i=1}^n S_k(\pi_\lambda(g_i)(\tau_i f)) \leqno(6.18)$$
for any choice of $g_1, \ldots, g_n\in G$.

Recall that for all smooth functions $h$ with support in $16B$ we have
$$S_k(h)\leq C S_{k, \oline N} (h)$$
for a constant $C$ depending only on $k$ and $\lambda$. 
If we choose our elements $g_i$ such that $\supp (\pi_\lambda(g_i)\tau_i)\subeq 
16 B$,
we get from (6.18) that 
$$ S_k^G (f)\leq C \sum_{i=1}^n S_{k,\oline N}(\pi_\lambda(g_i)(\tau_i f)). 
\leqno(6.19) $$

Next we make a good choice of functions $\tau_i$ and elements $g_i$. 
We start with the 
$\tau_i$.

Let $\phi(x)$ be the one variable function from the proof of Proposition 6.6. 
For any $j\in \N_0$ define a function on $\R^p$ by 
$$\phi_j(X)=\phi(2^j (1-u))$$
and notice that 
$$\supp \phi_j\subeq \{ X\in \R^p\:  2^{-j-1}\leq |u-1|\leq 2^{-j+1}\}\ 
.\leqno(6.20)$$

From now on we use the fact that $p=2$. Then elements $X\in \R^2$
are written in polar coordinates as $X=(u\cos 2\pi\theta, u\sin 2\pi\theta)$ 
with $\theta \in \R$. 
Now choose a nonnegative smooth function $\psi(\theta)$
with support in $[0,2]$ such that 
$$\sum_{m\in \Z}\psi(m+\cdot )=\1.$$
Fix $j\geq 2$. For $l\in\Z$ we define functions on $\R$ by 
$$\psi_{j,l} (\theta)=\psi(2^j \theta-l).$$
Notice that 
$$\supp \psi_{j,l}\subeq \Big[{l\over 2^j}, {l+2\over 2^j}\Big].\leqno(6.21)$$
As $j\geq 2$ these intervals have length at most ${1\over 2}$ and so 
these functions descend to smooth functions on the circle $\R/\Z$.  
In particular, we obtain that 
$$\sum_{l=0}^{2^j} \psi_{j,l}(\theta)=1$$
for all $\theta\in \R/\Z$.

Now we can define our partition of unity. We fix $m\geq 2$ and define for 
all
$2\leq j\leq m-1$ and $0\leq l\leq 2^j$ the functions 
$$\tau_{j,l}(X)=\phi_j(X) \psi_{j,l}(\theta)$$
where $X=(u\cos 2\pi\theta, u \sin 2\pi\theta)$.

Recall the one-variable functions $\tau_m$ from Proposition 6.6. We define 
functions $\tau_{m,l}$, $0\leq l\leq 2^m$ by 
$$\tau_{m,l}(X)=\tau_{m-1}(1-u) \psi_{m-1,l}(\theta)\ .$$

We claim that we have for all $\gamma\in \N_0^p$ with $|\gamma|\leq k$ 
$$\|\partial^\gamma \tau_{j,l}\|_\infty \leq C 2^{|\gamma|j}\leqno(6.22)$$
for a constant $C$ depending only on $k$, $\psi$ and $\phi$. In fact it follows from (6.20) 
that all the $\tau_{j,l}$ are supported away from the origin where the cartesian 
partial derivatives can be dominated by spherical partial derivatives 
${\partial \over \partial u}$ and ${\partial \over \partial \theta}$. The 
claim follows then from 
the construction of the $\tau_{m,j}$.

It also follows from the construction of the $\tau_{j,l}$  that we have for 
all choices of $j,l$ 
that 
$${\rm diam} (\supp \tau_{j,l})\leq 16\cdot  2^{-j}. \leqno(6.23) $$

For a fixed $f$ as in the statement of the lemma we may assume that 
$$\sum_{j=2}^m \sum_{l=0}^{2^j}\tau_{j,l}=\1 \qquad \hbox{on} \qquad  \supp f.$$
Then from (6.19) we get 
$$ S_k^G (f)\leq C \sum_{j=2}^m \sum_{l=0}^{2^j} S_{k,\oline 
N}(\pi_\lambda(g_{j,l})(\tau_{j,l} f)) \leqno(6.24)$$
for elements $g_{j,l}\in G$  to be specified. 

For each $j,l$ we now pick an element $Z_{j,l}\in \supp \tau_{j,l}$. 
Then define 
$$g_{j,l} = (-Z_{j,l})b_{2^{-j}} \in A\oline N \subeq G$$
and notice that 
$$(\pi_\lambda(g_{j,l})h)(X)=   2^{-{jp\over 2} +\lambda} h(2^{-j} 
X+Z_{j,l})\leqno(6.25)$$
for all functions $h$ and $X\in \R^p$.  Then (6.20), (6.21), (6.23) and (6.25) imply that 

$$\supp [\pi_\lambda(g_{j,l}) \tau_{j,l}] \subeq  16 B.$$

We choose $m$ large enough such that 
$$\sum_{j=0}^{2^m}  S_{k,\oline N} (\pi_\lambda(g_{m,j})\tau_{m,j} f) \leq 
\|f\|$$
holds. This is possible in view of (6.22) and (6.25). Also from (6.22) we can obtain  
\begin{eqnarray*}
 S_k^G(f)&&\hskip-.25in\leq \|f\|\\
 &&\hskip-.25in + C \sum_{|\gamma|\leq k} 
\sum_{j=0}^{m-1}\sum_{l=0}^{2^j}  \Big(\int_{
2^j(\supp \tau_{j,l} -Z_{j,l})} 
 2^{-2j|\gamma|} 2^{-jp} |(\partial^\gamma f) (2^{-j} X+Z_{j,l})|^2   
dX\Big)^{1\over 2}\\ && \hskip-.25in\leq \|f\| + C \sum_{|\gamma|\leq k} \sum_{j=0}^{m-1}\sum_{l=0}^{2^j}  
\Big(\int_{\supp \tau_{j,l}} 
 2^{-2j|\gamma|} |(\partial^\gamma f) (X)|^2 \ dX\Big)^{1\over 2}.
\end{eqnarray*}
Now on $\supp \tau_{j,l}$ we have $2^{-j}\leq 2|u-1|$ by (6.20)  and so it follows 
that 
\begin{eqnarray}\quad\qquad
S_k^G(f)&\nhs\leq\nhs & \|f\| + C \sum_{|\gamma|\leq k} 
\sum_{j=0}^{m-1}\sum_{l=0}^{2^j}  \Big(\int_{\supp \tau_{j,l}} 
 |u-1|^{2|\gamma|} |(\partial^\gamma f) (X)|^2 \ dX\Big)^{1\over 2}\speqnu{6.26}\\ &\nhs \leq\nhs& C \sum_{|\gamma|\leq k} \|\  |u-1|^{|\gamma|}
\partial^\gamma  f\|,\nonumber
\end{eqnarray}
as was to be shown. \enddemo

\demo{The case of $q>0$} Here we have $\oline N=\R^p\oplus\R^q$ and the 
singularity 
of $f_\eps$ lies on the sphere $u=0$ and $v=1$. The analog of Lemma 6.15 for 
this 
geometry is 
\enddemo

 \proclaim{Lemma} Suppose that $q>0$ and let $B=\{ Z\in \R^p\oplus \R^q\: \|Z\|\leq 
2\}$. Then for 
every $k\in \N_0$ there exists a constant $C>0$ depending on $k$ and $\lambda$ 
such that 
$$S_k^G(f)\leq C \sum_{|\gamma|\leq k} \|\  |u|^{|\gamma|} \partial^\gamma 
f\|\leqno(1)$$
and 
$$S_k^G(f)\leq C \sum_{|\gamma|\leq k} \|\  |v-1|^{|\gamma|} \partial^\gamma 
f\|\leqno(2)$$
for all smooth functions $f$ with support in $B$. 
\endproclaim

\demo{Proof} The proof is essentially the same as the one of Lemma 6.15 and we 
describe only the necessary modifications. 

The partition of unity is now by the truncators
$$\tau_{j,l}(X,Y)=\phi(2^j u) \phi(2^{2j}(1-v))\psi_{2j,l}(\theta)$$
with $\theta$ the spherical variable of $\R^q$. Next, as done earlier leading to (6.25),  choose  
$$g_{j,l}=(-Z_{j,l})b_{2^{-j}} \in A\oline N$$
with $Z_{j,l}\in \R^q$ and in the support of 
$\phi(2^{2j}(1-v))\psi_{2j,l}(\theta)$. 
Then in the last step (6.26) of Lemma 6.15 we can interpret $2^{-j}$ either as $u$ 
or $|1-v|$. 
Accordingly (1) and (2) follow. \enddemo

{\it Sharp estimates for invariant Sobolev norms
for real rank one}.

\proclaim{Theorem}  Let $(\pi_\lambda, {\cal H}_\lambda)${\rm ,} $\lambda\in i\a^*${\rm ,} 
be a unitary spherical principal series representation of a semisimple 
Lie group $G$ of real rank one. Let $k\in \N$. 
Then there exists  a constant $C>0$ depending on $\lambda$ and $k$ such that
$$(\forall a\in \Omega_i)\qquad S_k^G (\pi_\lambda(a)v_0)\leq C
\|\pi_\lambda(a)v_0\|.$$ 
\endproclaim

\demo{Proof} For $m,n\in \N_0$ define a function
on $\oline N =\R^p\times \R^q$ by  
$$f_\eps^{m,n}(X,Y)={1\over [(1+ e^{i({\pi\over 2}-\eps)} \|X\|^2)^2+
e^{i(\pi-2\eps)}\|Y\|^2]^{{m+2n\over 4}-{\lambda\over 2}}}.$$
Notice that $f_\eps=f_\eps^{p,q}$. 

In view of Lemma 6.14, it is sufficient to prove that there exists a
constant $C>0$ such that  

$$S_{k,\oline N}^G  (\tau f_\eps^{p,q}) \leq C \|\pi_\lambda 
(a_\eps)v_0)\|\leqno(6.27)$$
for all $0<\eps\leq 1$. Here $\tau\in C^\infty(\oline N)$ is a fixed function 
with 
support in the ball $\|Z\|\leq 2$, $Z\in \oline N$.

In order to keep track of  what differentiation does to the function 
$f_\eps^{p,q}$ 
we define a shift operator ${\cal L}$ by 
$${\cal L}(f_\eps^{m,n})=\left\{ \begin{array}{ll}f_\eps^{m,n+2} & \hbox{for $n>0$}\\ f_\eps^{m+2, 0} & \hbox{for $n=0$}.\end{array}\right.$$

Now for $Z\in\oline \n$ the elements $d\pi_\lambda(Z)$ 
are the usual differentiations on $\oline N=\R^p\times \R^q$. 
Consider the action of $d\pi_\lambda(Z)$ on $\tau f_\eps^{m,n}$. When applied to 
$\tau$ it essentially does not increase the function for $\eps\to 0$. 
Applied to $f_\eps^{p,q}$ it increases the exponent ${{p+2q\over 
4}-{\lambda\over 2}}$ by $1$ 
and multiplies with a function $P_\eps(X,Y)$, uniformly 
bounded in $\eps$, $X$ and $Y$, hence can be bounded by a constant. Thus, for 
every $Z\in \oline \n$ 
and $k\in \N_0$ we obtain the inequality 
$$|Z^k\tau f_\eps^{p,q}|\leq C |\tau {\cal L}^{k}(f_\eps^{p,q})|\leqno(2)$$
for all $0<\eps\leq 1$ and a constant $C>0$ independent of $\eps$.

We now distinguish the cases $q=0$ and $q>0$. 
We start with the $q=0$ case. In view of (6.27), (6.28) and Lemma 6.15 it 
is sufficient to show that 
$$||\  \tau |u-1|^k f_\eps^{p+2k,0}\| \leq C \|\tau   f_\eps\|$$
as $\|\tau   f_\eps\|\asymp \|\pi_\lambda(a_\eps)v_0\|$. 
Now this estimate is proved as in Lemma 5.5.

Finally the $q>0$ case is proved similarly by employing 
(6.27), (6.28), Lemma 6.16 (1)  and the proof of Lemma 5.6. \enddemo

\numbereddemo{{R}emark} Using the reduction results in Lemma 6.3 and 
Lemma 6.4 one can prove for an arbitrary $K$-finite vector that 
$$(\forall a\in \Omega_i)\qquad S_k^G (\pi_\lambda(a)v)\leq C
\|\pi_\lambda(a)v\|.$$ 
Since we do not need this result
in this paper and since the proof, albeit more complicated, is fundamentally the 
same, we omit the details.  \enddemo

However, calculations in the higher rank case prompt us to state the following 
conjecture.

\demo{Conjecture {\rm C}} Let $(\pi,{\cal H})$ be an irreducible 
unitary representation of $G$ and 
$v$ a $K$-finite vector. Then there 
exists a constant $C>0$ such that 
$$(\forall a\in \Omega_i)\qquad S_k^G(\pi(a)v)\leq C\|\pi(a)v\|.\
$$
\enddemo

\section{Applications to automorphic forms}

For completeness we recall a few notions from automorphic forms. We denote by 
$\|g\|$, $g\in G$, the operator norm of $\Ad g$, and we write 
${\cal Z}(\g)$ for the center of the universal enveloping algebra ${\cal U}(\g)$ 
of $\g$. 
 
\advance\theoremcount by 1

\demo{Definition {\rm 7.1. (cf.\ [Bo97, Ch.\ V])}} Let $\Gamma<G$ be a discrete subgroup of 
co-finite volume. 
A smooth function $f\: G\to \C$ is called an {\it automorphic form} for 
$\Gamma$ if the following conditions are satisfied: 
\hiha
 {(Aut1)} $f$ is left $\Gamma$-invariant, i.e., $f(\gamma g)=f(g)$ for all 
$\gamma\in \Gamma$, $g\in G$. 

\hiha
{(Aut2)} $f$ is right $K$-finite, i.e., $\span_\C \{ f(\cdot k)\: k\in K\}$ 
is a finite dimensional subspace in $C^\infty (G)$. 

\hiha
{(Aut3)} $f$ is ${\cal Z}(\g)$-finite, i.e., ${\cal Z}(\g)f$ is a 
finite-dimensional subspace of $C^\infty(G)$. 
\hiha
 {(Aut4)} $f$ is of polynomial growth, i.e., there exists an $n\in \N$ and a 
$C>0$ such that 
$$|f(g)|\leq C\|g\|^n$$
for all $g$ in a Siegel set  ${\cal S}\subeq  G$ for the group $\Gamma$.  \enddemo
 
If $(\pi, {\cal H})$ is a unitary representation of $G$, then we write ${\cal 
H}^\infty$
for the Fr\'echet submodule of smooth vectors. The space of distribution vectors 
${\cal H}^{-\infty}$
is by definition the strong antidual of ${\cal H}^\infty$. If $\Gamma<G$ is a 
subgroup, then we write 
$({\cal H}^{-\infty})^\Gamma$ for the $\Gamma$-invariants of ${\cal 
H}^{-\infty}$.

The following proposition is well known but, because it is crucial to the 
approach used, we include its short proof. 

 \proclaim{Proposition} \hskip-8pt Let $(\pi, {\cal H})$ be an irreducible unitary representation 
of~$G${\rm ,} 
$\eta\in ({\cal H}^{-\infty})^\Gamma$ and $v\in {\cal H}$ a $K$\/{\rm -}\/finite vector. 
Then the function 
$$\theta_{v,\eta}\: G\to\C, \ \ g\mapsto \la 
\pi(g)v,\eta\ra=\oline{\eta(\pi(g)v)}$$
is an automorphic form. 
\endproclaim

\demo{Proof} Since $K$-finite vectors are analytic, the function $\theta_{v,\eta}$ is 
defined. As  $\eta$
 is $\Gamma$-invariant,  $({\rm Aut}1)$ follows,  while $({\rm Aut}2)$ is a 
consequence of the\break 
$K$-finiteness of $v$. Since $(\pi, {\cal H})$ is irreducible, $({\rm Aut}3)$ is 
a consequence 
of Schur's Lemma. Finally, the fact that ${\cal H}^\infty$ has moderate growth 
(cf.\ [Wal92, 11.5.1])
implies $({\rm Aut}4)$.  \enddemo

\numbereddemo{{R}emark}   There is also a very useful converse   to Proposition 7.2, 
i.e., 
every automorphic form is a generalized matrix coefficient (cf.\ [Wal92, 
11.9.2]). \enddemo 

\proclaim{Proposition}   Let $G$ be a semisimple Lie group 
with $G\subeq G_\C$. Let $\Gamma<G$ be a 
co\/{\rm -}\/compact subgroup. Let $(\pi_\lambda, {\cal H}_\lambda)${\rm ,} $\lambda\in i\a^*${\rm ,} 
be a unitary spherical principal series representation. Let 
$\eta\in ({\cal H}_\lambda^{-\infty})^\Gamma$ define an embedding 
$${\cal H}_\lambda^\infty\to C^\infty(\Gamma\bs G), \ \ v\mapsto 
\theta_{v,\eta}; 
\ \theta_{v,\eta}(\Gamma g) =\la \pi_\lambda(g)v, \eta\ra.$$
Then for all $k> {1\over 2} \dim G$ there exists a constant $C>0$ such that 
$$(\forall v\in {\cal H}_\lambda^\infty)\qquad 
\|\theta_{v,\eta}\|_\infty\leq C S_k(v).$$
 In particular{\rm ,} since $\|\cdot\|_\infty$ is $G$\/{\rm -}\/invariant  
$$(\forall v\in {\cal H}_\lambda^\infty)\qquad 
\|\theta_{v,\eta}\|_\infty\leq C S_k^G(v).$$
\endproclaim

\demo{Proof} This is the content of [BeRe99, Lemma 3.3 and Prop.\ B.2]. \enddemo

Combining Theorem 6.17, Proposition 7.4 and Remark 6.18 we
obtain the following 
$L^\infty$ estimate on automorphic forms.

\proclaim{Theorem} Let $(\pi_\lambda, {\cal H}_\lambda)${\rm ,} $\lambda\in i\a^*${\rm ,}  be a 
unitary spherical principal series 
representation of a semisimple Lie group $G$ of real rank one. Let 
$\eta\in ({\cal H}_\lambda^{-\infty})^\Gamma$. For any 
$K$\/{\rm -}\/finite 
vector $v\in {\cal H}_{\lambda, K}$ there is a constant $C$ such that for all 
$a\in \Omega_i$  
$$\|\theta_{\pi_\lambda(a)v,\eta}\|_\infty\leq C \|\pi_\lambda(a)v\|.$$
 \endproclaim

{\it Triple products of Maa{\ss} forms for real rank one}.
Let $\Gamma$ be a co-compact discrete 
subgroup of $G$ and set $Y =\Gamma\bs G$ and $X =\Gamma\bs G/K$. 
For $\lambda\in \a_\C^*$ we have defined the $K$-spherical principal series 
representation 
$(\pi_\lambda, {\cal D}_\lambda)$.
If $\pi_\lambda$ is unitarizable, then ${\cal K}_\lambda$ denotes the 
Hilbert completion of ${\cal D}_\lambda$ in the compact realization. Denote by 
$\hat G$ the unitary dual of $G$ and by $\hat G_s\subset \hat G$ the subset 
corresponding to the $K$-spherical representations, i.e., corresponding 
to the unitarizable $K$-spherical principal series. It is then convenient to
consider $\hat G_s$ as a subset of $\a_\C^*$ by identifying the equivalence 
class of $\pi_\lambda$ with $\lambda$. 

For $\Gamma<G$ co-compact the Plancherel theorem for the right regular 
action 
of $G$ on $L^2(\Gamma\bs G)$ says 
$$L^2(\Gamma\bs G)\cong \hat\bigoplus_{\pi \in\hat G} m_\pi {\cal 
K}_\pi.\leqno(7.1)$$
Here $m_\pi= \dim ({\cal K}_\pi^{-\infty})^\Gamma<\infty$ is the multiplicity of 
$(\pi, {\cal K}_\pi)$ in 
$L^2(\Gamma\bs G)$. If $0\neq \eta\in ({\cal K}_\pi^{-\infty})^\Gamma$, then 
the $G$-equivariant map 
$${\cal K}_\pi^\infty\to C^\infty(\Gamma\bs G), \ \ v\mapsto (\Gamma g\mapsto 
\la \pi(g)v,\eta\ra)$$
extends (up to multiplication by a scalar) to an isometry ${\cal K}_\pi\to 
L^2(\Gamma\bs G)$. 

Write ${\cal K}_\pi^K$ for the subspace of $K$-fixed elements and recall 
that 
$\dim {\cal K}_\pi^K=1$ for $\pi\in \hat G_s$ and zero otherwise. 
Taking $K$-fixed vectors  in (7.1)
we obtain that 
$$L^2(\Gamma\bs G/K)\cong \hat\bigoplus_{\pi \in\hat G_s} m_\pi {\cal 
K}_\pi^K.\leqno(7.2)$$
We will identify $L^2(\Gamma\bs G/K)$ with  a subspace of $L^2(\Gamma\bs G)$. 

If $v_0\in {\cal K}_\pi^K$ and $\eta\in ({\cal K}_\pi^{-\infty})^\Gamma$, 
then  
$$\psi_{v_0,\eta}(\Gamma g K)=\la \pi(g)v_0, \eta\ra$$
defines an element in $C^\infty(\Gamma\bs G/K)$. The function $\psi_{v_0,\eta}$ 
is referred to as a 
{\it Maa{\ss} form}. 

Let $(\psi_n)_{n\in \N}$ be an orthonormal basis of Maa{\ss} forms of $ 
L^2(\Gamma\bs G/K)$. 
Then each $\psi_n(\Gamma  g K)$ is of the form $\la \pi_{\lambda_n}(g)v_0^n, 
\eta\ra$ for a specific 
$\lambda_n\in\hat G_s$, a unit vector $v_0^n\in {\cal K}_{\lambda_n}^K$ and a 
specific 
$\eta\in ({\cal H}_{\lambda_n}^{-\infty})^\Gamma$.

Fix $X_0\in \partial \Omega_i$ and for 
$v\in {\cal H}_{\lambda, K}$ and $0<\eps\leq 1$ set
$$v_\eps =\pi_\lambda(\exp(i(1-\eps)X_0))v.$$

If $\psi=\theta_{v_0,\eta}$ is a Maa{\ss} form, then we write 
$\psi_\eps =\theta_{(v_0)_\eps,\eta}\in L^2(Y)$ for all $0<\eps\leq 1$. 
Now, since $\psi$ is continuous, we have $\psi^2\in L^2(X)\subeq L^2(Y)$. 
Hence we get 
$$\psi^2=\sum_{i\in I} c_i \psi_i \qquad \hbox{with}\qquad c_i=\la \psi^2, 
\psi_i\ra.\leqno(7.3)$$
If one considers (7.3) as an identity in $L^2(Y)$, then analytic continuation 
yields
$$\psi_\eps^2=\sum_{i\in I} c_i \psi_{i,\eps},$$
for all $0<\eps\leq 1$. Taking norms, we get 
$$\|\psi_\eps^2\|^2=\sum_{i\in I} |c_i|^2 \|\psi_{i,\eps}\|^2=\sum_{i\in I} 
|c_i|^2 \|(v_0^i)_\eps\|^2.\leqno(7.4)$$

\proclaim{Theorem} Let 
$G$ be a simple Lie group of real 
rank one and $\Gamma<G$ a co\/{\rm -}\/compact discrete subgroup. Then for every Maa{\ss} 
form $\psi${\rm ,} the coefficients 
$c_i$ of the Fourier series of $\psi^2=\sum_{i\in I} c_i \psi_i$ satisfy the 
following estimates.
\begin{itemize}
\ritem{(i)} If $q=0$, then there exists a constant $C>0$ such that for all
$T>1${\rm ,}
$$\sum_{|\lambda_i|\leq T} |c_i|^2 e^{\pi|\lambda_i|}\leq C \cdot \left\{ \begin{array}{ll} 
T^{2p-2} & \hbox{if $p>1$}, \\ (\log T)^2 & \hbox{if $p=1$}. \end{array}\right.$$  
\ritem{(ii)} If $q>0${\rm ,} then there exists a constant $C>0$ such that for all
$T>1${\rm ,}
$$\sum_{|\lambda_i|\leq T} |c_i|^2 e^{{\pi\over 2}|\lambda_i|}\leq C \left\{ \begin{array}{ll} 
T^{2q-1} & \hbox{if $q>1$}, \\ T (\log T)^2 & \hbox{if $q=1$}. \end{array}\right.$$  
\end{itemize}
\endproclaim

\demo{Proof} We start the proof with the identity (7.4):
$$\|\psi_\eps^2\|^2=\sum_{i\in I} |c_i|^2 \|(v_0^i)_\eps\|^2.$$

Now $\|\psi_\eps^2\|^2\leq \|\psi_\eps\|_\infty^2 
\|\psi_\eps\|^2=\|\psi_\eps\|_\infty^2 \|(v_0)_\eps\|^2$. 
\vglue4pt
(i) If $q=0$, then we have by Theorem 7.5  and Theorem  5.1(ii) for 
$\eps\to 0^+$,
$$\|\psi_\eps^2\|^2\leq C \left\{ \begin{array}{ll} \eps^{2-2p} & \hbox{if $p>1$}\\ |\log\eps|^2 & \hbox{if $p=1$}\end{array}\right.$$ 
for some constant $C$. On the other hand from Proposition  5.7 we get a lower 
bound. Thus
$$\sum_{i\in I} |c_i|^2 e^{\pi|\lambda_i|}e^{-7\eps|\lambda_i|}\leq C  
\left\{ \begin{array}{ll} \eps^{-2p+2} & \hbox{if $p>1$}, \\ |\log \eps|^2 
& \hbox{if $p=1$}. \end{array}\right. $$  

Setting $\eps={1\over T}$ and collecting the $\lambda_i$ with $|\lambda_i|\leq 
T$, the assertion in (i) 
follows. 
\vglue4pt

(ii) By Theorem 7.5, Theorem  5.1(ii) and Proposition  5.7 the 
proof goes as in (i). \enddemo

\numbereddemo{{R}emark}   (a) For $q=0$ and $p=1$ the estimate in Theorem 7.6 (i) is a 
slight improvement of that obtained by  
Bernstein and Reznikov (cf.\ [BeRe99]), viz.\ $(\log T)^3$ compared to our 
$(\log T)^2$.  
\vglue4pt
(b) For $q=0$ and $p=2$ (this corresponds to $G/K\cong \H^3$), 
Sarnak proved in [Sa94] that 
$$|c_i|\leq C (|\lambda_i|^2+1)^{3\over 2} e^{-{\pi\over 2}|\lambda_i|}$$
for all $i$. 
Our estimate in Theorem 7.5(i) yields the slight improvement to 
$$|c_i|\leq C |\lambda_i| e^{-{\pi\over 2}|\lambda_i|}.$$
\vglue4pt
(c) With a more detailed analysis one can improve on the lower estimate 
in Proposition 5.7 in the case of $q=0$. One can show that 
$$\|(v_0)_\eps\|^2=|\phi_\lambda(\exp(-2i(1-\eps)X))|\geq C  e^{(\pi- 
\eps)|\lambda|}\left\{ \begin{array}{ll} |\log\eps| & \hbox{for $p=1$}, \\ \eps^{-p+1} & \hbox{for $ p>1$},\end{array}\right.$$
for all $0<\eps\leq 1$. In particular, this gives a small improvement 
on the estimates of the triple products.
\vglue4pt
(d) We have presented these techniques for co-compact $\Gamma$; 
however, they apply equally well to finite volume but not co-compact lattices. 
We illustrate this in the next section where we obtain estimates on triple 
products of cusp forms.  
 \enddemo

\section{$ G=\Sl(2,\R)$}

{\it Analytic continuation of the discrete series}.
Let $G=\Sl(2,\R)$ and choose $K=\SO(2,\R)$ as a maximal compact subgroup. For 
every $m\in \Z$ 
define a character $\chi_m$ of $K$ by setting 
$$\chi_m\Big(\left( \begin{array}{cccc}\cos \theta & \sin\theta\\ -\sin\theta 
&\cos\theta\end{array}\right)\Big)=e^{{\rm im}\theta} \qquad 
(\theta\in\R).$$
We will identify $\hat K$ with $\Z$ by means of the above isomorphism.

For every $k\in \N$ there exists a unitary highest weight representation of 
$G$ 
with highest weight $-k$ and $K$-weight spectrum 
$-k, -k-2, -k-4, \ldots$. For $k\geq 2$ we obtain the discrete series which can 
be realized in the holomorphic functions on the upper half-plane $X=\{z\in \C\: 
\Im z>0\}$. 
More precisely, a unitary highest weight representation $(\pi_k, {\cal H}_k)$ 
with 
highest weight $-k$, $k\geq 2$, is given by the Hilbert space 

$${\cal H}_k=\{ f\in {\cal O}(X)\: \int_H |f(z)|^2\  {dx \ dy \over 
y^{2-k}}<\infty\}$$   
and the action 

$$(\pi_k(g)f)(z)=(cz+d)^{-k} f\big({az+b\over cz+d}\big) \qquad 
(g^{-1}=\left( \begin{array}{cccc}a &b \\ c &d\end{array}\right)).$$

We recall from Example 4.3 the notation  $A$, 
$A_\C$, and $A_\C^1$. For $\eps>0$ small we define elements $a_\eps\in A_\C^1$ 
by 
$$a_\eps=\left( \begin{array}{cccc}e^{i{\pi\over 4}(1-\eps)} & 0 \\ 0 & e^{-i{\pi\over 
4}(1-\eps)}\end{array}\right).$$
If $v$ is a $K$- weight vector of $(\pi_k, {\cal H}_k)$, then we are interested 
in estimating   
$\|\pi_k(a_\eps)v\|^2$ for $\eps\to 0$. The estimates given thus far have been 
related to  principal series representations. We remark that the method we shall 
follow applies more generally to unitary highest weight modules of other groups.

We have found that estimates for $\|\pi_k(a_\eps)v\|^2$ are obtained more 
easily if
we switch to the realization on the positive real axis. 
For $k\in \N$ we define the Hilbert space:
$${\cal W}_k=L^2(\R^+, x^k{dx\over x})=\{ f\: \R^+\to \C\: \int_0^\infty 
|f(x)|^2 x^k\  {dx \over x}<\infty\}.$$
Then for $k\geq 2$ the mapping 
$$\Phi_k\: {\cal W}_k\to {\cal H}_k,\ \  f\mapsto \Phi_k(f);\ \Phi_k(f)(z)=
\int_0^\infty e^{ixz} f(x) x^k{dx\over x}$$
is, by the Paley-Wiener theorem, up to multiplication by a scalar, an isomorphism 
of Hilbert spaces.
Thus we can transport the $G$-action on ${\cal H}_k$ and obtain a unitary 
$G$-action
on ${\cal W}_k$. Call this representation $(\rho_k, {\cal W}_k)$. 
Of course, one also has a unitary highest weight representation 
$(\rho_1, {\cal W}_1)$ with highest weight~$-1$. 

 We have for all $k\geq 1$
$$(\forall a>0)(\forall f\in {\cal W}_k)\qquad  (\rho_k(\left(\begin{array}{cc} a & 0\\ 0 & 
a^{-1}\end{array}\right))f)(x)
=a^k f(a^2x).$$
For every $m\in \Z$ let ${\cal W}_k^m =\{ v\in {\cal W}_k\: (\forall k\in K)\  
\pi_k(k)v=\chi_m(k)v\}$. 
Then 
$${\cal W}_k=\hat \bigoplus_{j=0}^\infty {\cal W}_k^{-k-2j}$$
is the orthogonal decomposition of ${\cal W}_k$ into $K$-isotypical 
components. Moreover, one can show that 
$${\cal W}_k^{-k-2j}=\C\{ e^{-x} p_j(x)\} $$
with $p_j(x)$ a polynomial of degree $j$.  

For the discrete series of $\Sl(2,\R)$ we come now to their analog of Theorem~5.1.

\proclaim{Theorem}  Let $(\rho_k, {\cal H}_k)$ be a unitary highest weight 
representation 
of $G=\Sl(2,\R)$ with lowest weight $-k\in -\N$. Then for every $K$\/{\rm -}\/finite 
vector 
$v\in {\cal W}_k$ with weight $m\in -k-2\N_0${\rm ,}
$$\|\rho_k(a_\eps)v\|^2\asymp \eps^{-|m|}.$$
\endproclaim

\demo{Proof} Let $m=-k-2j$. Then $v$ is a multiple of $f(x)=e^{-x}p_j(x)$. Hence it 
suffices to show that 
for $g(x)=e^{-x} x^n$ one has 
$$\|\rho_k(a_\eps)g\|^2\asymp \eps^{-2n-k}$$
for every $n\in \N_0$. 

Clearly for fixed $g$  the map
$$F\: A_\C^1\to {\cal W}_k; \  F(a)(x)\ =a^k g(a^2x)$$
is analytic. Since $F\res_A=\rho_k\res_A(\cdot)g$, we have that 
$$(\rho_k(a)g)(x)=a^kg(a^2x)$$
for all $a\in A_\C^1$. Therefore
\begin{eqnarray*}
 \|\rho(a_\eps)g\|^2 &=&\int_0^\infty |e^{-a_\eps^2 x} (a_\eps^2 
x)^n|^2 x^k\ {dx\over x}\\ &=&|a_\eps|^{4n}\int_0^\infty e^{-2\Re (a_\eps^2) x} x^{2n+k}\
{dx\over x}\\ &\asymp& (\Re a_\eps^2)^{-2n-k}= (\Re e^{i{\pi\over 
2}(1-\eps)})^{-2n-k}\asymp\eps^{-2n-k}.\\
\noalign{\vskip-36pt}\end{eqnarray*}
  \enddemo

\vglue12pt
{\it Automorphic forms associated to the discrete series}.
Let  $\Gamma<G$ be a lattice in $G$ but not co-compact. Let $P<G$ be a parabolic 
subgroup of $G$ and $P=MA_PN_P$ with $M=\{\pm I\}$ its Langlands decomposition. 
Call $P$ {\it cuspidal} for $\Gamma$ if $\Gamma\cap N_P\neq \{\1\}$. 

\numbereddemo{Definition} An automorphic form $f\: G\to \C$  is called a {\it cusp form} 
if for all 
cuspidal parabolic subgroups $P<G$,
$$(\forall g\in G) \qquad \int_{(N_P\cap \Gamma)\bs N_P} f(ng) \ dn =0.
$$
\enddemo

Recall from [Bo97, Cor.\ 8.7] that every cusp form $f\: G\to \C$ belongs to 
$L^p(\Gamma\bs G)$ for all $1\leq p\leq \infty$. 

If $(\pi, {\cal H})$ is an irreducible unitary representation of $G$, with space 
of\break $K$-finite vectors ${\cal H}_K$, set 
$$({\cal H}^{-\infty})_c^\Gamma= \{ \eta\in ({\cal H}^{-\infty})^\Gamma\: 
\theta_{v,\eta}
\ \hbox{is a cusp form for all}\  v\in {\cal H}_K\}.$$
Note that if $v\in {\cal H}_K$, $v\neq 0$, then the irreducibility of $(\pi, 
{\cal H})$ implies 
that 
$$({\cal H}^{-\infty})_c^\Gamma=\{ \eta\in ({\cal H}^{-\infty})^\Gamma\: 
\theta_{v,\eta}
\ \hbox{is a cusp form}\}.$$

Let now $\eta\in ({\cal H}^{-\infty})_c^\Gamma$. Then the map 
$${\cal H}_K\to L^2(\Gamma\bs G), \ \ v\mapsto \theta_{v,\eta}$$
gives rise (up to scalar multiple) to an isometric embedding 
$${\cal H}\to L^2(\Gamma\bs G).$$
For $v$ a $K$-weight vector, as usual, we set $v_\eps =\pi(a_\eps)v$ for all 
$\eps>0$ small.

\proclaim{Theorem}  Let $(\pi_k, {\cal H}_k)${\rm ,} $k\geq 2${\rm ,}  be a holomorphic discrete 
series 
representation. Let $\eta\in ({\cal H}_k^{-\infty})_c^\Gamma$ and $v\in {\cal 
H}_k$ a 
$K$\/{\rm -}\/weight vector of weight $m$. Then for $\eps$ small the following assertions 
hold\/{\rm :}
\begin{itemize} 
\ritem{(i)} $\|\theta_{v_\eps, \eta}\|_{L^2(\Gamma\bs G)}\asymp\eps^{-{|m|\over 
2}}${\rm ;}
\ritem{(ii)} $\|\theta_{v_\eps,\eta}\|_\infty\leq C\eps^{-{|m|\over 2}}$ for a 
constant $C$ depending only on $v$. 
\end{itemize}
\endproclaim

\demo{Proof} (i) is just a restatement of Theorem 8.1, since the embedding ${\cal 
H}_k\to L^2(\Gamma\bs G)$
is (up to scalar multiple) isometric, see Proposition 7.4.
\vglue4pt
(ii) First we need a little notation. Let $X_1, X_2, X_3$ be a basis 
of $\g$. The $n^{\rm th}$ Sobolev norm of the representation 
$(\pi_k, {\cal H}_k)$ is given by the equivalent norm
$$(\forall v\in {\cal H}_k^\infty)\qquad S_n(v) =\| (\1 -\sum_{i=1}^3 
d\pi_k(X_i)^2)^{n\over 2}v\|.$$ 
It follows from [BeRe99, Prop.\ 4.1] that 
$$\|\theta_{v,\eta}\|_\infty \leq C S_3(v)  \quad \hbox{for all $v\in {\cal 
H}_k^\infty$}.$$ 
Moreover since $\|\cdot\|_\infty$ is $G$-invariant we obtain 
$$\|\theta_{v,\eta}\|_\infty \leq C S_3^G(v)  \quad \hbox{for all $v\in {\cal 
H}_k^\infty$},$$ 
where $S_k^G(\cdot)=\inf_{g\in G} S_k(\pi_k(g)\cdot)$ is the infimum seminorm 
(cf. Definition 6.1). 

We will work with the realization ${\cal W}_k$ on $\R^+$ of 
$\pi_k$. Then 
$$d\pi_k({\cal U}(\g_\C))\subeq \C [x, {d\over dx}]$$
which can be seen either by direct computation or as a special case of a general 
result on unitary highest weight representations (cf.\ [KrNe00]). Set 
$D^N=(\1 -\sum_{i=1}^3 d\pi_k(X_i)^2)^N$, $N\in \N_0$. It suffices to show that 
$$\inf_{g\in G}\|D^N\pi_k(g)v_\eps\|\leq C\|v_\eps\|$$
for any $K$-weight vector $v$ and a constant $C$ depending only on $v$ and not 
on $\eps$. 
We do this only for a highest weight vector $v$; the computation for the 
other $K$-types is similar.
Write $D^N=\sum_{j,l=0}^N a_{jl}\  x^j {d^l\over dx^l}$ and note that 
$v$ is given by the function $v(x)=e^{-x}$. 
For $t>0$, $\in \R$ we define elements $b_t, n_s\in G$ by 
$$b_t 
=\left( \begin{array}{cccc}\sqrt{t} &\ 0\\ 0 & {1\over \sqrt{t}}\end{array}\right),\quad  \hbox {and}
\quad n_s=\left( \begin{array}{cccc} 1 & s \\ 0 & 1\end{array}\right).$$ 
Then for all $f\in {\cal W}_k$
$$(\pi_k(n_s b_t)f)(x) =e^{-isx} t^{k\over 2} f(tx).$$
Thus,
\begin{eqnarray*}
\big(D^N\pi_k(n_sb_t)v_\eps\big)(x)&
=&\sum_{j,l=0}^N a_{jl}\  x^j {d^l\over dx^l} \big(a_\eps^{k} t^{k\over 2} 
e^{-a_\eps^2 tx}e^{-isx}\big)\\ &=&a_\eps^k t^{k\over 2}
\sum_{j,l=0}^N (-is-a_\eps^2 t )^la_{jl}\ x^{j} e^{-a_\eps^2 
tx}e^{-isx},\end{eqnarray*}
and
\begin{eqnarray*}
\|D^N\pi_k(n_sb_t)v_\eps\|^2&=&t^k\sum_{j,j',l,l'=0}^N a_{jl}  
\oline{a_{j'l'}} (is+a_\eps^2 t)^l 
\oline {(is+a_\eps^2 t)}^{l'}\\
&&\times\ \int_0^\infty
x^{j+j'} e^{-2\Re (a_\eps^2) tx} x^k \  {dx\over x}\\ &=&\sum_{j,j',l,l'=0}^N a_{jl}  \oline{a_{j'l'}} (is+a_\eps^2 t)^l 
\oline {(is+a_\eps^2 t)}^{l'} t^{-(j+j')}\\ &&\times\ (2\Re (a_\eps^2)) ^{-(j+j'+k)}\int_0^\infty
x^{j+j'} e^{-x} x^k \  {dx\over x}\\ &\leq & C \sum_{j,j',l,l'=0}^N |is+a_\eps^2 t|^{l+l'} t^{-(j+j')} 
\eps^{-(j+j'+k)}.\end{eqnarray*}
Taking $t={1\over \eps}$ we get 
$$\|D^N\pi_k(n_sb_t)v_\eps\|^2\leq C \eps^{-k} \sum _{l,l'=0}^N 
|is+{a_\eps^2\over \eps}|^{l+l'}.$$
Finally for $s=-{1\over \eps}$ the expression $|is+{a_\eps^2\over \eps}|$ is 
bounded for all 
$0<\eps\leq 1$ and so we see that $\inf_{g\in G}\|D^N\pi_k(g).v_\eps\|\leq C 
\eps^{-{k\over 2}}$, 
completing  the proof of (ii).  \enddemo

{\it Triple products in the co\/{\rm -}\/compact case}.
 We recall how to relate automorphic forms of weight $m$ to automorphic 
functions on $G/K$.
 An automorphic form $f\: G\to \C$ is called of {\it weight} $m\in\Z$ if 
$$(\forall k\in K) \qquad f(gk)=\chi_m(k)f(g).$$

We identify $X$ with $G/K$ by means of the isomorphism 
$G/K\to X,\break   gK\mapsto g.i$. 
For every $z\in X$ and $g=\left(\begin{array}{cc} a& b\\ c& d\end{array}\right)$ set $\mu(g,z)=cz+d$
and recall that $\mu$ satisfies the cocycle relation $\mu(g_1g_2,z)=\mu(g_1, 
g_2\cdot z) \mu(g_2, z)$.
For $m\in \Z$ we set $\mu_m=\mu^{-m}$ and note that $\mu_m(k,i)=\chi_m (k)$ for 
all
$k\in K$. If $f\: G\to \C$ is of weight $m$, then the function 
$$ F(gK)=\oline{\mu_m(g,i)} f(g)\leqno(8.1)$$
defines an analytic function on $G/K$ which satisfies 
$$(\forall \gamma\in \Gamma)\qquad  F(gK) =\oline{\mu_m(\gamma, 
gK)}^{-1}F(\gamma g K)$$
for all $gK\in G/K$.

We say that $f\:G\to \C$ is an {\it anti\/{\rm -}\/holomorphic automorphic form} if\break $F\: 
G/K\to \C$ 
is an anti-holomorphic function. 
In the usual notation we write $M_k^0(\Gamma)$ for the anti-holomorphic cusp 
forms on $X=G/K$ of weight 
$k\in \N$. 
If $(\pi_k, {\cal H}_k)$ is a unitary highest weight representation, then we 
write 
$v_k$ for a normalized highest weight vector. If $\Theta_{v_k, \eta}$ is 
the function on $G/K$ associated to  $\theta_{v_k, \eta}$ via (8.1), then it can 
be deduced 
with the help of Proposition 7.2  and [Wal92, 11.9.2] that 
the mapping 
$$({\cal H}_k^{-\infty})_c^\Gamma\to M_k^0(\Gamma), \ \ \eta\mapsto \Theta_{v_k, 
\eta}$$
is an isomorphism of (finite dimensional) vector spaces. 

Let $f\: G\to\C$ be an automorphic form of weight $m$. Then $|f|$ factors 
to a function 
on $G/K$ which we also denote by $|f|$. Then (8.1) gives 

$$(\forall z=x+iy\in X)\qquad |f|^2(z)= y^{-m} |F(z)|^2.\leqno(8.2)$$

Set 
$$\quad X_0=\left( \begin{array}{cccc} 1 & 0 \\ 0 & -1\end{array}\right),$$
and define $\rho\in \a^*$ by $\rho(X_0)=1$. In a slightly different way,   we 
will identify $\a_\C^*$ with $\C$ by means of the isomorphism 
$$\C\mapsto \a_\C^*, \ \ z\mapsto z2\rho.$$

Let $(\pi_k, {\cal H}_k)$, $k\in \N_0$, be a unitary highest weight 
representation of 
$G$ and $\theta_{v,\eta}(\Gamma g)=\la \pi_k(g)v,\eta\ra$ an automorphic  form 
of weight $m$. Then  $|\theta_{v,\eta}|^2\in L^2(\Gamma\bs G/ K)$ and as in 
(7.3) we have 
$$|\theta_{v,\eta}|^2=\sum_{n\in \N} c_n \psi_n$$
with 
$$c_n=\la |\theta_{v,\eta}|^2, \psi_n\ra.$$
As before, we are interested in estimating the coefficients $c_n$, the  so 
called 
{\it triple products}. For this we use the strategy introduced by Bernstein and 
Reznikov as used in the proof of Theorem 7.6. Of course, here we must employ the 
estimates just obtained related to discrete series.

\proclaim{Theorem}  Let $\Gamma$ be a discrete co\/{\rm -}\/compact subgroup of $G=\Sl(2,\R)$. 
Let 
$\theta_{v,\eta}$ be an automorphic form of weight $m\in Z$ associated to a 
unitary highest weight representation $(\pi_k, {\cal H}_k)${\rm ,} $k\in \N_0${\rm ,} of 
$G$. 
If $|\theta_{v,\eta}|^2=\sum_{n\in \N} c_n \psi_n$ is the orthogonal expansion 
of $|\theta_{v,\eta}|^2$ in Maa{\ss} forms{\rm ,} then there exists a constant $C>0$ such 
that 
for all $T>0${\rm ,}
$$\sum_{|\lambda_n|\leq T} |c_n|^2 e^{\pi |\lambda_n|}\leq C T^{2|m|}.$$
\endproclaim

\demo{Proof} We start with the identity (7.4),
$$\|\theta_{v_\eps,\eta}^2\|^2=\sum_{n\in \N} |c_n|^2 \|v_{0,\eps}^n\|^2.$$
In view of Theorem 8.3, the left-hand side can be estimated as 
$$\|\theta_{v_\eps,\eta}^2\|^2\leq \|\theta_{v_\eps,\eta}\|_\infty^2 
\|\theta_{v_\eps,\eta}\|^2
\leq C_1 \eps^{-|m|}\eps^{-|m|}=C_1 \eps^{-2|m|}$$
for a constant $C_1>0$ independent of  $\eps>0$. On the other hand, from  
Proposition 5.7,
$$\|v_{0,\eps}^n\|^2 \geq C_2 e^{(\pi -7\eps)|\lambda_n|}$$
for a constant $C_2>0$ independent of $n$ and $\eps$. Thus 

$$\sum_{n\in \N}|c_n|^2 e^{(\pi -7\eps) |\lambda_n|}\leq C \eps^{-2|m|}.$$
Taking $\eps={1\over T}$ and collecting all $c_n$ with $|\lambda_n|\leq T$  
prove 
the theorem. \enddemo

{\it Triple products for the noncocompact case}.
Let $P_1, \ldots,P_N$ be a set of representatives of the $\Gamma$-conjugacy 
classes
of cuspidal parabolic subgroups. Every $P_j$ admits a Levi decomposition 
$P_j=MA_jN_j$. We choose $A_j$ such that ${\rm Lie}(A_j)$ is orthogonal to $\k$ 
with respect to 
the Cartan-Killing form on $\g$ . Write 
$a_j\: N_j A_j K\to A_j$ for the middle projection in the Iwasawa decomposition 
$G=N_jA_jK$.  

For every $1\leq j\leq N$ and $s\in \C$ with $\Re s >1$ one defines the  
{\it Eisenstein series}
by 
$$E_j(\Gamma g, s)\:=\sum_{\gamma\in (\Gamma\cap P_i)\bs\Gamma} 
e^{(1+s)\rho(\log a_j(\gamma g))}\qquad 
(g\in G).$$
 For convenience we summarize the properties of Eisenstein series needed here, 
as well as the familiar structure of the Plancherel theorem for $L^2(\Gamma\bs 
G)$. [Bo97] is a convenient reference. The Eisenstein series admit a meromorphic 
continuation in the variable $s$ to the 
entire complex plane (cf.\ [Bo97, 11.9]). The meromorphic continuation of $E_j$ 
will be also denoted 
by $E_j$.   
We note that $E_j(\cdot, s)$, when defined, is an automorphic form in the sense 
of Definition 7.1
(cf.\ [Bo97, Th.\ 10.4]); $E_j(\cdot, s)$ has no poles on $i\R$ (cf.\ [Bo97, 
Th.\ 11.13]).

We normalize the inner product on ${\cal K}_\lambda$ such that 
the spherical vector $v_0^\lambda(g)=a(g)^{-\lambda-\rho}$, $g\in G$, has norm 
$1$.     
Then the fact that $E_j(\cdot, s)$ is an automorphic form together with [Wal92, 
11.9.2] implies 
the existence of an $\eta_{j,s}\in ({\cal K}_s^{-\infty})^\Gamma$ such that 
$$\la \sigma_s(g)v_0^s, \eta_{j,s}\ra=E_j(\Gamma g, 2s)
\quad (g\in G).\leqno(8.3)$$

Let $(\rho, L^2(\Gamma\bs G))$ denote the right regular representation of 
$G$ on $L^2(\Gamma\bs G)$. 
 We write 
$L^2(\Gamma\bs G)_s$ for the $G$-invariant subspace generated  by $L^2(\Gamma\bs 
G/K)$, i.e., 
$$L^2(\Gamma\bs G)_s=\oline{\span\{\rho(G)L^2(\Gamma\bs G/K)\}}.$$ One has
$$L^2(\Gamma\bs G)_s=L^2(\Gamma\bs G)_{s,d}\oplus L^2(\Gamma\bs 
G)_{s,c}\leqno(8.4)$$
where 
$$L^2(\Gamma\bs G)_{s,d}=\hat\bigoplus_{\pi \in \hat G_s} m_\pi {\cal 
K}_\pi\leqno(8.5)$$
with $m_\pi<\infty $ (cf.\ [Bo97, Th.\ 16.2, Th.\ 16.6]), 
and 
$$L^2(\Gamma\bs G)_{s,c}=\sum_{j=1}^N\int_{i\R}^\oplus {\cal K}_{s,j}\ 
ds\leqno(8.6)$$
(cf.\ [Bo97, Th.\ 17.7]). In (8.6) the module ${\cal K}_{s,j}$ is isometrically 
equivalent 
to ${\cal K}_s$ and ${\cal K}_{s,j}^\infty$ is realized 
as the image of the $G$-equivariant embedding 
$${\cal H}_s^\infty\to C^\infty(\Gamma\bs G), \ \ v\mapsto (\Gamma g\mapsto \la 
\sigma_s(g)v, \eta_{j,s}\ra).
\leqno(8.7)$$   

Let $(\psi_n)_{n\in\N}$ be an orthonormal basis of $L^2(\Gamma\bs 
G)_{s,d}\cap L^2(\Gamma\bs G/K)$
of Maa{\ss} cusp forms. Then $\psi_n(\Gamma g)$ equals $\la 
\sigma_{\lambda_n}(g)v_0^n,\eta\ra $ for $v_0^n$ a
normalized $K$-fixed vector in ${\cal K}_{\lambda_n}$ and some element
$\eta\in ({\cal K}_\lambda^{-\infty})^\Gamma$. 

If $f\in L^2(\Gamma\bs G/K)$, then (8.3)--(8.7) imply that 
$$f=\sum_{n=1}^\infty \la f,\psi_n\ra\psi_n +\sum_{j=1}^N\int_\R \la f, 
E_j(\cdot, 2is)\ra E_j(\cdot, 
2is)\  ds \leqno(8.8)$$   
and 
$$\|f\|^2=\sum_{n=1}^\infty |\la f,\psi_n\ra|^2 +\sum_{j=1}^N\int_\R |\la f, 
E_j(\cdot, 2is)\ra|^2 \ ds.
\leqno(8.9)$$

\proclaim{Theorem} Let $(\pi_k, {\cal H}_k)$ be a unitary highest weight 
representation of $G$ and 
$f=\theta_{v,\eta}$ an associated cusp form of weight $m$. Then $|f|^2\in 
L^2(\Gamma\bs G/K)$ and 
there exists a constant $C>0$ such that 
for all $T>0${\rm ,}
$$\sum_{|\lambda_n|\leq T} |\la |f|^2, \psi_n\ra|^2 e^{\pi|\lambda_n|} 
+\sum_{j=1}^N 
\int_{-T}^T |\la |f|^2, E_j(\cdot, 2is)\ra|^2 e^{\pi |s|} ds \leq C T^{2|m|}.$$ 
\endproclaim

\demo{Proof} It is clear that $|f|$ is right $K$-invariant. Moreover, since 
cusp forms are rapidly decreasing (cf.\ [Bo97, Th.\ 7.5]), we 
have $|f|^2\in L^2(\Gamma\bs G/K)$. Hence (8.8) gives
$$|f|^2=\sum_{n=1}^\infty \la |f|^2,\psi_n\ra\psi_n +\sum_{j=1}^N\int_\R \la 
|f|^2, E_j(\cdot, 2is)\ra E_j(\cdot, 
2is)\  ds.\leqno(8.10)$$   
We analytically continue (8.10) as in the proof of Theorem 8.4. First notice 
that 
$E_j(\cdot, 2s)$ corresponds via the Plancherel theorem to $v_0^s\in {\cal 
K}_s$, our 
unit\break $K$-spherical vector in ${\cal K}_s$ (cf.\ (8.3), (8.7)). So we define 
$E_{j,\eps}(\cdot, 2s)$ as the element corresponding to $v_{0,\eps}^s$. 
Now analytic continuation of (8.10) in $L^2(\Gamma\bs G)_s$ gives 
$$|f_\eps|^2=\sum_{n=1}^\infty \la |f|^2,\psi_n\ra\psi_{n,\eps} +
\sum_{j=1}^N\int_\R \la |f|^2, E_j(\cdot, 2is)\ra E_{j,\eps}(\cdot, 
i 2s)\  ds.\leqno(8.11)$$   
Taking norms in (8.11) we arrive at 
$$\||f_\eps|^2\|^2=\sum_{n=1}^\infty |\la |f|^2,\psi_n\ra|^2 \|v_{0,\eps}^n\|^2 
+
\sum_{j=1}^N\int_\R |\la |f|^2, E_j(\cdot, 2is)\ra|^2 \|v_{0,\eps}^s\|^2\ 
ds.\leqno(8.12)$$   
In the proof of Theorem 8.4 we showed that $\||f_\eps|^2\|^2\leq 
C_1\eps^{-2|m|}$. It follows from  
Proposition 5.7  that there exists a constant $C_2>0$ such that 
$\|v_{0,\eps}^s\|^2 \geq C_2 e^{(\pi -7\eps)|s|}$ and $\|v_{0,\eps}^n\|^2 \geq 
C_2 
e^{(\pi -7\eps)|\lambda_n|}$ for all $s\in i\R $ and $n\in \N$. 
Thus 
$$\sum_{n=1}^\infty |\la |f|^2,\psi_n\ra|^2 e^{(\pi -7\eps)|\lambda_n|}+
\sum_{j=1}^N\int_\R |\la |f|^2, E_j(\cdot, 2is)\ra|^2 e^{(\pi -7\eps)|s|}\ ds
\leq C \eps^{-2|m|}.\leqno(8.13)$$
Setting $\eps={1\over T}$ in (8.13) and collecting the appropriate terms prove 
the\break theorem. \pagebreak \enddemo

\numbereddemo{{R}emark}  In [Go81a,b] Good proved a special case of Theorem 8.5 with  
methods from 
analytic number theory. To be more specific, for\break $k>2$, $k$ even, and 
$f=\theta_{v_k,\eta}$, i.e., $|f|^2(z)=y^k|F(z)|^2$, $z\in X$, for some\break
(anti-)holomorphic automorphic form 
$F$ on the upper half-plane $X$, Good proved (cf.\ [Go81a, Th.\ 1]) the estimate 
given in Theorem 8.5 for such an~$f$. 
A comparison of the proofs shows the effectiveness of the  representation 
theoretic approach.
We should point out that the number theory normalization of the 
Eisenstein series $E_j(\cdot, u)$ used in [Go81a] differs by a change of 
parameters $u=2s-1$
from the representation theory notation used in this paper. \enddemo

{\it Estimating Fourier coefficients of holomorphic cusp forms}.
Let $F\in M_k^0(\Gamma)$ be an anti-holomorphic cusp form on the upper half-plane $X$. Let
$P$ be a cuspidal parabolic subgroup for $\Gamma$. Replacing 
$\Gamma$ with a certain $G$-conjugate we  may assume that 
$$\Gamma\cap N=\{\left( \begin{array}{cccc} 1& n\\ 0 & 1\end{array}\right)\: n\in \Z\}.$$
Then $\oline F$ is a holomorphic cusp form on $X$ and as such admits
a Fourier expansion at infinity 
$$\oline F(z)=\sum_{n=1}^\infty c_n e^{2\pi i nz}.$$
As before we have 
$$\forall z=x+iy\in X\qquad |f|^2(z)= y^k |\oline F(z)|^2$$ 
for some cusp form $f=\theta_{v_k, \eta}$ associated to a unitary highest weight 
representation $(\pi_k, {\cal H}_k)$ of $G$. 

If $h\in C^\infty (\Gamma\bs G)$ is rapidly decreasing, then we  define its 
constant term 
$$h_P\: N\bs G\to \C, \ \ Ng\mapsto \int_{(\Gamma\cap N)\bs N} h(ng)\ dn.$$
Since $|f|^2$ is a cusp form and hence rapidly decreasing, the Rankin-Selberg 
convolution theorem
(cf.\ [Bo97, Prop.\ 10.10]) yields for $s\in \C$ with $\Re s>1$ that 
$$\la |f|^2, E(\cdot, s)\ra_{L^2(\Gamma\bs G)}=\la (|f|^2)_P, e^{(1+s)\rho\log 
a(\cdot)}\ra_{L^2(N\bs G)}
\leqno(8.14)$$
with $E$ the Eisenstein series associated to $P$. 
{}From (8.14) and  a straightforward calculation one gets the familiar 
Rankin-Selberg zeta function (cf.\ [Bu97, p.~71--72])
$$\frac12(4\pi)^{{1\over 2}-k - {s\over 2}} \Gamma(k-\frac12+{s\over 2}) \sum_{n=1}^\infty 
|c_n|^2 n^{{1\over 2}-k - {s\over 2}
} =\la |f|^2, E(\cdot, s)\ra\leqno (8.15)$$
for all $s$ with $\Re s>1$. 
{}From Theorem 8.4 we obtain the estimate
  $$\int_{-T}^T |\la |f|^2, E(\cdot, 2is\ra|^2 e^{\pi |s|}\ ds \leq C T^{2k}. 
\leqno(8.16)$$

It was shown by Good in [Go81a, p.\ 544-547] (see also [Pe95, p.\ 121--122]) that 
a combination of 
(8.15) and (8.16) yields the following result: 

\proclaim{Theorem} Let $\oline F(z)=\sum_{n=1}^\infty c_n e^{2\pi i nz}$ be a 
holomorphic 
cusp form of weight $k\in \N$ with respect to an arbitrary discrete subgroup 
$\Gamma< G$ of co\/{\rm -}\/finite 
volume and the property $\Gamma\cap N=\{\left( \begin{array}{cccc} 1& n\\ 0 & 1\end{array}\right)\: n\in 
\Z\}$. 
Then 
$$|c_n|<< n^{{k\over 2}-{1\over 6}+\eps}$$
for every $\eps>0$. \endproclaim

{\it Fourier coefficients of Maa{\ss} forms}.
We shall give yet another application of holomorphic extension, but with a 
different technique. Here we view\break
Whittaker functions  as eigenfunctions of the invariant differential operators 
on the locally symmetric space and use their holomorphic extension. It is 
classical that Whittaker functions have such extensions but this seems to be a 
new use of it.
To avoid technicalities we suppose that $G=\Sl(2,\R)$ and $\Gamma<G$ is a 
discrete subgroup with co-finite volume. We assume that 
$\Gamma$ admits at least one cusp and that 
$$\Gamma\cap N =\{ \left( \begin{array}{cccc} 1 & nc\\ 0 & 1\end{array}\right)\: n\in \Z\}$$
for some $c>0$. Let $\theta_{v,\eta}$ be a Maa{\ss} form on $\Gamma\bs G/ K$. For 
an element 
$z=x+iy\in \C$ with $y>0$ we define $g_z =\left( \begin{array}{cccc}\sqrt y & {x\over \sqrt y}\\ 0 & {1\over \sqrt y}\end{array}\right)$
and we note that $g_z\cdot i=z$. Then $\theta_{v,\eta}$ admits a partial Fourier 
series 
$$\theta_{v,\eta} (\Gamma g_z)=a_0(y)+\sum_{m\in \Z\atop m\neq 0} a_m \sqrt{y} 
K_s
({2\pi\over c} |m| y) e^{2\pi i {m\over c} x}$$
where $s\in \C$ and 
$$K_s(y)={1\over 2} \int_0^\infty e^{-y (t+{1\over t})/ 2} t^s\  {dt\over t} 
\qquad (y>0)$$
is the $K$-Bessel function.

\proclaim{Theorem} Let $\theta_{v,\eta}$ be a Maa{\ss} cusp form for $G=\Sl(2,\R)$. 
Then the Fourier 
coefficients of $\theta_{v,\eta}$ satisfy 
$$(\forall N\geq 2)\qquad \sum_{|m|\leq N\atop m\neq 0} {|a_m|^2\over m} \leq C 
\log N.$$
\endproclaim

\demo{Proof} For every $m\in \Z$ define a unitary character $\chi_m$ of the circle 
group 
$\Gamma\cap N\bs N$ by
$$\chi_m\Big(\left( \begin{array}{cccc}1 & x \\ 0 & 1\end{array}\right)\Big) =e^{2\pi i {m\over c} x}.$$ 
Then we have for all $m\in \Z$, $m\neq 0$, the identity 
$$a_m \sqrt{y} K_s ({2\pi\over c} |m| y) e^{2\pi i {m\over c} 
x}=\int_{\Gamma\cap N\bs N}\theta_{v,\eta}(ng_z)
\chi_{-m}(n) \ dn \leqno (8.17)$$
for $z=x+iy$ in the upper half-plane. 
Now we holomorphically extend both sides of (8.17). For $0<\eps\leq 1$ recall 
that 
$$a_\eps=\left(\begin{array}{cc}  e^{i{\pi\over 4}(1-\eps)} & 0 \\ 0 &  e^{-i{\pi\over 
4}(1-\eps)}\end{array}\right)\in A_\C^1$$
and $y_\eps =e^{i{\pi\over 2}(1-\eps)}$ and note that $g_{y_\eps}=a_\eps$. Then 
by analytic continuation,
$$a_m \sqrt{y_\eps} K_s ({2\pi\over c} |m| y_\eps) =\int_{\Gamma\cap N\bs 
N}\theta_{v,\eta}(na_\eps)
\chi_{-m}(n) \ dn $$
for all $0<\eps\leq 1$. For every $\eps$ define 
$$ f_\eps\: \Gamma\cap N\bs N\to \C, \ \ n\mapsto \theta_{v,\eta} (na_\eps).$$
It follows from Theorem 7.5 that $\|\theta_{v_\eps, 
\eta} \|_\infty 
\leq C |\log \eps|^{1\over 2}$ and so in particular $f_\eps$ is bounded and 
belongs to 
$L^2(\Gamma\cap N\bs N)$. Therefore we get that 
$$f_\eps(n)=\sum_{m\in\Z} b_m \chi_m(n)$$
and 
$$\sum_{m\in \Z} |b_m|^2 =\|f_\eps\|_2^2\leq \|f_\eps\|_\infty^2 \leq C 
|\log\eps|.$$
Since $b_m=a_m \sqrt{y_\eps} K_s ({2\pi\over c} |m| y_\eps) $ we thus obtain 
that 
$$\sum_{m\in \Z\atop m\neq 0} |a_m|^2 |y_\eps| \big|K_s ({2\pi\over c} |m| 
y_\eps)\big|^2 \leq C |\log \eps|.$$
In particular for all $N\in \N$ we get that 
$$\sum_{|m|\leq N\atop m\neq 0} |a_m|^2 \big|K_s ({2\pi\over c} |m| 
y_\eps)\big|^2 \leq C |\log \eps|.$$
Now choose $\eps={1\over N}$. Then for all $|m|\leq N$, 
$$|m|y_\eps=|m|ie^{-i{2\over N}}\approx |m| i +{2|m|\over N}$$ 
and so by the asymptotic expansions of the Bessel functions there exists $C'>0$
such that 
$$\big|K_s ({2\pi\over c} |m| y_\eps)\big|\geq {C'\over \sqrt{|m|}}$$
for $|m|$ large. This proves the theorem. \enddemo

\numbereddemo{{R}emark}  The Ramanujan conjecture for Maa{\ss} forms says that the 
coefficients 
$|a_n|$ grow more slowly than $n^\eps$ for all $\eps>0$. Comparison with the result 
in Theorem 8.8 
shows that our result is consistent with this conjecture and that our result is 
essentially sharp. A little more care with obtaining an asymptotic estimate of C 
log N, instead of a bound of C log N,  together with a Tauberian theorem for 
logarithmic means give  the equivalence of such a result with the order of the 
pole of the Rankin-Selberg zeta function. We thank Wenzi Luo for an informative 
conversation on this topic.  \enddemo

\section{$ G=\Sl(3,\R)$}

Now we are going to apply our techniques to a group of higher rank, namely 
$G=\Sl(3,\R)$. This group is low dimensional enough for 
explicit computations to be possible, yet it also illustrates that the technique 
 works in higher rank  as well. 
In particular, we will verify part of Conjecture B and give a complete 
answer to 
the 
boundary behaviour of the analytically continued spherical functions in the 
direction of the extremal rays. Finally, with these estimates available we can 
give an application to triple 
products.

Let us briefly summarize the notation for this special case:
$$\a =\{\diag(x_1, x_2, x_3)\: x_i\in \R; \ \sum_{j=1}^3 x_j=0\}$$
and 
$$A =\{\diag(a_1, a_2, a_3)\: a_i>0; \ \prod_{j=1}^3 a_j=1\}.$$ 
The positive  system is $\Sigma^+=\{\eps_1-\eps_2, \eps_1-\eps_3, 
\eps_2-\eps_3\}$ and the associated  simple roots are given by 
$$\Pi=\{\eps_1-\eps_2, \eps_2-\eps_3\}.$$
Here $\omega_1=\eps_1$ and $\omega_2=\eps_1+\eps_2$. 
The Weyl group, ${\cal W}_\a\cong S_3$, in this case acts 
as the permutation group of $\{\eps_1, \eps_2,\eps_3\}$.  

Define $3\times 3$ matrices $E_{ij}$ by $E_{ij}=(\delta_{k-i, l-j})_{k,l}$. 
Note that $\g^{\eps_i-\eps_j}=\R E_{ij}$. Now $\oline \n_1 =\g^{\eps_2-\eps_1}$ 
and 
$\oline\n_2 =\g^{\eps_3-\eps_2}\oplus \g^{\eps_3-\eps_1}$ are 
subalgebras of $\oline \n$ satisfying $\oline \n=\oline \n_1 +\oline \n_2$, with 
$\oline \n_2$  abelian. The map

$$\Phi\: \R^3\to \oline N, \ \ (x,y,z)\mapsto \exp(xE_{21})\exp(y 
E_{32})\exp(zE_{31})=
\left( \begin{array}{cccc} 1 & & \\ x & 1 & \\ z & y & 1\end{array}\right)$$
is a diffeomorphism. We take a choice of Haar measure 
$d\oline n$ on 
$\oline N$ so that its pullback  under $\Phi$ is the product of Lebesgue 
measures $dx\ dy \ dz$.

 \proclaim{Lemma}  For all $$a=\diag(a_1,a_2, a_3)\in A$$ and $$\oline n= 
\exp(xE_{21})\exp(y E_{32})\exp(zE_{31})\in \oline N$$ the following assertions 
hold\/{\rm :} 
\vglue4pt
{\rm\ (i)} $a(a\oline n)^{2\omega_1}=a_1^2 +a_2^2 x^2 +a_3^2 z^2;$
\vglue4pt {\rm (ii)} $a(a\oline n)^{2\omega_2}=a_1^2 a_2^2 +a_1^2
a_3^2 y^2+a_2^2 a_3^2 (z+xy)^2 =a_3^{-2} +a_2^{-2}y^2+a_1^{-2}(z+xy)^2.$
\vglue4pt\endproclaim

\demo{Proof} This is elementary. \enddemo

We can improve Theorem 1.8 (ii) slightly.

 \proclaim{Lemma}   For $G=\Sl(3,\R)${\rm ,}
$$B_\C^1 G\subeq K_\C A_\C^{0,\leq} N_\C.$$
\endproclaim

\demo{Proof} By the usual argument with the Bruhat decomposition it is enough to 
show that $B_\C^1 \oline N\subeq K_\C A_\C^{0,\leq} N_\C$. 
This follows readily from Lemma 9.1 and Lemma 2.1(i). \enddemo

In view of the discussion leading to Theorem 4.2(i),
we can say that all spherical functions $\phi_\lambda$, $\lambda\in \a_\C^*$, 
extend to $K_\C B_\C^0 K_\C$.

Since $\g$ is split, $\rho=\omega_1+\omega_2$. Then Theorem 
4.2(iii)
together with Lemma~9.1 gives the following expression 
for spherical functions on $G=\Sl(3,\R)$, 
valid for all  $a\in \exp i\b^0$ and all 
$\lambda=\lambda_1\omega_1+\lambda_2\omega_2\in i\a^*$ (the only 
reason we 
restrict ourselves to $\lambda$ imaginary is to keep the following formula 
manageable).
$$\phi_\lambda(a)=\int_{\R^3} {dx\ dy \ dz\over\big| (a_1 +a_2 x^2 +a_3 
z^2)^{1-\lambda_1}
(a_1 a_2 +a_1 a_3 y^2+a_2 a_3 (z+xy)^2)^{1-\lambda_2}\big|}\leqno (9.1)$$

 For each $\alpha\in \Sigma$ we write $H_\alpha\in \a$ for the  co-root
of $\alpha$. 
Notice that ${\pi\over 2 }H_\alpha\in \partial \b^0$ for all $\alpha\in \Sigma$.

We consider the radial limits of $\phi_\lambda$ in the directions of 
the co-roots and the fundamental weights. 
We begin with the co-root directions.

 \proclaim{Lemma} Let $\alpha\in \Sigma$ and $H_\alpha$ be the corresponding co\/{\rm -}\/root. 
Then{\rm ,}
for all $\lambda\in \a_\C^*$ for which $(\pi_\lambda, {\cal H}_\lambda)$ is 
unitarizable{\rm ,}
$$\phi_\lambda (\exp(i{\pi\over 2}(1-\eps)H_\alpha))\geq C|\log\eps|$$
for a constant $C=C(\lambda)$.
\endproclaim

\demo{Proof} This is a special case of Corollary 4.6. \enddemo

Next, for the fundamental weights we have
$$ X_{\omega_1}=\diag({2\over 3}, -{1\over 3}, -{1\over 3})\quad \hbox{and}\quad 
X_{\omega_2}=\diag({1\over 3}, {1\over 3}, -{2\over 3})$$
and $\pm \pi X_{\omega_j}\in \partial \b^0$.

For these weights we use the following interesting splitting of the spherical 
functions in extremal directions. 

 \proclaim{Proposition}   Let $G=\Sl(3,\R)$ and $X\in \R X_{\omega_1}\cap \b^0$. Put $a 
=\exp(iX)$. Then{\rm ,}
for all spherical functions $\phi_\lambda${\rm ,} $\lambda\in \a_\C^*${\rm ,}   
\begin{eqnarray*}
\phi_\lambda(a)&=&C \int_{\R^2} {1\over\big| (a_1 +a_2 x^2 
+a_3 z^2)^{1-\lambda_1}
(a_1 a_2 +a_2 a_3 z^2)^{1-\lambda_2}\big|}\\ &&\times\ {1\over (\oline{a_1}+\oline{a_2} x^2+\oline{a_3} z^2)^{\Re\lambda_1} (\oline{ 
a_1}\oline {a_2}+
\oline {a_2} \oline{a_3}z^2)^{\Re\lambda_2}}\ dx \ dz
\end{eqnarray*}
with 
$$C=\int_{\oline N'} a(\oline n')^{-\eps_2+\eps_3} \ d\oline n'$$
and $\oline N' =\exp(\g^{\eps_3-\eps_2})$. 

\demo{Proof} Set $b =\exp(i{1\over 2}X)$. Let $\oline\n' =\g^{\eps_3-\eps_2}$ and 
$\oline \n'' =\g^{\eps_2-\eps_1}\oplus\g^{\eps_3-\eps_1}$
and note that $\oline \n=\oline \n'\oplus \oline \n''$ is a direct sum of 
subalgebras. Hence 
$\oline N=\oline N'\oline N''\cong \oline N'\times \oline N''$ with $\oline 
N'=\exp(\oline\n')$ and 
$\oline N''=\exp(\oline\n'')$. Hence we get from Theorem 4.2(iii) that  
\begin{eqnarray*}
\phi_\lambda(b^2)&=&\int_{\oline N} |a(b\oline n)^{2(\lambda 
-\rho)}|\cdot  
 \oline{a(b\oline n)^{-2\Re\lambda}} \ d\oline n\\ &=& \int_{\oline N'}\int_{\oline N''} |a(b\oline n'\oline n'')^{2(\lambda 
-\rho)}|\cdot  
 \oline{a(b\oline n'\oline n'')^{-2\Re\lambda}} \ d\oline n'\ d\oline n''.\end{eqnarray*}
According to the Iwasawa decomposition,  $\oline n'=k' a'n'$. Since 
$k'$ commutes with~$b$,
\begin{eqnarray*}
\phi_\lambda(b^2)&=& \int_{\oline N'}\int_{\oline N''} 
|a(bk'a'n'\oline n'')^{2(\lambda -\rho)}|\cdot  
 \oline{a(bk'a'n'\oline n'')^{-2\Re\lambda}} \ d\oline n'\ d\oline n''\\ &=& \int_{\oline N'}\int_{\oline N''} |a(ba'n'\oline n''n'^{-1})^{2(\lambda 
-\rho)}|\cdot  
 \oline{a(ba'n'\oline n''(n')^{-1})^{-2\Re\lambda}} \ d\oline n'\ d\oline 
n''.\end{eqnarray*}
Since $N'$ normalizes $\oline N''$ in a unipotent way, 
that 
\begin{eqnarray*}
\phi_\lambda(b^2) &=& \int_{\oline N'}\int_{\oline N''} |a(ba'\oline 
n'')^{2(\lambda -\rho)}|\cdot  
 \oline{a(ba'\oline n'')^{-2\Re\lambda}} \ d\oline n'\ d\oline n''\\ &  = &\int_{\oline N'}|a(\oline n')^{2(\lambda-\rho)}|a(\oline 
n')^{-2\Re\lambda}\\
&&\times\ \int_{\oline N''} 
|a(ba'\oline n''(a')^{-1})^{2(\lambda -\rho)}|\cdot  
 \oline{a(ba'\oline n''(a')^{-1})^{-2\Re\lambda}} \ d\oline n'\ d\oline n''\\ &  =&\int_{\oline N'}a(\oline n')^{-2\rho}\int_{\oline N''} 
|a(ba'\oline n''(a')^{-1})^{2(\lambda -\rho)}|\cdot  
 \oline{a(ba'\oline n''(a')^{-1})^{-2\Re\lambda}} \ d\oline n'\ d\oline 
n''.
\end{eqnarray*}
Now the Jacobian of the map $\oline N''\to \oline N'', \ \oline n''\mapsto 
a'\oline n''(a')^{-1}$
is given by $(a')^\gamma$ with $\gamma=-((\eps_1-\eps_2)+(\eps_1-\eps_3))$. 
Hence,
$$\phi_\lambda(b^2)= \int_{\oline N'}a(\oline n')^{-(\eps_2-\eps_3)}\ d\oline n'
\int_{\oline N''} 
|a(b\oline n'')^{2(\lambda -\rho)}|\cdot  
 \oline{a(b\oline n'')^{-2\Re\lambda}} \ d\oline n''.$$
Finally Lemma  9.1 gives
\begin{eqnarray*}
&&\hskip-.75in \int_{\oline N''} |a(b\oline n'')^{2(\lambda -\rho)}|\cdot  
 \oline{a(b\oline n'')^{-2\Re\lambda}} \ d\oline n''\\
&&=\int_{\R^2} {1\over\big| 
(a_1 +a_2 x^2 +a_3 z^2)^{1-\lambda_1}
(a_1 a_2 +a_2 a_3 z^2)^{1-\lambda_2}\big|}\\ &&\phantom{\int}\times {1\over (\oline{a_1}+\oline{a_2}x^2+\oline
{a_3}z^2)^{\Re\lambda_1}  (\oline{ a_1}\oline {a_2}
+\oline {a_2} \oline{a_3}z^2)^{\Re\lambda_2}} \ dx\ dz,\end{eqnarray*}
concluding the proof of the proposition.  \enddemo

\proclaim{Lemma}   Let $\lambda\in \a_\C$. Then for $j=1,2${\rm ,}
$$|\phi_\lambda(\exp(\pm i\pi(1-\eps)X_{\omega_j}))|\leq C {1\over 
\eps}$$
for a constant $C>0$. Furthermore if $\lambda\in i\a^*${\rm ,} then  
$$|\phi_\lambda(\exp(\pm i\pi(1-\eps)X_{\omega_j}))|\asymp {1\over 
\eps}.$$
\endproclaim

\demo{Proof} Notice that $s_{\eps_1-\eps_3}(X_{\omega_1})=-X_{\omega_2}$ (i.e.
for $\Sl(3,\R)$ $\pi_1$ is contragredient to $\pi_2$) so that 
we have to consider only the radial limits in the direction of $\pm X_{\omega_1}$. 
We restrict ourselves to the case of $-X_{\omega_1}$. Then Proposition 9.4 
gives
\begin{eqnarray*}
&&\hskip-.3in\phi_\lambda  (\exp(i\pi(1-\eps) X_{\omega_1}))\\[4pt]
&&\hskip-.25in=C(\lambda) 
\int_{\R^2} {1\over |(e^{i{2\over 3}\pi(1-\eps)}+e^{ -i{1\over 3}\pi(1-\eps)}x^2 
+ 
e^{- i{1\over 3}
\pi(1-\eps)}z^2)^{1-\lambda_1}|} \\ &&\hskip-.25in\times {1\over |( e^{i{1\over 3}\pi(1-\eps)}+  e^{-i{2\over 
3}\pi(1-\eps)}z^2)^{1-\lambda_2}|
(e^{ -i{2\over 3}\pi(1-\eps)}+e^{i{1\over 3}\pi(1-\eps)}x^2 + e^{i{1\over 3}
\pi(1-\eps)}z^2)^{\Re\lambda_1}}  \\ &&\hskip-.25in\times {1\over (e^{-i{1\over 3}\pi(1-\eps)}+  e^{i{2\over 3}\pi(1-\eps)}z^2)^{\Re 
\lambda_2}}\ dx\ dz.\end{eqnarray*}

Hence we get 
\begin{eqnarray*}&&\hskip-.35in\phi_\lambda(\exp(i \pi(1-\eps) X_{\omega_1}))=C(\lambda) 
\int_{\R^2} {|e^{i{1\over 3}\pi(1-\eps)(\lambda_2-\lambda_1)}| e^{i{2\over 
3}\pi(1-\eps)(\Re \lambda_2-\Re \lambda_1)}
 \over |(e^{ i\pi(1-\eps)}+x^2 +z^2)^{1-\lambda_1}|}\\[4pt] &&\hskip-.25in\times {1\over |( 1+  e^{-i\pi(1-\eps)}z^2)^{1-\lambda_2}|
(e^{ -i\pi(1-\eps)}+x^2 + z^2)^{\Re\lambda_1}( 1+  e^{i\pi(1-\eps)}z^2)^{\Re 
\lambda_2}}\ dx \ dz.\end{eqnarray*}

We see that the singularity of the integral is located at $x^2+z^2=1$ and 
$z^2=1$. 

Let us assume now that $\lambda\in i\a^*$. The upper estimate 
for general $\lambda\in\a_\C^*$ will be proved in the same way. 

So for fixed $\lambda\in i\a*$ we get
\begin{eqnarray*}
\phi_\lambda (\exp(-i \pi(1-\eps) X_{\omega_1}))&\asymp&
\int_0^2\int_0^2  {dx\ dz\over |(e^{ i\pi(1-\eps)}+x^2 +z^2)( 
e^{i\pi(1-\eps)} + z^2)|}\\ &\asymp & \int_0^2\int_0^2  {dx\ dz\over |(-1+i\eps +x^2 +z^2)( -1+ i\eps  
+ z^2)|}\\ &\asymp&  \int_0^2\int_0^2  {dx\ dz\over (\eps +|1-x^2 -z^2|)( \eps +
|1- z^2|)}\\ &\asymp & \int_0^2\int_{-1}^1  {dx\ dz\over (\eps +|x^2 +z^2+2z|)( \eps +
|z|)}\\ &\asymp&  \int_0^2\int_0^1  {dx\ dz\over (\eps +x^2 +z^2+2z)( \eps +
z)}\\ &\asymp&  \int_0^1  {dz\over (\eps +z^2+2z)( \eps +
z)}\\ &\asymp&  \int_0^1  {dz\over (\eps +z)^2}\\ &\asymp&  {1\over \eps}.\\
\noalign{\vskip-36pt}
\end{eqnarray*}
 \enddemo
\vglue12pt

Putting these results together we have

\proclaim{Theorem}  Let $G=\Sl(3,\R)$ and $(\pi_\lambda, {\cal H}_\lambda)$
be a unitarizable principal series representation of $G$. Then for the 
associated spherical functions\/{\rm :}
\begin{itemize}
\ritem{(i)} If $X\in \partial \b^0$, then there exists a constant $C\geq 0$ 
such that 
$$|\phi_\lambda(\exp(\pm i(1-\eps)X))|\leq  C{1\over \eps}; $$ 
\ritem{(ii)} If $X=\pm \pi X_{\omega_j}\in \a${\rm ,} $j=1,2${\rm ,} and $\lambda\in i\a^*${\rm ,} 
then 
$$|\phi_\lambda(\exp(i(1-\eps) X))|\asymp {1\over \eps};$$
\ritem{(iii)} If $X \in \partial \b^0${\rm ,} then for a constant $C>0$ 
$$|\phi_\lambda(\exp(i(1-\eps)X))|\geq C|\log\eps|;$$ 
\ritem{(iv)} For all $\lambda\in i\a^*$ the domain $\b_\C^0=\a+i\b^0$ is 
the maximal connected tube domain $\a+i\omega\subeq \a_\C${\rm ,} $\omega\subeq \a${\rm ,} 
containing $0$ such that $\phi_\lambda\res_A$ extends holomorphically to 
$\exp(\a+i\omega)$. 
\end{itemize}

\endproclaim

\demo{Proof} For the proof it is helpful to keep Figure 1 in mind. 
\vglue4pt
(i) For $X$ an extreme point, the assertion  is Lemma 9.5. The 
general case follows from  Proposition 4.4 (cf.\ Figure 1).
\vglue4pt
(ii) Lemma 9.5.
\vglue4pt
(iii) For a co-root direction this inequality is Lemma 9.3. For  
$X_{\omega_j}$ this follows from (ii). An arbitrary $X \in \partial \b^0$ is a 
convex combination of co-roots and $X_{\omega_j}$ (cf.\ Figure 1) so that the result 
follows from Phragmen-Lindel\"of as in the proof of Proposition 4.4.
\vglue4pt
(iv) In view of (iii), this follows from Proposition 4.4
(cf.\ Figure 1).  \enddemo

\numbereddemo{{R}emark}  Theorem 9.6 gives an almost   complete
description of the boundary behaviour for the spherical functions
for $G=\Sl(3,\R)$. Generally we expect 
the following for generic $\lambda\in \a_\C^*$:
$$|\phi_\lambda(\exp(i(1-\eps)X))|\asymp
\left\{ \begin{array}{ll}{1\over \eps} & \hbox{for $X\in \partial {\b^0}$ extremal} \\ |\log \eps|  & \hbox{for $X\in \partial {\b^0}$ not
extremal} .\end{array}\right.$$
 \enddemo

{\it Lower estimates}.
For $G=\Sl(3,\R)$ the best lower estimates one can get for 
spherical functions are in the 
direction of the fundamental weights. We illustrate only this case.

 \proclaim{Proposition}  For $G=\Sl(3,\R)$ and $X\in \oline{\b^0}$ an extreme point{\rm ,} 
i.e.{\rm ,} 
$X=\pm \pi X_{\omega_j}${\rm ,} there exists a constant $C>0$ such that 
$$|\phi_\lambda(\exp(i(1-\eps)X))|\geq C e^{(1-\eps) \sup_{w\in {\cal W}_\a} 
w.\Im\lambda(-iX)}$$
for all $\lambda\in \a_\C^*$. 

\demo{Proof} We may assume that $X=\pi X_{\omega_1}=\pi\diag({2\over 3}, -{1\over 3}, 
-{1\over 3})$.
In what follows we will also see that the assumption $\lambda\in i\a^*$ is 
justified. Fix $0<\eps\leq 1$ and set 
$a=\diag(a_1,a_2,a_3)=\exp(i(1-\eps)X)$. 
Recall from (9.1) that 
$$\phi_\lambda(a)=\int_{\R^3} {dx\ dy \ dz\over\big| (a_1 +a_2 x^2 +a_3 
z^2)^{1-\lambda_1}
(a_1 a_2 +a_1 a_3 y^2+a_2 a_3 (z+xy)^2)^{1-\lambda_2}\big|}.$$
Hence the fact that 
$a_1, a_2, a_3$ are $\R$-collinear in $\C$ as well as that $a_1 a_2, a_1 a_3, 
a_2 a_3$ are 
$\R$-collinear in $\C$ implies that 
\begin{eqnarray*}
\phi_\lambda(a)&\nhs\geq\nhs& \int_0^{1\over 2} \int_0^{1\over 2} \\
&\nhs\nhs&\times\
\int_0^{1\over 2}
{dx\ dy \ dz\over\big| (a_1 +a_2 x^2 +a_3 z^2)^{1-\lambda_1}
(a_1 a_2 +a_1 a_3 y^2+a_2 a_3 (z+xy)^2)^{1-\lambda_2}\big|}\\ &\nhs\geq\nhs& C a_1^{\lambda_1} (a_1 a_2)^{\lambda_2} =C
e^{\lambda(i(1-\eps)X)}.\end{eqnarray*}
 Now the assertion of the proposition follows from the Weyl group invariance of 
$\phi_\lambda$
in the $\lambda$-variable. \enddemo

{\it Estimates on automorphic forms.} 
Here $\Gamma$ denotes 
a co-compact discrete subgroup of $G=\Sl(3,\R)$ and $(\pi_\lambda, {\cal 
H}_\lambda)$ 
a unitary spherical principal series representation of $G$. 
Recall from Proposition 7.4 that each $\eta\in ({\cal 
H}_\lambda^{-\infty})^\Gamma$ defines an embedding 
$${\cal H}_\lambda^\infty\to C^\infty(\Gamma\bs G ), \ \ v\mapsto 
\theta_{v,\eta}; 
\ \theta_{v,\eta}(\Gamma g) =\la \pi_\lambda(g)v, \eta\ra$$
with 
$$(\forall v\in {\cal H}_\lambda^\infty)\qquad 
\|\theta_{v,\eta}\|_\infty\leq C S_k^G(v)$$
for $k> 4={1\over 2}\dim G$. 
As before $v_0$ denotes the normalized $K$-spherical vector of $(\pi_\lambda, 
{\cal H}_\lambda)$. 
Fix an extremal element $X\in\oline{\b^1}$, i.e., $X=\pm {\pi\over 
2}X_{\omega_j}$, and set 

$$(v_0)_\eps =\pi_\lambda(\exp(i(1-\eps)X))v_0$$
for all $0<\eps\leq 1$. As in Section~6 we are interested in the $G$-invariant 
Sobolev norms $S_k^G((v_0)_\eps)$ for $\eps\to 0$.

 \proclaim{Proposition} Let $(\pi_\lambda, {\cal H}_\lambda)${\rm ,} 
$\lambda\in i\a^*${\rm ,}  be a unitary principal 
series representation of $G=\Sl(3,\R)$ and $\eta\in ({\cal 
H}_\lambda^{-\infty})^\Gamma$ for a co\/{\rm -}\/compact 
discrete subgroup $\Gamma<G$. Then there exists a constant $C>0$ such that 
$$\|\theta_{(v_0)_\eps, \eta}\|_\infty\leq C \eps^{-{11\over 2}}
$$
for all $0<\eps\leq 1$. 
\endproclaim

\demo{Proof} We may assume that the extreme point $X\in\oline{\b^1}$ is 
$-{\pi\over 2}X_{\omega_1}$. Since $\dim G=8$, it follows 
from Proposition 7.4 that it suffices to estimate $S_5^G((v_0)_\eps)$. We shall show that 
$$S_5((v_0)_\eps)\leq C  \eps^{-{11\over 2}}$$
for all $0< \eps\leq 1$.

More generally, we will consider $S_k( (v_0)_\eps)$
for an arbitrary $k\in \N_0$. Let $X_1,\cdots,X_8$ be 
a basis of $\g$. Then 
$$S_k(v)=\sum_{m_1+\cdots+m_8\leq k} \| X_1^{m_1}\cdots X_8^{m_8} v\|$$
for all smooth vectors $v$. 

We will work with the noncompact realization of $\pi_\lambda$, i.e., 
${\cal H}_\lambda=L^2(\oline N)$. 
Then the vector $(v_0)_\eps$ is of the form
$c(\lambda, \eps)F_\eps$, where 
$$F_\eps(x,y,z) ={1\over (e^{i(\pi-\eps)}+x^2 +z^2)^{{1-\lambda_1}\over 2}
(1 +y^2+ e^{-i(\pi-\eps)} (z+xy)^2)^{{1-\lambda_2}\over 2}}$$
and $c(\lambda,\eps)$ is a constant uniformly bounded in $\eps$.

We will use results from the proof of Proposition  9.4. 
We denote the change 
of variables in Proposition  9.4 that yielded the splitting by 
$(x,y,z)\mapsto \phi(x,y,z)$. Recall that $\phi$ was given by a composition 
of two maps  $(x,y,z)\mapsto (x +yz,y,z)$ followed by $(x,y,z)\mapsto 
(g(y)x,y, g(y)^{-1}z)$ where $g(y)=1+y^2$. 
Then 
$$U\: L^2(\oline N)\to L^2(\oline N), \ \ f\mapsto 
f\circ \phi$$
defines a unitary operator. Set 
$$ f_\eps(x,y,z) ={1\over (e^{i(\pi-\eps)}+x^2 +z^2)^{{1-\lambda_1}\over 2}
(1 + e^{-i(\pi-\eps)} z^2)^{{1-\lambda_2}\over 2}} 
\cdot g(y)^{-{1\over 2}+
(\lambda_1-\lambda_2)}\ .$$
Then the change of variable formula in Proposition  9.4 implies that
$$ U(F_\eps)(x,y,z)=f_\eps(x,y,z).$$

Define a differential operator $U_j =U\circ X_j\circ U^{-1}$
for $j=1,\cdots,8$. For every 
$k$ we define a seminorm $N_k$ for functions $f\in C^\infty(\R^3)$ 
by 
$$ N_k(f)= \sum_{m_1+m_2\leq 
k} \| U_1^{m_1} \cdots U_8^{m_8} f\|,$$
whenever $N_k(f)<\infty$. 
In particular,  we have 
$$S_k(F_\eps)=N_k(f_\eps)\leqno(9.2)$$
Let $B=[-2,2]^3$ and let $\tau\in C_c^\infty (\R^3)$ with $\tau\res_B
=1$. We claim that 
$$N_k(f_\eps)\leq N_k(\tau f_\eps) +C\leqno(9.3)$$
for all $\eps>0$ and a constant $C>0$ independent of $\eps$. 
Write 
$$f(x,y,z) ={1\over (1+x^2 +z^2)^{{1-\lambda_1}\over 2}
(1 + z^2)^{{1-\lambda_2}\over 2}} \cdot g(y)^{-{1\over 2}+ 
(\lambda_1-\lambda_2)}$$
and note that $f$ corresponds to the vector $v_0$ after 
the change of variables. Clearly we have 
$$N_k(f)=S_k(v_0)<\infty.$$
Now, outside the ball $B$ we can compare $f_\eps$ with $f$ and 
our claim (9.3) follows.

Now define the two variable function $g_\eps(x,z)$ by 
$$g_\eps(x,z)={1\over (e^{i(\pi-\eps)}+x^2 +z^2)^{{1-\lambda_1}\over 2}
(1 + e^{-i(\pi-\eps)} z^2)^{{1-\lambda_2}\over 2}}\ .$$
Let $\tau_2\in C_c(\R^2)$ be such that $\tau_2\res_{[-2,2]^2}=1$. 
Then it follows from (9.2) and (9.3) that there exists a 
constant $C>0$ such that 
$$ S_k((v_0)_\eps)\leq  C + C \sum_{j+l \leq k} \|{\partial^j
\over \partial x^j} {\partial^l
\over \partial z^l} \tau_2 g_\eps\|_{L^2(\R^2)}.$$
The same computation as in Lemma 9.5 then leads to 
$$S_k((v_0)_\eps)\asymp \eps^{-{1+2k\over 2}}.$$
Specializing to $k=5$ proves the proposition. \enddemo

\numbereddemo{{R}emark} (a) The estimate in Proposition 9.9 is certainly 
not optimal. The conjectured optimal estimate would be 
$$\|\theta_{(v_0)_\eps, \eta}\|_\infty\leq C\|(v_0)_\eps\|
\asymp  {1\over \sqrt\eps}.$$ The major technical difficulty
arises from the singularities of the spherical
vector $(v_0)_\eps$ for $\eps\to 0$  considered as a function on 
$\oline N$. In the rank one case the singular locus always is a {\it compact} 
variety, whereas in this case the singular locus lies on a complicated unbounded 
variety. 
This is reminiscent of the theory of intertwining operators in rank one versus 
higher rank. The technique there (essentially due to Gindikin-Karpelevic and 
later Schiffmann)
is the origin of the change of variables used in Proposition 9.4 but here it did not 
lead to a simple product structure.
 \enddemo

{\it Triple products}.
For a co-compact discrete 
subgroup $\Gamma<G$ recall that 
$Y=\Gamma\bs G$, $X=\Gamma\bs G/K$ 
and there is the Plancherel decomposition 
$$L^2(Y)\cong \hat\bigoplus_{\pi\in \hat G} m_\pi {\cal H}_\pi$$
and 
$$L^2(X)\cong\hat\bigoplus_{\pi\in \hat G_s} m_\pi {\cal H}_\pi^K$$
where $\hat G_s$ denotes the subset of $\hat G$ which corresponds 
to the unitary spherical representations. 

As before we let $(\psi_i)_{i\in I}$ be an orthonormal 
basis of Maa{\ss} forms of $L^2(X)$, also considered as an orthogonal
system in $L^2(Y)$. Note that $\psi_i(\Gamma g)=\theta_{v_0^i,\eta}=\la 
\pi_{\lambda_i}(g)v_0,\eta\ra$
for some unitary principal series representation $(\pi_{\lambda_i}, {\cal 
H}_{\lambda_i})$
and $\eta\in ({\cal H}_{\lambda_i}^{-\infty})^\Gamma$.

Let $X\in \oline{\b^1}$ be an extreme point; if $\psi=\theta_{v_0,\eta}$ is a Maa{\ss} form, then  
$\psi_\eps =\theta_{(v_0)_\eps,\eta}\in L^2(Y)$ for all $0<\eps\leq 1$. 
As before we have the 
identity 
$$\psi_\eps^2=\sum_{i\in I} c_i \psi_{i,\eps},$$
for all $0<\eps\leq 1$. Taking the norms, we have again 
$$\|\psi_\eps^2\|^2=\sum_{i\in I} |c_i|^2 \|\psi_{i,\eps}\|^2=\sum_{i\in I} 
|c_i|^2 \|(v_0^i)_\eps\|^2.
\leqno(9.4)$$

We now introduce new coordinates on $\a_\C^*$ by means of the simple roots;
i.e., 
we will write $\lambda=\lambda_1'\alpha_1+\lambda_2'\alpha_2$. 
As norm on $\a_\C^*$ we use the maximum norm $\|\lambda\| =\max\{|\lambda_1'|, 
|\lambda_2'|\}$.

\proclaim{Theorem} Let $G=\Sl(3,\R)$ 
and $\Gamma<G$ be a co-compact discrete subgroup. Then for every Maa{\ss} form 
$\psi$ 
corresponding to a
unitary principal series  there 
exists a constant $C>0$ such that for all
$T>1${\rm ,} for the coefficients 
$c_i$ of the Fourier series of $\psi^2=\sum_{i\in I} c_i \psi_i${\rm ,} one has
$$\sum_{\|\lambda_i\|\leq T} |c_i|^2 e^{\pi\|\lambda_i\|}\leq C  T^{12} 
.$$ \endproclaim

\demo{Proof} Given the previous estimates, the pattern of proof follows that of 
Theorem 7.6. and is left to the reader. \enddemo

\numbereddemo{{R}emark}  Friedberg, in [Fr87], determines precise gamma factors for the 
Rankin-Selberg convolution for $\Gamma=\Sl(3,\Z)$. By Stirling's 
approximation it is easily seen that the exponential growth of these gamma 
factors differs from the exponential term in Theorem 9.11. A likely explanation 
for this is that the classical Rankin-Selberg integral is computed using 
Eisenstein series off a maximal parabolic subgroup. Our results have used 
Eisenstein series or Maa{\ss} forms associated to the minimal parabolic. We think 
that any relationship of the exponential factors must involve the embedding 
parameters of the representation into the principal series off the minimal 
parabolic. We thank the referee for bringing this point to our attention.  \enddemo

{\it Generalizations to $\Sl(n,\R)$}.
Much of what has been said so far for $\Sl(3,\R)$ can be 
generalized easily to $G=\Sl(n,\R)$, 
$n\geq 3$.  
We will be interested in the radial limits of the 
spherical functions $\phi_\lambda$, $\lambda\in \a_\C^*$,  in the imaginary 
direction 
of the first fundamental weight 
$$X_{\omega_1}=\diag({n-1\over n}, -{1\over n}, \ldots, -{1\over n}).$$
First we need the splitting formula for the spherical functions analogous to 
Proposition  9.4. 
The splitting will be accomplished for the choice of subalgebras 
$\oline \n'=\bigoplus_{\alpha\in \Sigma^+\atop 
\alpha(X_{\omega_1})=0}\g^{-\alpha}$
and  $\oline \n''=\bigoplus_{\alpha\in \Sigma^+\atop 
\alpha(X_{\omega_1})>0}\g^{-\alpha}$. 
Then $\oline \n=\oline \n'\oplus\oline \n''$ and we note that $\oline \n''\cong 
\R^{n-1}$
is abelian. The splitting formula in Proposition~9.4 then generalizes to 

$$\phi_\lambda(\exp(i\pi(1-\eps)X_{\omega_1}))=C(\lambda, \eps) \int_{\R^{n-1}}
\prod_{j=1}^{n-1} {1\over f_j(\eps, x_j,\ldots,\  x_{n-1})} \ dx_1\ldots\ 
dx_{n-1},\leqno(9.5) $$
where 
$$f_j(\eps, x_j,\ldots,\  x_{n-1})=|e^{i(\pi-\eps)}+x_j^2+\ldots+x_{n-1}^2|
(e^{-i(\pi-\eps)}+x_j^2+\ldots+x_{n-1}^2)^{\Re \lambda_j}$$
for all $1\leq j\leq n-1$ and $C(\lambda, \eps)$ is bounded when $\eps\to 0$ for 
a fixed $\lambda$. With (9.5) one easily concludes as in Lemma  9.5 that 
$$\phi_\lambda(\exp(i\pi(1-\eps)X_{\omega_1}))\leq C{1\over 
\eps^{n-2}}\leqno(9.6)$$
for all $0<\eps\leq 1$ and a constant $C$ only depending on $\lambda$. 
Also, the proof of the lower estimate in Proposition  9.8 immediately 
generalizes to 
$$\phi_\lambda(\exp(i\pi(1-\eps)X_{\omega_1}))\geq Ce^{\pi (1-\eps)\sup_{w\in 
{\cal W}_\a}
w.\lambda(X_{\omega_1})}\leqno(9.7)$$
now for a constant $C$ independent of $\lambda$. 
Also, the proof of Proposition  9.9 generalizes easily to $G=\Sl(n,\R)$. For a 
co-compact 
subgroup $\Gamma<G$ we get that 

$$\|\theta_{(v_0)_\eps, \eta}\|_\infty^2\leq C 
\eps^{-(n+\big[{n^2-1\over 2}\big])}\leqno(9.8)$$
for all $0<\eps<1$. 
Clearly (9.5)--(9.8)  also hold for $\pm \pi w.X_{\omega_1}$, 
$w\in {\cal W}_\a$, 
instead of $\pi X_{\omega_1}$. 
We introduce now a norm on $\a_\C^*$ by setting $\|\lambda\|=\sup_{w\in {\cal 
W}_\a}|\lambda(w.X_{\omega_1})|$. 
Then the method in the proof of Theorem  9.11 together with (9.5)--(9.8) 
give the 
following estimates on triple products:

\proclaim{Theorem} Let $G=\Sl(n,\R)$ 
and $\Gamma<G$ be a co\/{\rm -}\/compact discrete subgroup. Then for every Maa{\ss} form $\psi$ 
corresponding to a
unitary principal series  the coefficients 
$c_i$ of the Fourier series of $\psi^2=\sum_{i\in I} c_i \psi_i$ there 
exists a constant $C>0$ such that for all
$T>1${\rm ,} one has
$$\sum_{\|\lambda_i\|\leq T} |c_i|^2 e^{\pi\|\lambda_i\|}\leq  C 
T^{2n +\big[ {n^2 -1\over 2}\big]-2}. $$ \endproclaim

\vglue16pt \centerline{\bf Appendix A: The case $G=\Sl(2,\R)$} 
\vglue12pt
In this appendix we deal with the group $G=\Sl(2,\R)$ where we can perform the 
most explicit computations. We think this is still of interest since it is the 
guiding 
example  on which one can make conjectures for the general 
case. In particular, some of the ideas of the proofs in Section 1 are   
present 
in the explicit computations  below.

For $G=\Sl(2,\R)$ note that $G\subeq G_\C=\Sl(2,\C)$ and $G_\C$ 
is simply connected. We let $\k=\so(2)$,  
$$\a =\{ \left(\begin{array}{cc} x & 0\\ 0 & -x\end{array}\right)\: x\in \R\}\quad \hbox{and}\quad \n =\{ 
\left( \begin{array}{cccc}0 & n\\ 0 & 0\end{array}\right)\: n\in \R\}.$$ 
For $z\in \C^*$, $x\in \C$ and  $\theta\in\C$  we set 
$$a_z =\left( \begin{array}{cccc}z & 0\\ 0 & z^{-1}\end{array}\right)\in A_\C, \quad 
n_x =\left( \begin{array}{cccc} 1& x\\ 0 & 1\end{array}\right)\in N_\C$$
and
$$
k_\theta =\left( \begin{array}{cccc}\cos \theta & \sin \theta\\ -\sin \theta & \cos\theta\end{array}\right)\in 
K_\C.$$
Note that 
$$A_\C^0=\{ a_z\: \Re(z)>0\}\quad\hbox{and}\quad A_\C^1=\{ 
a_z\:|\arg(z)|<{\pi\over 4}\}.$$
\enddemo

 \nonumproclaim{Proposition A.1} Let $G=\Sl(2,\R)$. Then the following assertions hold\/{\rm :}
\begin{itemize}
\ritem{(i)} For all $a_z\in A_\C^1$ 
and $\theta \in \R${\rm ,}
$$a_zk_\theta\in K_\C a_{z'} N_\C$$
with $a_z'\in A_\C^1$ and $z'$ defined by
$$z'=\sqrt{z^2+\sin^2\theta({1\over z^2} -z^2)}.$$
\ritem{(ii)} $A_\C^1 K\subeq K_\C A_\C^1 N_\C.$
\end{itemize}
\endproclaim

\demo{Proof} (i) Set $I =\{ z\in \C\: |\arg (z)|<{\pi\over 4}\}$ and fix $\theta\in 
\R$. 
Let $\Omega =\break \{ z\in I\: a_zk_\theta\in K_\C A_\C^1 N_\C\}$. By Lemma 1.4
the set $\Omega$ is open and not empty. We have to show that 
$\Omega=I$. For $z\in \Omega$ define 
analytic functions $z'(z)$, $\phi(z)$, $x(z)$ such that 
$$a_z k_\theta= k_{\phi(z)} a_{z'(z)} n_{x(z)}.$$
Writing this identity in matrix form yields
$$\left( \begin{array}{cccc}z\cos\theta  & z\sin\theta \\ -{\sin\theta\over z} & 
{\cos\theta\over z}\end{array}\right)=
\left( \begin{array}{cccc}z'\cos\phi & z'x\cos\phi   + {\sin\phi\over z'}\\ -z'\sin\phi &  -z'x\sin \phi  + {\cos\phi\over z'}\end{array}\right).$$
Thus we get 
$$z^2\cos^2\theta  +{\sin^2\theta\over z^2}= (z')^2\cos^2\phi  
+(z')^2\sin^2\phi$$
or equivalently 
$$(z')^2= z^2 +\sin^2 \theta ({1\over z^2} -z^2).\leqno({\rm A.1})$$
Taking real parts in (A.1) yields 
$$\Re(z')^2= ( 1-\sin ^2\theta) \Re (z^2) +\sin^2 \theta \Re ({1\over 
z^2})>0.\leqno({\rm A.2})$$ 
We conclude that 
$$z'\: I\to I, \ \ z\mapsto z'(z)\:=\sqrt{z^2+\sin^2\theta({1\over z^2} -z^2)}$$
is a well-defined holomorphic map.

Assume that $\Omega\neq I$.  Then there exists a sequence $(z_n)_{n\in \N}$ 
in
$\Omega$ such that $z =\lim z_n\in I\bs \Omega$. 

{}From (A.2) we now conclude that $z' =\lim_{n\to \infty} z'(z_n) $ exists in 
$I$. 
Now (A.1)
implies that  the limits $\phi =\lim_{n\to \infty}\phi(z_n)$ and 
$x =\lim_{n\to \infty} x(z_n)$ exist. Thus
$$a_zk_\theta=\lim_{n\to \infty} a_{z_n} k_\theta= \lim_{n\to \infty}
k_{\phi(z_n)} a_{z'(z_n)} n_{x(z_n)}=k_\phi a_{z'} n_x \in K_\C A_\C^1 N_\C, $$
contradicting $a_zk_\theta\not\in K_\C A_\C^1 N_\C$. This proves (i). 
\vglue4pt
(ii) This follows from (i).  \enddemo

\demo{{R}emark  {\rm A.2. (convexity theorems)}} For a linear semisimple Lie group $G$ 
with Iwasawa decomposition $G=KAN$ 
Kostant proved two convexity theorems (cf.\ [Kos73]): 
the {\it linear convexity theorem} which asserts
$$(\forall X\in\a) \qquad p_\a(\Ad(K).X)=\conv({\cal W}_\a.X)\leqno({\rm A.3})$$
with $p_\a\: \p\to\a$ the orthogonal 
projection with respect to the Cartan-Killing form, and 
the {\it nonlinear convexity theorem} which can be stated as 
$$(\forall a\in A) \qquad a(aK)=\exp\big(\conv({\cal 
W}_\a.\log(a))\big).\leqno({\rm A.4})$$

For $G=\Sl(2,\R)$ a simple calculation shows that 
(A.3) extends to   
$$(\forall X\in\a_\C) \qquad p_{\a_\C}(\Ad(K).X)=\conv({\cal W}_\a.X)$$
with $p_{\a_\C}\: \p_\C\to\a_\C$ the complex linear extension of $p_\a$
and we conjecture that (A.3) holds for all semisimple Lie groups 
$G$. 
Also it is natural to ask whether the nonlinear convexity theorem (A.4) 
generalizes to elements $a\in A_\C^1$ (this makes sense in view of 
$A_\C^1 K\subeq K_\C A_\C^{1,\leq} N_\C$ (cf.\ Proposition A.1(ii))). The answer 
is 
no and we can already see this for $G=\Sl(2,\R)$. Here we have ${\cal 
W}_\a\cong\Z_2=\{\1,s\}$
 with $s.a=a^{-1}$ for all $a\in A_\C$. For $z,w\in \C$ let 
$$l_{z,w} =\{\lambda z+(1-\lambda)w\: 0\leq \lambda\leq 1\}$$
be the line segment in $\C$ connecting $z$ and $w$.   Then Proposition A.1 shows 
that for all $z\in\C^*$ with $|\arg(z)|<{\pi\over 4}$ that 

$$a(a_zK)=\{ a_w\: w\in (l_{z^2, z^{-2}})^{1\over 2}\}.$$
Note that for $z=e^{i\phi}$, $0<|\phi|<{\pi\over 4}$ we have 
$\1\not\in\{ a_w\: w\in (l_{z^2, z^{-2}})^{1\over 2}\}$, but 
$\1\in\exp(\{ a_w\: w\in l_{-i\phi, i\phi}\})$. Thus we see that 
(A.4) usually does not hold for elements $a\in A_\C^1\bs A$.  \enddemo

\AuthorRefNames [AkGi99]

\end{document}